\documentclass{amsart}
\usepackage{melwun1}
\usepackage{amssymb,verbatim}
\usepackage{amscd}
\datver{Date: Wed Nov 21 08:57:39 2001
}
\begin{document}
\title[Wave equation on conic manifolds]
{Propagation of singularities for the wave equation on conic manifolds}
\author{Richard Melrose}
\address{Department of Mathematics, Massachusetts Institute of Technology}
\email{rbm@math.mit.edu}
\author{Jared Wunsch}
\address{Department of Mathematics, SUNY at Stony Brook}
\email{jwunsch@math.sunysb.edu}
\dedicatory{\datverp}
\thanks{First author supported in part by the National Science Foundation
under grant \#DMS-9622870, second author supported in part by the National
Science Foundation under grant \#DMS-0100501.}
\begin{abstract} For the wave equation associated to the Laplacian on a
  compact manifold with boundary with a conic metric (with respect to which
  the boundary is metrically a point) the propagation of singularities
  through the boundary is analyzed. Under appropriate regularity
  assumptions the diffracted, non-direct, wave produced by the boundary is
  shown to have Sobolev regularity greater than the incoming wave.
\end{abstract}
\maketitle

\intro
\section*{Introduction}

Solutions to the wave equation (for the Friedrichs extension of the
Laplacian) associated to a conic metric on a compact manifold with boundary
exhibit a diffractive, or ``ringing,'' effect when singularities strike the
boundary.  The main results of this paper describe the relationship between
the strength of the singularities incident on the boundary and the strength
of the diffracted singularities.  We first show that if \emph{no}
singularities arrive at the boundary at a time $\bar t$ then the solution
is smooth near the boundary at that time, in the sense that it is locally
in the intersection of the domains of all powers of the Laplacian. We then
show that if there are singularities incident on the boundary at time $\bar
t$ and in addition the solution satisfies an appropriate \emph{nonfocusing}
condition with respect to the boundary, then the strongest singularities
leaving the boundary at that time are on the geometric continuations of
those incoming bicharacteristics which carry singularities, whereas on the
diffracted, \ie\ not geometrically continued, rays the singularities are
weaker. If the incident wave satisfies a conormality condition, then the
singularity on the diffracted front is shown to be conormal. Applying this
analysis to the forward fundamental solution gives an extension of results
of Cheeger and Taylor (\cite{Cheeger-Taylor2, Cheeger-Taylor1}) from the
product-conic to the general conic case.  The results
contained in this paper represent a refinement of those previously
announced in \cite{Melrose-Wunsch1A}.

The problem of diffraction is an old one, with rigorous treatment
stretching back to the work of Sommerfeld, who, in 1896, discussed
diffraction around edges in the plane \cite{Sommerfeld1}; this includes the
case of Dirichlet boundary conditions at the slit $[0,\infty)\subset \CC.$
Using the method of images, this problem may be reduced to the study of the
wave equation on the cone over the circle of circumference $4\pi.$ An
overview of this and many other problems of diffraction around obstacles in
$\RR^n$ is given in Friedlander \cite{MR20:3703}.

Diffractive effects were extensively studied by Cheeger and Taylor
\cite{Cheeger-Taylor2, Cheeger-Taylor1} in the special case of
\emph{product-conic metrics}, \ie\ $\bbR_+$-homogeneous metrics on
$(0,\infty)\times Y$ of the form
$$
g= dx^2 + x^2 h(y, dy).
$$
Separation of variables can be used to give an explicit description of the
fundamental solution in terms of functions of the tangential Laplacian. See
also the discussion by Kalka and Menikoff \cite{Kalka-Menikoff1}. Rouleux
\cite{MR88g:58182} obtained a version of the results of Cheeger-Taylor
in the analytic category.  Lebeau \cite{Lebeau4, Lebeau5} has also obtained
a diffractive theorem in the setting of manifolds with corners in the
analytic category, and G\'erard-Lebeau \cite{MR93f:35130} have explicitly
analyzed the problem of an analytic conormal wave incident on an analytic
corner in $\RR^2.$ Estimates in $L^p$ spaces have been obtained for
product-conic metrics by M\"uller-Seeger \cite{MR96e:35093}.

Let $X$ be an $n$-dimensional conic manifold, that is to say, a compact
manifold with boundary, with a Riemannian metric $g$ on the interior which,
near the boundary, takes the degenerate form
\begin{equation}
g = dx^2 + x^2 h,
\label{intro:conicmetric}
\end{equation}
where $x$ is a boundary defining function and $h \in \CI(X; \Sym^2(T^*X))$
restricts to be a metric on $\pa X.$ Each boundary component of $X$ is thus
a ``cone point'' in the metric sense.  A trivial example of a conic metric
is obtained by blowing up a point $p$ in a Riemannian manifold; near $p$ we
can take $x$ to be the distance function to $p$, and
\eqref{intro:conicmetric} is simply the expression of the metric in
Riemannian polar coordinates. Such conic metrics exist on any manifold with
boundary.

The subject of this paper is the wave equation on a conic manifold, and in
particular the propagation of singularities for its solutions.  Let $\Lap$
be the Friedrichs extension of the (non-negative) Laplace-Beltrami operator
on $X.$ We consider solutions to
\begin{equation}
(D_t^2-\Lap) u=0
\label{intro:waveequation}
\end{equation}
on $\RR\times X^{\circ}$ which are \emph{admissible} in the sense that 
\begin{equation*}
u\in\mathcal{C}^{-\infty}(\bbR;\Dom(\Lap^{s/2}))\Mforsome s\in\bbR
\label{melwun1.316}\end{equation*}
and the equation holds in
$\mathcal{C}^{-\infty}(\bbR;\Dom(\Lap^{s/2-1})).$

H\"ormander's theorem \cite{MR48:9458} on the propagation of singularities
for operators of real principal type yields rather complete information
about the location of the singularities of any solution of
\eqref{intro:waveequation} away from $\pa X.$ Namely,
$\WF^r(u)$ (the wavefront set computed with respect
to the scale of Sobolev spaces) is contained in the characteristic variety
and is a union of maximally extended null bicharacteristics. Since the null
bicharacteristics are essentially time-parametrized geodesics over
$X^\circ$, this means, somewhat loosely speaking, that the singularities
travel with unit speed along geodesics.

H\"ormander's theorem does not address the question of what happens when
singularities reach $\pa X.$ That is, it leaves open the question of how
singularity-carrying geodesics in $\pa X$ \emph{terminating} at $\pa X$ are
connected with those \emph{emanating} from $\pa X$.  The answer involves
both the geometric propagation of singularities that one would expect from
the limiting behavior of geodesics that come close to the tip of the cone
without striking it, and an additional, diffractive, effect.

By the finiteness of the propagation speed for the wave equation it
suffices to consider a single component $Y$ of $\pa X$ (\ie\ a single
metric cone point). We first note that for every $y \in Y,$ there is a
unique (maximal, unparametrized) geodesic in $X^\circ=X\setminus\pa X$
with $y$ in its closure; conversely there is a neighborhood of $Y$ in which
each point is the end point of a unique, short, geodesic segment
terminating at $Y$.  This corresponds to a product decomposition $[0,
\ep)_x \times Y$ of a neighborhood of $Y$ in which projection to $Y$
is the map to the end-point and $x$ is the length of this segment. We
henceforth take $x$ to be the boundary defining function; the metric still
has the form \eqref{intro:conicmetric} but now with $h(\pa_x,\cdot)=0,$ so
$h$ is a family of metrics on $Y$ parameterized by $x.$

It follows that the surfaces $t\pm x=c,$ for any constant $c,$ which are
well-defined near each boundary component, are characteristic for the wave
equation. These radial surfaces carry the extra, diffractive, singularities,
the existence and regularity of which is our main object of study.

Our results are best illustrated by the fundamental solution
itself. Consider an initial point $\bar m\in X^{\circ},$ which we shall
take to be close to some boundary component $Y$. There is then a unique short
geodesic interval from $\bar m$ to the boundary; its length is $d(\bar
m,Y).$ The wave cone emanating from $\bar m$ is smooth for small
positive times but for longer times generally becomes singular. However it
is the projection of a smooth Lagrangian submanifold $I_{\bar m}\subset
T^*(\bbR\times X^\circ)$ consisting of the union of those maximally
extended null bicharacteristics of the wave operator which pass above $\bar
m$ at $t=0.$ In addition we consider the radial surface mentioned above,
$D_{\bar m}=\{t=x+d(\bar m,Y)\},$ which is well defined if $d(\bar
m,Y)$ is small enough and which emanates from the boundary at the time
of arrival of the geodesic from $\bar m.$

\begin{theorem}[Fundamental solution]\label{melwun1.374} Let $E _{\bar
m}$ be the fundamental solution of the wave equation for the
Friedrichs extension of the Laplacian of a conic metric on a compact
manifold with boundary $X,$ with pole at $\bar m\in X^\circ.$ If $d(\bar
m,Y)$ is sufficiently small then, regarding $E_{\bar m}$ as a
distribution on $\bbR\times X^\circ,$
\begin{equation}
\WF(E_{\bar m})\subset I_{\bar m}\cup N^*D_{\bar m}\Min |t|<2d(\bar m,Y)  
\label{melwun1.375}\end{equation}
and $E_{\bar m}$ is conormal with respect to $D_{\bar m}$ away
from the projection of the closure of $I_{\bar m}$ and is of Sobolev order
$\frac12-\delta$ there for any $\delta >0.$
\end{theorem}
A comparison with the explicit results of Cheeger-Taylor
\cite{Cheeger-Taylor2, Cheeger-Taylor1} in the special case of product
cones shows that this result is optimal as far as Sobolev regularity is
concerned.

The first step in proving this is to show that no singularities arise from
the boundary spontaneously.

\begin{theorem}[Diffractive Theorem]\label{intro:diffthm} 
If an admissible solution to the wave equation on a conic manifold has
singularities, of Sobolev order $r,$ on at least one null
bicharacteristic hitting a boundary component at time $\bar t$ then it
has singularities, of Sobolev order $r,$ on at least one null
bicharacteristic leaving that boundary component at time $\bar t.$
\end{theorem}
\noindent
This diffractive theorem does not give any localization of singularities in
the boundary component $Y.$ It in no way distinguishes between the
different ``outgoing'' null bicharacteristics leaving the boundary
component at $t=\bar t$ although in fact it is easily strengthened to show
that singularities in the two components of the characteristic set,
corresponding to the sign of the dual variable to $t,$ do not interact at
all.

In order to state a more refined theorem, which distinguishes between
different points in $Y,$ we need first to consider the possible \emph{geometric
continuations} of a geodesic terminating at $y \in Y.$ Geometric
continuation of geodesics corresponds to the relation 
\begin{equation}
G(y)=\{y'\in Y; y\Mand y'\text{ are endpoints of a geodesic segment of length }
\pi\Min Y\},
\label{melwun1.329}\end{equation}
where $Y$ is endowed with the metric $h\restrictedto_{x=0}.$  The geodesic
segments of length $\pi$ in $Y$ arise naturally as limits of geodesics in
$X^\circ$ which narrowly miss $Y.$

Our second main result shows that under certain circumstances, a
singularity arriving at a point $y \in Y$ at $t=\bar t$ will produce only
\emph{weaker} singularities along rays emanating from points in $Y\setminus
G(y)$ than on the geometrically continued rays, which are those emanating
from $G(y).$ In other words, we show that diffracted singularities are
weaker than geometrically propagated ones.

The further hypothesis that we make is one of nonfocusing of
singularities. Suppose that the solution $u$ is microlocally of Sobolev
order $r$ on all incoming null bicharacteristics arriving at $Y$ at $t=\bar t.$ We
require that \emph{tangential smoothing} improves this regularity; for our
purposes the amount of tangential smoothing required to gain regularity is
irrelevant. Let $\Lap_Y$ denote the Laplace-Beltrami operator on $Y$ with
respect to the metric $h\restrictedto_{x=0}$, extended to operate on a
collar neighborhood of $Y$ by the metric product decomposition described
above. We assume that for $N$ sufficiently large $(1+\Lap_Y)^{-N} u$ is
more regular (microlocally) than $u$ in the sense that
\begin{multline}
p\notin\WF^{r'}\left((1+\Lap_Y)^{-N} u\right)\Mfor p\Mon\\
\text{ the incoming null bicharacteristic reaching }Y\text{ at }y\in Y,\
t=\bar t.
\label{intro:nonfocusing}
\end{multline}

\begin{theorem}[Geometric Propagation]\label{melwun1.318} If $u$ is an
admissible solution of the wave equation, then the outgoing null
bicharacteristic emanating from a point $y\in Y$ at $t=\bar t$ is
disjoint from $\WF^{R}(u)$ provided all the incoming null bicharacteristics
reaching $G(y)$ at $t=\bar t$ are outside $\WF^{R}(u)$ and 
\eqref{intro:nonfocusing} holds with $r'>R$ microlocally for all other incoming
null bicharacteristics meeting the boundary at time $\bar t.$
\end{theorem}

\noindent This theorem may not be strengthened by dropping the nonfocusing
assumption \eqref{intro:nonfocusing} as is shown by a counterexample in
\S\ref{section:conormal}. On the other hand the two components of the
characteristic variety are again completely independent.

At a time $\bar t$ the \emph{background regularity,} $r,$ for an admissible
solution near a boundary component is the Sobolev regularity that holds on
all incoming null bicharacteristics arriving at that boundary component at
$t=\bar t.$ Theorem~\ref{intro:diffthm} shows that this background
regularity propagates to all outgoing null bicharacteristics leaving the
boundary at time $\bar t.$ Theorem~\ref{melwun1.318} above shows
that there is additionally a gain of any number, $l,$ of Sobolev
derivatives on any one of these outgoing bicharacteristics, say emanating
from $\bar y\in Y,$ provided that $u$ has Sobolev regularity $r+l$ on all
incoming bicharacteristics which are geometrically related to this one and
that $u$ has the property that sufficient tangential smoothing increases
its regularity microlocally near all other (and hence all) incoming null
bicharacteristics, reaching the boundary at time $\bar t,$ to greater than
$r+l.$

Although there is in principle no upper limit to the gain of regularity
compared to background available from Theorem~\ref{melwun1.318}, there are
practical limits. The most important case of nontrivial tangential
smoothing is when the solution is, for $t<\bar t$ and away from the
boundary, a conormal distribution associated to a hypersurface that is
simply tangent to the incoming radial surface $x=\bar t-t.$ For such a
distribution, the lemma of stationary phase shows that tangential smoothing
gains almost $(n-1)/2$ derivatives, hence Theorem~\ref{melwun1.318}
guarantees that the diffracted wave is almost $(n-1)/2$ derivatives
smoother than the main singularity.  In addition, the diffracted wave for
such a solution will itself be conormal.

More precise statements of the theorems above as well as further
results are given in Section~\ref{section:statement}.  We now briefly
describe the ingredients in the proofs.

We make extensive use of the calculus of edge pseudodifferential operators
as developed by Mazzeo \cite{Mazzeo4}. A related calculus of
pseudodifferential operators adapted to edge structures has also been
constructed by Schulze \cite{Schulze2} but the lack of ``completeness'' in
this calculus makes it much less applicable; in particular it does not seem
to have associated with it a useful notion of wavefront set. We develop and
use just such a notion for Mazzeo's calculus measuring microlocal
regularity with respect to the intrinsic weighted Sobolev spaces. We use
this to obtain a result on the propagation of singularities at $\pa X$
analogous to H\"ormander's interior result.  The technique of the proof is
a positive commutator construction similar to that used in
\cite{MR48:9458}, but with the crucial distinction, familiar from
scattering theory, that the bicharacteristic flow now has radial points,
thus necessitating a more subtle construction.  At the radial surfaces the
propagation arguments are only valid for a limited range of weighted (edge)
Sobolev spaces, corresponding to a linear estimate on the weight in terms
of the regularity. Such an estimate amounts to a divisibility property for
the solution in terms of powers of the boundary defining function.

The proof of the diffractive theorem, Theorem~\ref{intro:diffthm}, relies on
the extraction of a leading part, the non-commutative normal operator in
the edge calculus, via an appropriately rescaled FBI
(Fourier-Bros-Iagolnitzer) transform $S$ near $Y$.  The model operator is
$\Lap_0-1$, where $\Lap_0$ is the Laplacian on the tangent cone
$\RR_+\times Y$ for the product-conic metric $dx^2+ x^2(h\restrictedto_{x=0}).$
Iterative application of the outgoing resolvent for
this model operator, corresponding to the scattering structure at the
infinite volume end of the cone, combined with the microlocal propagation results
discussed above, yields the regularity of $u.$ We use estimates from
\cite{Melrose43}, although in this product case we could instead rely on
direct methods, as used for instance in \cite{Cheeger-Taylor2}.

To prove the geometric propagation theorem, we begin by establishing a
\emph{division theorem} to the effect that the nonfocusing assumption
\eqref{intro:nonfocusing} implies that the solution has better decay in $x$
than would be predicted merely by energy conservation. This allows us to
apply the propagation results to prove a special case of
Theorem~\ref{melwun1.318}.

We also show that conormality on a radial surfaces persists for solutions.
That is, if $(u(0), D_t u(0))$ are conormal distributions at a surface
$x=\bar x$, supported away from a boundary component $Y$, and if $\bar x$
is sufficiently small then by standard interior regularity results $(u(t),
D_t u(t)),$ is conormal with respect to $\{x=\bar x-t\}\cup \{x=\bar x
+t\}$ for small positive $t<\bar x;$ we show that it continues to be
conormal with respect to $\{x=\abs{\bar x - t}\}\cup \{x=\bar x+t\}$ for
small $t>\bar x,$ \ie\ after the wave leaves the boundary. Such conormal solutions
provide the counterexamples mentioned following Theorem~\ref{melwun1.318}.

For solutions, such as the fundamental solution, that are initially
conormal at a surface meeting the radial surfaces $x=\bar x-t$
with at most simple tangency, we are additionally able to prove that the
diffracted front is conormal; this is in effect a microlocal version of the
radial conormality argument discussed above. The conormality of the
diffracted front and the special case of Theorem~\ref{melwun1.318}
mentioned above suffice to establish conormality and (sharp) regularity for
the diffracted wave of the fundamental solution, which is then used to
prove Theorem~\ref{melwun1.318} in its full generality.

The outline of the paper is as follows.  
In \S\ref{Normal.form} we prove the existence of a product decomposition
for a conic metric in a neighborhood of $\pa X$; this is equivalent to
reducing the metric to the normal form
$$
g=dx^2+ x^2 h(x,y,dy).
$$
Then in \S\ref{section:mapping} and \S\ref{section:domains}, we discuss the
mapping properties and domains of the Laplacian and its powers on a conic
manifold.  This enables us to give precise statements of the main theorems
in \S\ref{section:statement}.
Sections~\ref{section:calculus}--\ref{section:propagation} build up the
machinery of the edge calculus and culminate in the proof of the
propagation theorem for the edge wavefront set.  Then in
\S\ref{section:FBI} and \S\ref{section:reduced} we discuss the rescaled FBI
transform and the normal operator in the edge calculus, which we use to
prove Theorem~\ref{intro:diffthm} in \S\ref{section:diffractive}.  In
\S\ref{section:tangential} we demonstrate the conservation of tangential
regularity, \ie\ iterated regularity under powers of $\Lap_Y$.  This is
then used in two ways.  In \S\ref{section:radial} we use conservation of
tangential regularity together with conservation of iterated regularity
under vector fields of the form $(x D_x + t D_t)$ to prove conservation of
radial conormality.  In \S\ref{section:divthm1}, we use the conservation of
tangential regularity to prove the division theorem; the proof of the
sharper division result for conormal initial data is closely related to
that of radial conormality.  In \S\ref{section:genprop} we prove a
preliminary version of Theorem~\ref{melwun1.318} and establish a theorem on
conormality of the diffracted front.  In \S\ref{section:conormal} we
discuss a class of examples including the fundamental solution and prove
Theorem~\ref{melwun1.374}; as a consequence we then prove
Theorem~\ref{melwun1.318} in full generality, and discuss a counterexample
to that theorem when the nonfocusing condition is omitted.

The authors thank Daniel Grieser and Andr\'as Vasy for comments on the manuscript.

\paperbody
\section{Conic metrics; geodesics and normal form\label{Normal.form}}

Let $X$ be an $n$-dimensional manifold with compact boundary, $\pa X$ its
boundary and $X^\circ=X\setminus\pa X$ its interior.
\begin{definition}
  A \emph{conic metric} on $X$ is a Riemannian metric $g$ on $X^\circ$ such
  that in a neighborhood of any boundary component $Y$ of $X$, there exists
  a boundary defining function $x$ ($x\geq 0$, $\{x=0\}=\pa X$,
  $dx\restrictedto_{\pa X} \neq 0$) in terms of which
\begin{equation}
  g= dx^2 + x^2 h
\label{melwun1.148}\end{equation}
  where $h \in \CI (\Sym^2 T^* X)$ and $h\restrictedto_Y$ is a metric.  A
\emph{conic manifold} is a compact manifold with boundary endowed with a
conic metric. 
\end{definition}

Let $Y$ be a compact manifold without boundary.  The \emph{product-conic}
metrics on $[0,\infty)_x \times Y$ are the conic metrics of the form
\begin{equation}
g=dx^2+x^2h_0
\label{melwun1.149}\end{equation}
where $h_0$ is a metric on $Y$.  These are the model cases which form the
basis of our analysis below. They also motivate the basic normal form for
conic metrics which we discuss next.

\begin{theorem}\label{thm:normalform}
  Let $g$ be a conic metric on $X.$ There exists a collar neighborhood
  $\ocal$ of $\pa X$ and an isomorphism $(x,\Pi):\ocal\longrightarrow
  [0,\epsilon) \times \pa X,$ $\epsilon >0,$ such that in terms of this product
  decomposition
\begin{equation}
  g= dx^2+x^2 \Pi^*h(x)\Mnear \pa X,
\label{melwun1.150}\end{equation}
where $h(x)\in \CI([0, \ep); \text{metrics on }\pa X)$.
\end{theorem}

\begin{proof} The existence of such a normal form distinguishes a vector
  field $V$ by 
\begin{equation}
Vx=1\Mand \iota_{V} h(x) =0\Mnear \pa X.  \label{melwun1.151}\end{equation}
The integral curves of this vector field are geodesics. Our main task is
thus to show that through any point in a neighborhood of $\pa X$, there is a
unique short geodesic reaching $\pa X;$ the length of this short geodesic
segment will then furnish the desired defining function for $\pa X$.

As our construction is local near a boundary component, we assume without
loss of generality that $\pa X$ is connected.  We begin with any collar
neighborhood and associated projection $\pi_{\pa X}:\ocal \to \pa X$, and
choose a boundary defining function $\rho$ for $\pa X$ and coordinates
$\upsilon=\pi_{\pa X}^* \bar y$ with $\bar y$ any coordinates on $\pa X$.
We may assume $\rho$ to have been chosen in accordance with our definition
of a conic metric, so that
$$
g=d\rho^2+\rho^2 h
$$
where $h$ is a smooth, symmetric two-tensor.  Let
$$
h=\bullet d\rho^2+
\bullet d\rho\, d\upsilon + k_{ij}(\upsilon) d\upsilon^i d\upsilon^j +
\rho\bullet d\upsilon^i d\upsilon^j;
$$
then
\begin{equation}
\ang{\cdot,\cdot}_g^*=
    (1+\ocal(\rho^2))\pa_\rho^2 + \ocal(1)\pa_\rho\pa_\upsilon+ (\rho^{-2}
    k^{ij}(\upsilon)+ \ocal(\rho^{-1})) \pa_{\upsilon_i} \pa_{\upsilon_j},
\label{new.dualmetric}
\end{equation}
where, here and henceforth, $\ocal (\rho^k)$ means $\rho^k$ times a smooth
function of $(\rho,\upsilon).$

We now regard the dual metric as a function on a rescaled version of the
cotangent bundle.  Let $\bT X$ denote the \emph{b-tangent bundle} whose sections
are the vector fields tangent to $\pa X$.  Let $\bTstar X$ be its dual.  Then
sections of $\bTstar X$ are $\CI$-linear combinations of $d\rho/\rho$ and
$d\upsilon_i$'s.  Writing the canonical one-form on $\bTstar X$ as 
\begin{equation}
\xi \frac{d\rho}\rho + \eta\cdot d\upsilon,
\label{bform}
\end{equation}
we can now write the dual metric as
\begin{equation}
\ang{(\xi,\eta),(\xi,\eta)}_g^* =  \frac{\xi^2+ k^{ij} \eta_i \eta_j}{\rho^2}
 + \ocal(1) \xi^2 + \ocal(\rho^{-1}) \xi \eta + \ocal(\rho^{-1}) \eta^2.
\label{dual:metric}
\end{equation}

Over the interior $\bTstar X^\circ \equiv T^* X^\circ$ and the canonical
symplectic form on $T^*X$ lifts to a form on $\Tbstar
X,$ singular at the boundary, given in canonical coordinates by $d (\xi
d\rho/\rho + \eta\cdot d\upsilon).$  Associated with this symplectic form 
and the energy function $\ang{\cdot,\cdot}_g^*$ is the Hamilton vector field
\begin{align*}
H_g  &= \frac 2{\rho^2} (\xi \rho \pa_\rho + (\xi^2+k^{ij} \eta_i\eta_j)
\pa_\xi + H_{\pa X})+ P, \\ \intertext{where}
    P &= \frac{1}{\rho^2} ((\ocal(\rho^3)\xi+\ocal(\rho^2)\eta)\pa_\rho
    + (\ocal(\rho^3) \xi^2+ \ocal(\rho)\xi\eta+\ocal(\rho)\eta^2)\pa_\xi \\
    &+  (\ocal(\rho) \xi+ \ocal(\rho)\eta) \pa_\upsilon
    + (\ocal(\rho^2)\xi^2+\ocal(\rho)\xi\eta+\ocal(\rho)\eta^2)\pa_\eta,
\end{align*}
and $H_{\pa X}$ is geodesic spray in the $(\upsilon,\eta)$ variables with
respect to the metric $k.$ The projections of the integral curves of $H_g$
to $X$ are geodesics.

Let $H_0=H_g-P$ be the main term in the Hamilton vector field.  The vector
field $H_0$ is thus tangent to the submanifold $N=\{ \eta=0 \}$ of ``normal''
directions to ${\pa X}$.  Since on $N$, the flowout of $(\rho^2/2) H_0$ is
$\rho=\rho_0/(1-\xi_0 s)$ and $\xi=\xi_0/(1-\xi_0 s)$, the map
$$
p \mapsto \lim_{s\to+\infty} \exp \left[-s \sgn\xi \frac{\rho^2}{2} H_0\right]
$$
maps $N$ to ${\pa X}$ smoothly in a neighborhood of the boundary, taking $p
\mapsto \upsilon(p) \in {\pa X}$.  The projections of integral curves of
$H_0$ on $N$ would be geodesics reaching the cone point were the metric of
the desired form \eqref{melwun1.150}.

We now show that an analogous \emph{normal manifold,} given by a perturbation of
$N$, exists for the full vector field $H_g$.  The flow along this manifold
to the boundary will provide the desired geodesics to the cone point.

To simplify the discussion of the vector field, $H_g$, we scale away its
homogeneity in $(\xi, \eta)$ as follows.  Consider the smooth function
$$
\sigma=\frac 1{\rho(\ang{\cdot,\cdot}_g^*)^{1/2}}
$$
on the complement of the zero-section in $\bTstar X$.  Note that $\sigma$
is approximately equal to $(\xi^2+k^{ij} \eta_i \eta_j)^{-1/2}$, and is
homogeneous of degree $-1$ in the fibers of $\bTstar X$.  We further set
\begin{equation}
\bar\xi=\sigma \xi,\ \bar \eta = \sigma \eta.
\label{vf:1}
\end{equation}
Since $\ang{\cdot,\cdot}_g^*$ is preserved under the flow
$H_g$, $(x^2/2) H_g \sigma = -\bar\xi + \ocal(\rho^2) \bar\xi +
\ocal(\rho) \bar \eta$.  Note that the error terms in this expression
are at least quadratic in $(\rho,\bar\eta).$
Thus we can write
\begin{equation}
\h \rho^2 \sigma H_g = \bar\xi \sigma  \pa_\sigma + \bar\xi \rho \pa_\rho +k^{ij}
\bar\eta_j\pa_{\upsilon_i} + k^{ij} \bar\eta_i \bar\eta_j \pa_{\bar\xi} -
\left(\bar\xi_k +\h \frac{\pa k^{ij}}{\pa \upsilon_k} \bar\eta_i \bar\eta_j\right)
\pa_{\bar\eta_k} + P'
\label{vf:2}
\end{equation}
where the perturbation term has the form
$$
P' = \ocal(\rho) \bar\xi \pa_\upsilon + \ocal(\rho^2+\bar\eta^2)\left(\CI\text
{ vector field}\right);
$$
the single non-quadratic error term comes directly from the corresponding
term in $P$.  By homogeneity, the error terms above are independent of
$\sigma$, \ie\ $(1/2) \rho^2 \sigma H_g$ pushes forward to a vector field
on the unit sphere bundle $\bSstar X\equiv \bTstar X/\RR_+$.  On $\bSstar
X$, $(\rho,\upsilon,\bar\eta)$ are coordinates near $\bar\eta=0$, since
$\ang{(\bar\xi,\bar\eta),(\bar\xi,\bar\eta)}_g^*=1$ at $\bar\eta=0$,
$\bar\xi=\pm 1$.  Moreover we easily see that $\bar\xi\mp 1$ vanishes
to second order at $\rho=\bar\eta=0$; of crucial importance is the fact
that the power of $\rho$ multiplying the $\xi^2$ error term in
\eqref{dual:metric} is larger than those multiplying the $\xi \eta$ and
$\eta^2$ terms.

In these coordinates, the linearization of the vector field $(1/2) \rho^2
\sigma H_g$ in $(\rho, \bar \eta)$ near $\rho=\bar\eta=0$ and $\bar\xi=\pm
1$ is simply
$$
L=\pm \rho\pa_\rho+ (k^{ij}\bar\eta_j+\ocal(\rho))\pa_{\upsilon_i}\mp
\bar\eta_i \pa_{\bar\eta_i}
$$
Since $(1/2) \rho^2 \sigma H_g$ vanishes identically at $\rho=\bar\eta=0$,
and its linearization has eigenvalues $\pm 1$ in the normal directions,
$(1/2)\rho^2 \sigma H_g$ is $r$-normally hyperbolic near $\rho=\bar\eta=0$
for all $r$, in the notation of \cite{Hirsch-Pugh-Shub1}.  Hence by the
Stable/Unstable Manifold Theorem as stated in Theorem~4.1 of
\cite{Hirsch-Pugh-Shub1}, near $\{\rho=\bar\eta=0,\ \bar \xi= 1\}$ there
exists a stable invariant manifold $N$ for the flow of
$(1/2) \rho^2 \sigma H_g$, with $N$ tangent to $\{\bar\eta=0\}$ at
$\{\rho=\bar\eta=0\}$.  Because $N$ is tangent to
$\{\bar\eta=0\}$, the projection map
$$
\pi: \bSstar X \to X
$$
restricts to give a diffeomorphism $N \cong X$ (locally, for $\rho$
sufficiently small).

On $N$, $\lim_{s\to-\infty} \rho=0$, where $s$ parametrizes the flow along
$(1/2)\rho^2\sigma H_g$; using $\rho$ as a parameter along the flow on
$N$ yields, by \eqref{vf:1} and \eqref{vf:2},
$$
\frac{d\upsilon}{d\rho}=\ocal(1) + \ocal(\bar\eta/\rho);
$$
the latter term is in fact $\ocal(1)$ since $N$ is tangent to
$\{\bar\eta=0\}$ at $\rho=0$.  Thus, $\lim_{s\to-\infty} \upsilon$ is in
fact a smooth map.  Hence
\begin{equation}
\Pi: p \mapsto \lim_{s\to -\infty} \exp_{\pi^{-1}(p)} \h \rho^2\sigma H_g
\label{Pi}
\end{equation}
is a smooth map from a neighborhood $\ocal'$ of ${\pa X}$ in $X$ to ${\pa
X}$.  Indeed $\Pi$ is a fibration, with fibers $\RR_+$ given by the
projections of integral curves of $(1/2) \rho^2 \sigma H_g$, \ie\ by
geodesics hitting the ``cone point'' $\pa X$.

On $\ocal'$, set $y_i(p)=\upsilon_i(\Pi(p))$ and $x(p) = d_g(p,\Pi(p))$.
To first order at $\pa X$, $x=\rho$ and $y=\upsilon$.  Hence
$(x,y_1,\dots,y_{n-1})$ form a coordinate system on a neighborhood
$\ocal''$.  In these coordinates, $\ang{\pa_x,\pa_x}_g=1$ since $\exp
t\pa_x$ is unit speed flow along geodesics reaching ${\pa X}$.
Furthermore, $\ang{\pa_x,\pa_{y_i}}_g=0$ for all $i$ by Gauss's
Lemma\footnote{The usual proof of Gauss's Lemma using the first variation
formula works even at a cone point, \ie\ geodesics with one endpoint on
${\pa X}$ are orthogonal to the hypersurfaces $d_g(p,{\pa X})=\epsilon$.
We simply let $y(u)$ be a curve in ${\pa X}$, and apply the first variation
formula to the family of geodesics connecting $x=x_0$, $y=y(u)$ to $x=0$,
$y=y(u)$.}  and $\ang{\pa_{y_i},\pa_{y_j}}_g = \ocal (x^2)$.  Hence in the
coordinates $(x,\ y)$, $g$ takes the form \eqref{melwun1.150}.
\end{proof}

Since $x$ is the distance along the normal geodesics it is uniquely
determined by \eqref{melwun1.150}; the vector field $V$ determining the
product decomposition is also fixed geometrically by \eqref{melwun1.151}.
The choice of a conic metric additionally induces a metric on ${\pa X},$
namely
$$
h_0=h(0).
$$
\emph{Henceforth, $x$ will always denote this distance function for the
given conic metric.}

The proof of Theorem~\ref{thm:normalform} used in a crucial way the
existence of a unique normal geodesic starting at each point of the
boundary. These geodesics foliate a neighborhood of the boundary, and
indeed, the existence of this foliation characterizes conic metrics among
more general nondegenerate forms in $dx$ and $x dy_i$.  On $\RR\times
X^\circ$, we thus obtain a foliation of a neighborhood of $\RR \times \pa
X$ by projections of null bicharacteristics for the symbol of the
d'Alembertian, $D_t^2-\Lap$.

We shall fix notation for various sets corresponding to the geodesic
segments that hit the boundary and are within small distance $\epsilon>0$
of it.
\begin{definition}
Let $N^\ep=\{\alpha dt + \beta dx: \abs\alpha=\abs\beta,\ x<\ep \} \subset
T^*(\RR\times X^\circ)$.  For $\bar t \in \RR$, $\bar y \in Y$, let
\begin{align*}
R_{\pm, I}^\ep(\bar t, \bar y) &= N^\ep \cap \{\sgn \alpha\beta =1,\
\sgn \alpha = \pm,\ t=\bar t-x,\ y=\bar y\},\\
R_{\pm, O}^\ep(\bar t, \bar y) &= N^\ep \cap \{\sgn \alpha\beta =1,\
\sgn \alpha = \pm,\ t=\bar t+x,\ y=\bar y\}.
\end{align*}
If any of the parameters $\bullet = I\text{ or } O$, $\pm$, $\bar t$, or
$\bar y$ is omitted, the resulting set is defined as the union over all
possible values of that parameter.  If $Y$ is a single component of $\pa
X$, let
$$
R_{\pm, \bullet}^\ep(\bar t, Y) = \bigcup_{y \in Y} R_{\pm, \bullet}^\ep(\bar t, y).
$$
\end{definition}
\noindent The set $N^\ep$ is the \emph{normal set} near the boundary of $X$ of
points in the cotangent bundle which lie along geodesics (projections of
bicharacteristics) entering and leaving $\pa X$.  Here ``I'' and ``O''
stand for the ``incoming'' and ``outgoing'' components of $N$, on which
$dx/dt$ is respectively negative and positive.  The additional sign $\pm$
microlocalizes in the sign of the dual variable to $t$.  The sets $R_{\pm,
\bullet}^\ep (\bar t, \bar y)$ are the points through which the short
geodesic to the boundary arrives at or departs from $\pa X$ at time $t=\bar
t$ and at the point $\bar y.$

If $X$ is a compact conic manifold then every geodesic starting at an
interior point can be extended maximally in both directions until and
unless it terminates at the boundary. This naturally suggests the question of
how such geodesics can, or should, be further extended.

\begin{definition}\label{melwun1.152} By a \emph{limiting geodesic} in a
  conic manifold we mean a continuous piecewise smooth curve
  $c:I\longrightarrow X,$ where $I=\cup_jI_j$ is decomposed as a locally
  finite union of relatively closed subintervals on each of which $c$
  restricts to be a smooth curve $c_j$ and such that
\begin{enumerate}
\item Each $c_j$ is either a geodesic in $X$ or a geodesic (for $h_0)$ in
$\pa X$ and such segments alternate.
\item Boundary segments are of length at most $\pi$ and if such a boundary
segment is not the first or last segment then its length is exactly $\pi.$
\end{enumerate}
\end{definition}

\begin{lemma}\label{melwun1.153} If $\xi _i$ is a sequence of geodesics in
$X^\circ$ which converges uniformly as curves in $X$ then its limit is a
limiting geodesic and conversely an open neighborhood of each boundary
segment of any limiting geodesic arises as such a limit.
\end{lemma}
\noindent This result will not be used except as motivation, hence we omit
its proof.

In view of this behavior of the geodesics we define a singular
relation, interpreted as a set-valued map
$$
\Gamma^\ep \subset R_{O}^\ep \times R_{I}^\ep\\
$$
such that
\begin{equation}
\text{for } p \in R_{\pm, O}^\ep (\bar t),\ \Gamma^\ep (p) = \{ q\in
R_{\pm, I}^\ep(\bar t) : \Pi(q) \in G(\Pi(p)) \}, \label{15.5.2000.9}
\end{equation}
where $\Pi$ is the projection map from a collar neighborhood of $\pa X$ to
$\pa X$ given by \eqref{Pi}, and for $y \in \pa X$, $G(y)$ is the geodesic
relation defined in \eqref{melwun1.329}. Clearly,
\begin{equation}
R_{\pm,I}^\epsilon (\bar t)=\bigcup_{p \in R_{\pm, O}(\bar t)}\Gamma^\epsilon
(p).
\label{15.5.2000.10}\end{equation}
In fact $\Gamma^\epsilon (p)$ is generically of codimension one
in $R_{\pm,I}^\epsilon (\bar t).$

To the given conic metric and a boundary component $Y$ we associated the
limiting product metric
\begin{equation}
g_0= dx^2 + x^2 h_0 \label{melwun1.163}
\end{equation}
on the normal bundle to $Y$, which we may identify with $[0,\infty)\times
Y$ using the decomposition associated with \eqref{melwun1.150}. Let $\Lap$
denote the (nonnegative) Laplace-Beltrami operator with respect to the
metric $g$ and $\Lap_0$ that with respect to $g_0$. Then, letting
$(w_0,\dots,w_{n-1}) = (x,y_1,\dots,y_{n-1})$,
$$
\Lap=\sum_{j,k=0}^n \frac{1}{\sqrt{g}} D_{w_j} g^{jk} \sqrt{g} D_{w_k}
$$
becomes, near $Y$,
\begin{multline}
\Lap= D_x^2 -\frac{i[(n-1)+xe]}{x} D_x -\frac i2 
D_x+ \frac{\Lap_{h}}{x^2},\\ \text{where } e=\frac12 \frac{\pa\log\det h(x)}{\pa x}=
\frac12\tr\left(h^{-1}(x)\frac{\pa h(x)}{\pa x}\right).
\label{melwun1.36}\end{multline}
Here $\Lap_{h}$ is the Laplacian on ${\pa X}$ with respect to the
($x$-dependent) metric $h(x)$ on ${\pa X}.$  Similarly
\begin{equation}
\Lap_0 = D_x^2 -\frac{i(n-1)}{x} D_x + \frac{\Lap_{h_0}}{x^2},
\label{modellap}
\end{equation}
hence if we identify $[0,\infty)\times Y$ with the metric product
decomposition near $Y,$
\begin{equation}
\Lap-\Lap_0 \in x^{-1} \left( \CI(x,y) (x D_x) + \sum \CI(x,y) D_{y_i}
  D_{y_j} + \sum \CI(x,y) D_{y_i} \right)
\label{lap:perturbation}
\end{equation}
We will also use the notation $\Lap_{Y}$ for $\Lap_{h_0}$, when restricting
our attention to a single boundary component.

If $g$ is a product-conic metric, it is easy to check that
\begin{align*}
[\Box, \Lap_Y]&=0,\\ [\Box, (x D_x + (t-\bar t) D_t)] &= -2i \Box.
\end{align*}
In the general conic case, these ``symmetries'' are broken.  It is crucial
for our purposes, though, that perturbed versions of the above identities
still hold.  In particular, if for brevity we set
\begin{equation}
R=x D_x + (t-\bar t) D_t,
\label{melwun1.175}\end{equation}
then
\begin{multline}
[\Box, \Lap_Y] =Q D_x +  x^{-1} P,\text{ with}\\ Q \in \CI([0,\ep); \Diff
1(Y)),\Mand P \in \CI([0, \ep); \Diff 3(Y))
\label{killing}
\end{multline}
and
\begin{multline}
[\Box, R] = -2 i \Box + a D_x + x^{-1} P,\text{with }\\ a\in \CI([0, \ep)\times
Y)\Mand P \in \CI([0, \ep); \Diff 2 (Y)).
\label{conformalkilling}\end{multline}
More generally,
\begin{lemma}\label{lemma:commutators}
For any $q,k\in \NN$,
\begin{multline}
\Box R^k=\sum\limits_{j=0}^kc_{k,j}R^j\Box+\sum\limits_{j=0}^{k-1}
\left(a_{k,j}D_x+\frac1xP_{k,j}\right)R^j,\Mwhere\\
c_{k,j}\in\bbC,\ a_{k,j}\in\CI([0,\ep)\times Y)\Mand
P_{k,j}\in\CI([0,\ep);\Diff2(Y)),
\label{melwun1.369}\end{multline}
\begin{multline}
\Box\Lap_Y^q =
\Lap_Y^q\Box+\sum\limits_{r=0}^{q-1}(Q_{q,r}D_x+\frac1xP_{q,r})\Lap_Y^r\\
\Mwith 
Q_{q,r}\in\CI([0,\ep);\Diff{q-r}(Y)),\ P_{q,r}\in\CI([0,\ep);\Diff{q-r+2}(Y)),
\label{melwun1.371}\end{multline}
and
\begin{multline}
\Box\Lap_Y^qR^k=\sum\limits_{j=0}^kc_{k,j}\Lap_Y^qR^j\Box
+\sum\limits_{\substack{0\le r\le q,\ 0\le j\le k\\ r+j<k+q}}
(Q_{k,l,q,r}D_x+\frac1xP_{k,l,q,r})\Lap_Y^rR^j\\
\Mwhere
Q_{k,l,q,r}\in\CI([0,\ep);\Diff{q-r}(Y)),\ P_{k,l,q,r}\in\CI([0,\ep);\Diff{q-r+2}(Y)).
\label{melwun1.372}\end{multline}
Furthermore, all the differential operators $P_{k,j}$, $P_{q,r}$,
$P_{k,l,q,r}$ have vanishing constant terms.
\end{lemma}
\begin{proof}
For $k=1$ \eqref{melwun1.369} follows directly from \eqref{conformalkilling}. For general
$k$ we may use this together with \eqref{melwun1.369}, as an inductive
hypothesis, to see that
\begin{equation*}
\Box R^{k+1}=R\Box R^k-[R,\Box]R^k
\label{melwun1.370}\end{equation*}
is of the stated form.  Equation \eqref{melwun1.371} follows similarly, and
\eqref{melwun1.372} from combining \eqref{melwun1.369} and
\eqref{melwun1.371}.
\end{proof}

It is also convenient to record a version of \eqref{melwun1.371} which holds
for real powers of $\Lap_Y$ in the tangential pseudodifferential calculus.
Let $Y_s=(1+\Lap_Y)^{s/2}$.
\begin{lemma}\label{lemma:morecommutators}
For all $s \in \RR$,
\begin{align*}
[\Box, Y_s] Y_{-s} &= P D_x + \frac 1 x Q, \\
 Y_{-s} [\Box, Y_s] &= P' D_x + \frac 1 x Q', \\
\end{align*}
with $P, P' \in \CI([0,\ep); \Psi^{-1}(Y))$ and $Q, Q' \in\CI([0,\ep);
\Psi^1(Y))$ and where $Q, Q'$ annihilate constants at $x=0$.
\end{lemma}

\section{Mapping properties of the Laplacian}\label{section:mapping}

By definition, $\Lap$ is symmetric as an operator on $\CI_c(X)$ with
respect
to the Riemannian volume form 
\begin{equation}
dg=x^{n-1}dx\, dh(x).
\label{melwun1.37}\end{equation}
To keep track of the weighted $L^2$ and Sobolev spaces which necessarily
appear here, we shall refer all weights to the intrinsic boundary
weights. These correspond to a non-vanishing positive smooth density on
$X^\circ$ with is of ``logarithmic'' form near the boundary 
\begin{equation}
\nu _b=a\frac{dx}xdh_0,\ 0<a\in\CI(X)\Mnear \pa X.
\label{melwun1.38}\end{equation}
Thus we set 
\begin{equation}
L^2_{\bo}(X)=\left\{u\in L^2_{\loc}(X^\circ);\int_X|u|^2\nu _b<\infty\right\}.
\label{melwun1.39}\end{equation}
Since the conic metric volume form is $dg=x^{n}a'\nu _b$ with $0<a'\in\CI(X),$
the metric $L^2$ space is
\begin{equation}
L^2_g(X)=\left\{u\in L^2_{\loc}(X);\int_X|u|^2dg<\infty\right\}=x^{-\frac n2}L^2_{\bo}(X).
\label{melwun1.40}\end{equation}

From \eqref{melwun1.36} it follows that 
\begin{equation}
\Lap=x^{-2}\left[(xD_x)^2 -\frac{i(n-1)+ xe}xD_x -i(n-2)xD_x+\Lap_h\right]\in x^{-2}\Diffb2(X).
\label{melwun1.41}\end{equation}
Here, $\Diffb*(X)$ is the filtered algebra of differential operators on $X$
which is the enveloping algebra of the Lie algebra $\Vb(X)$ of all smooth
vector fields on $X$ which are tangent to the boundary.

The weighted b-Sobolev spaces $x^pH_{\bo}^l(X)$ are essentially defined by
the mapping properties of these b-differential, and the corresponding
b-pseudodifferential, operators. They may also be defined directly using the
Mellin transform. Any weighted b-differential operator defines a continuous
linear map:
\begin{equation}
P\in x^r\Diffb m(X)\Longrightarrow P:x^pH^l_{\bo}(X)\longrightarrow
x^{p+r}H^{l-m}_{\bo}(X)\ \forall\ p,l\in\bbR.
\label{melwun1.42}\end{equation}

An \emph{elliptic} b-differential operator, \ie\ one for which the
characteristic polynomial in tangential vector fields is invertible
off the zero section of $\bTstar X,$ has the inverse property with respect
to regularity
\begin{multline}
P\in x^r\Diffb m(X),\ x^{-r}P \text{ elliptic },\ u\in
x^pH^{-\infty}_{\bo}(X),\ Pu\in x^{p+r}H^{l-m}_{\bo}(X)\\ \Longrightarrow u\in
x^pH^l_{\bo}(X).
\label{melwun1.43}\end{multline}
(See \cite{MR96g:58180} for a detailed discussion of b-differential and
-pseudodifferential operators.)

Such an elliptic operator is Fredholm as an operator \eqref{melwun1.42} for
all but a discrete set of $p\in\bbR.$ These correspond to the indicial
roots, those values of the complex parameter for which the indicial operator
is not invertible. The indicial operator is defined in general by 
\begin{equation}
P\in x^r\Diffb m(X)\Longrightarrow P(x^{is}v)=x^{is+r}(I(P,s)v+\cO(x)).
\label{melwun1.44}\end{equation}
This definition depends on the differential, at the boundary, of the
defining function chosen. Rather than carry the normal bundle information
to make this invariant we shall simply choose $x$ to be the defining
function in \eqref{melwun1.150}.

From \eqref{melwun1.36}, the indicial family of the Laplacian at a boundary
component $Y$ is 
\begin{equation}
I(\Lap,s)=\Lap_Y-i(n-2)s+s^2.
\label{melwun1.45}\end{equation}
If $0=\lambda_0<\lambda _1<\lambda _2\le\lambda _3\le \dots$ is the
sequence of eigenvalues of $\Lap_Y,$ repeated with multiplicity, then the
indicial roots of the Laplacian are
\begin{equation}
s^\pm_j=i\frac{(n-2)}2\pm\frac i2\sqrt{(n-2)^2+\lambda _j^2}.
\label{melwun1.46}\end{equation}

In general, assuming $P\in x^r\Diffb m(X)$ to be elliptic,  
\begin{multline}
P\text{ is Fredholm as an operator \eqref{melwun1.42}}\\
\Longleftrightarrow -p\text{ is not the imaginary part of an indicial root.}
\label{melwun1.47}\end{multline}
Thus, $0$ and $-n+2$ are always singular values of $p$ for the
Laplacian. For $n>2$ there is a gap in between:
\begin{equation}
(-n+2,0)\text{ is free of singular values.}
\label{melwun1.48}\end{equation}

Such a gap corresponds to boundary regularity for solutions. Thus 
\begin{multline}
P\in x^r\Diffb m(X)\text{ elliptic,}\ u\in x^pH^l_{\bo}(X),\ Pu\in
x^{q+r}H^{l-m}_{\bo}(X)\\
\Longrightarrow u\in x^qH^l_{\bo}(X)\Mprovided(p,q]\text{ is free of
singular values.}
\label{melwun1.49}\end{multline}

The conclusion of \eqref{melwun1.49} (where we assume that $p<q$ to avoid
triviality) does not follow if $(p,q]$ contains a singular value. For our
purposes it is enough to consider the special case that
\begin{equation}
p^*\in(p,q) \text{ is the unique singular value in }(p,q)\Mand
p^*<p+1,\ q<p^*+1.
\label{melwun1.50}\end{equation}
In fact let us further suppose that there is only one singular value $s\in\bbC$
such that  $p^*=-\Im s,$ and $I(P,s)$ is not invertible. Then
\begin{multline}
u\in x^pH^l_{\bo}(X),\ Pu\in x^{r+q}H_{\bo}^{l-m}\Longrightarrow u=u'+u'',\\
Pu'\in\dCI(X),\ u''\in x^qH^l_{\bo}(X),\
u'=\sum\limits_{0\le j\le k-1}x^{is}(\log x)^jv_j\chi+\tilde u'\Mand\\
\tilde u'\in x^qH^l_{\bo}(X)\Mwhere
I(P,xD_x)\left(\sum\limits_{0\le j\le k-1}x^{is}(\log x)^jv_j\right) =0.
\label{melwun1.51}\end{multline}
Here, $k$ is the order of $s$ as a pole of $I(P,s)^{-1}$ and $\chi =\chi
(x)$ is a cutoff which is identically $1$ near the boundary. If there are
several values of $s$ with $\Im s=-p$ at which $I(P,s)^{-1}$ is singular
then it is only necessary to add corresponding sums to \eqref{melwun1.51}.

\section{Domains and powers}
\label{section:domains}
Applying the general results above for b-differential operators to the
Laplacian we find
\begin{proposition}\label{melwun1.52} If $n\ge4$ then 
\begin{equation}
\Dom(\Lap)=\left\{u\in x^wL^2_{\bo}(X); \Lap u\in L^2_g(X)\right\}
\label{melwun1.53}\end{equation}
is independent of $w$ in the range 
\begin{equation}
-n+2<w<-\frac n2+2.
\label{melwun1.54}\end{equation}
If $n=3$ the same is true for $w$ in the range $-1<w<0$ and for $n=2$ 
\begin{equation}
\Dom(\Lap)=\left\{u\in L^2_g(X);u=c+u',\ c\in\bbC,\ u'\in x^wL^2_{\bo}(X),\ \Lap
u'\in L^2_g(X)\right\}
\label{melwun1.55}\end{equation}
is independent of $w$ for $w>0$ sufficiently small. In all cases, $\Lap$ is
an unbounded self-adjoint operator 
\begin{equation}
\Lap:\Dom(\Lap)\longrightarrow L^2_g(X)
\label{melwun1.56}\end{equation}
and $\Dom(\Lap)$ coincides with the domain of the Friedrichs extension; if
$n\ge3,$ then $\Lap$ is essentially self-adjoint.
\end{proposition}

\begin{proof} For $n\ge3$ the constancy of $\Dom(\Lap)$ in terms of $w$
follows from \eqref{melwun1.48} and \eqref{melwun1.49}. Thus, by
hypothesis, $u\in 
x^wL^2_{\bo}(X),$ with $w$ in the gap $(-n+2,0)$ and $Pu=x^2\Lap u\in
x^{-\frac n2+2}L^2_b.$ If $n>4$ then $-\frac n2+2$ also lies in the gap,
\eqref{melwun1.48}, so by \eqref{melwun1.43} and \eqref{melwun1.49} 
\begin{equation*}
n>4,\ u\in \Dom(\Lap)\Longrightarrow u\in x^{-\frac n2+2}H^2_{\bo}(X)
\end{equation*}
so then 
\begin{equation}
\Dom(\Lap)=x^{2-\frac n2}H^2_{\bo}(X),\ n>4. 
\label{melwun1.62}\end{equation}

For $n=4,$ $-\frac n2+2=0$ is the top of the gap whereas for $n=3$ the
gap is $(-1,0).$ In these cases we deduce
only that 
\begin{equation}
\Dom(\Lap)=\left\{u\in \bigcap_{w<0}x^wH^2_{\bo}(X):\ \Lap u\in L^2_g(X)\right\},\
n=3,4.
\label{melwun1.57}\end{equation}

On the other hand, for $n=2$ there is no gap. The hypothesis $\Lap u\in
L^2_g(X)$ is equivalent to $Pu=x^2\Lap u\in xL^2_g(X).$ The collapsed gap,
at $p=0,$ corresponds to a double root of $I(P,s)=\Lap_Y+s^2.$ Thus,
\eqref{melwun1.51} becomes 
\begin{multline}
u\in x^{-\epsilon}H^2_{\bo}(X),\ \Lap u\in L^2_g(X)\Longrightarrow\\
u=c+c'\log x+u'',\ u''\in x^{\epsilon}H^2_{\bo}(X),\ \epsilon >0\text{
sufficiently small.}
\label{melwun1.59}\end{multline}
In particular $\epsilon$ must be smaller than the smallest non-zero
eigenvalue  of $\Lap_Y$, $\lambda _1(\Lap_Y).$  The hypothesis
\eqref{melwun1.55} on $u\in\Dom(\Lap)$ is therefore just the vanishing of the
coefficient $c'$ in \eqref{melwun1.59}. It follows that \eqref{melwun1.55}
is also independent of $w$ provided $w<\lambda _1(\Lap_Y).$ 

We now demonstrate selfadjointness by showing that the unbounded
operators with the domains described above coincide with the Friedrichs
extension of $\Lap$.
By definition, $\Lap$ is associated to the Dirichlet form
\begin{equation}
  F(u,v) = \int_X \langle du,dv \rangle_g \, dg \,,\,\,
  u,v \in \CIc (X^\circ) 
\label{melwun1.181}\end{equation}
where $dg$ is the metric volume form \eqref{melwun1.37}.
The inner product in \eqref{melwun1.181} is that induced, by duality, by
the metric on $T^* X^\circ.$ Following Friedrichs we define
\begin{equation}
\dcal=\Dom(\Lap^{\frac12}) = \clos\left\{\CIc(X^\circ) \text{ w.r.t. }
F(u,u)+\|u\|^2_{L^2_g} \right\},
\label{Friedrichs}
\end{equation}
whenever $X$ is a compact conic manifold with boundary of dimension $n \geq
2.$
Then the Friedrichs extension of $\Lap$ is the unbounded operator
with domain 
\begin{equation}
\Dom(\Lap_{\text{Fr}})=\left\{u\in \dcal;\Lap u \in L^2_g (X)\right\}.
\label{melwun1.184}\end{equation}
where $\Lap$ is the bounded operator $\dcal \longrightarrow \dcal'$ and
$L^2_g(X)\subset \cD'$ is a well-defined subspace since $\dcal \subset
L^2_g(X)$ is dense.

The space $\Dom(\Lap^{\frac12})$ is independent of which conic metric on
$X$ is used to define it, since different conic metrics give equivalent
norms in \eqref{melwun1.181}. Moreover it can be localized by use of a
partition of unity. Thus it is the same, locally, as the same space on a
compact manifold blown up at a single point. This is well-known and easily
leads to the characterization
\begin{equation}
\Dom(\Lap^{\frac12})=
\begin{cases}
x^{-\frac n2+1}H^1_{\bo}(X)+\rho (x)\tilde \dcal& n=2,\\
x^{-\frac n2+1}H^1_{\bo}(X),& n\ge3,
\end{cases}
\label{melwun1.182}\end{equation}
where $\rho (x)$ is a cutoff near the boundary in terms of a boundary
defining function and $\tilde \dcal$ is the analogous domain of the
$1$-dimensional operator $D_x$, acting on $x^{-1} L^2_{\bo}([0, \infty))$.
Using Paley's theorem it is straightforward to characterize $\tilde
\dcal\subset x^{-1}L^2_{\bo}([0,\infty))$ in terms of the Mellin transform:
\begin{multline}
u\in \tilde \dcal\Longleftrightarrow u\in x^{-\delta}L^2_{\bo}([0,\infty),\
\delta >0,\Mand\\
u_M(x)=\int x^{is}u(x)\frac{dx}x\text{ is holomorphic in }-\infty<\Im s<0\Mwith\\
\int_{\Im s=r}|su_M(s)|^2d\Re s\text{ uniformly bounded in }r \in (-\infty,
0).
\label{melwun1.183}\end{multline}

In the definition of the Friedrichs domain, \eqref{melwun1.184}, the action
of $\Lap$ is distributional. Thus from \eqref{melwun1.53} and
\eqref{melwun1.54} it follows that
\begin{equation}
\Dom(\Lap_{\text{Fr}})=\left\{u\in x^{-\frac n2+1}H^1_{\bo}(X);\Lap u\in
x^{-\frac n2}L^2_{\bo}(X)\right\}=\Dom(\Lap),\ n>2.
\label{melwun1.186}\end{equation}

For $n=2$ the argument only needs to be modified slightly. It follows
directly from \eqref{melwun1.55} that
$\Dom(\Lap)\subset\Dom(\Lap_{\text{Fr}}).$ Moreover \eqref{melwun1.51}
shows that $u\in\Dom(\Lap_{\text{Fr}})$ has an expansion as in
\eqref{melwun1.55} except for the possibility of a logarithmic term. This
however is excluded by \eqref{melwun1.183} since it would correspond to a
double pole of the Mellin transform at $s=0.$

Thus in all case we have shown that the Friedrichs extension has domain
$\Dom(\Lap)$ as given in Proposition~\ref{melwun1.52}.
\end{proof}

We also need to describe the domains of the complex powers of $\Lap.$ For
integral powers it is straightforward to do so. Since
\begin{equation}
I(\Lap^k,s)= I(\Lap, s) I(\Lap, s+2i) \dots I(\Lap, s+2(k-1)i),
\label{melwun1.190}\end{equation}
it follows that the singular values of $\Lap^k$ are just the unions of the
shifts of those of $\Lap,$ with the possibility of accidental multiplicity
to be borne in mind. We are particularly interested in the domains of the
small powers with real part up to $n/4.$ Note that it follows from
arguments directly analogous to those above that 
\begin{equation}
\Dom(\Lap^p)=x^{-\frac n2+2p}H^{2p}_{\bo}(X)\Mfor p<\frac n4,\ 2p\in\bbN.
\label{melwun1.197}\end{equation}
For later applications we need to find the largest real $p$ for which this
remains true. To do so it is convenient to use complex interpolation.

\begin{lemma}\label{melwun1.187} For real $0\le p< n/4$ the
identification \eqref{melwun1.197} remains true.  The domain of
$\Lap^{n/4}$ is independent of the conic metric defining the Laplacian and is
given explicitly by
\begin{equation}
\Dom(\Lap^{\frac n4})=H^{\frac n2}_{\bo}(X)+\rho (x)\tilde \dcal_n
\label{melwun1.188}\end{equation}
where $\tilde \dcal_n$ reduces to \eqref{melwun1.183} in case $n=2$ and in
general is defined by 
\begin{multline}
u\in \tilde \dcal_n\Longleftrightarrow u\in
x^{-\delta}L^2_{\bo}([0,\infty)),\ \delta >0\Mand\\
u_M(x)=\int x^{is}u(x)\frac{dx}x\text{ is holomorphic in }
-\infty<\Im s<0\Mwith\\
\int_{\Im s=r}(|s|^2+1)^{\frac n2-1}
|su_M(s)|^2d\Re s\text{ uniformly bounded in } r\in (-\infty,0).
\label{melwun1.189}\end{multline}
There exists $\delta_0>0,$ depending on the metric and manifold, such that
for all $\delta \in (0, \delta_0)$,
\begin{equation}
\Dom(\Lap^{\frac n4+\delta})=x^{2\delta}H_{\bo}^{\frac n2+2\delta}(X)+\bbC
\label{melwun1.338}\end{equation}
and
\begin{equation}
\Dom(\Lap^p) \supset x^{-\frac n2+2p} \bH{2p} (X),\ \forall\ p\geq n/4;
\label{melwun1.331}\end{equation}
by duality there is a restriction map 
\begin{equation}
\Dom(\Lap^{-p})\longrightarrow  x^{-\frac n2-2p} \bH{-2p} (X).
\label{melwun1.330}\end{equation}
\end{lemma}
\noindent The map \eqref{melwun1.330} is \emph{not} injective for $p\ge\frac
n4$ since the inclusion \eqref{melwun1.331} does not then have dense range,
as follows from the computation of the domain of $\Lap^{n/4}.$

\begin{proof} We first use complex interpolation from the characterization of the
domain of $\Lap^{\frac 14}$ above. 
For two Banach spaces $X$ and $Y$, let $[X,Y]_{\theta}$ denote the complex
interpolation space at parameter $\theta$ as discussed, for example, in
\S1.4 of \cite{MR82i:35172}.
For any positive self-adjoint operator the
complex powers satisfy complex interpolation in the sense that 
\begin{equation}
\Dom(A^\theta)=[\Dom(A),L^2_g]_\theta.
\label{melwun1.198}\end{equation}
Furthermore the usual arguments with Sobolev spaces show that the weighted
b-Sobolev spaces exhibit the same interpolation property
$$
[x^{t} H^k_{\bo}(X), \bL(X)]_{\theta} = x^{t\theta}H^{k\theta}_{\bo}(X),\ 0\le\theta\le1.
$$
Applying this to \eqref{melwun1.197} we conclude that it remains true for
$p$ smaller than the greatest half-integer smaller than $n/4.$ For $n>2$ we
may apply the same argument again by noting that  
\begin{equation*}
\Dom(\Lap^q)=\left\{u\in\Dom(\Lap);\Lap u\in\Dom(\Lap^{q-1})\right\},\
1\le q<n/4.
\end{equation*}
This proves \eqref{melwun1.197} for all real $p<n/4.$

In fact essentially the same method applies to $\Dom(\Lap^{\frac n4})$ since
we have computed $\Dom(\Lap^{\frac n4-1}).$ The condition on the
Mellin transform in \eqref{melwun1.189} just represents decay at infinity
like $|s|^{-n/2}$ in a uniform $L^2$ sense except for the single factor of
$s$ which allows more general behavior at $s=0.$

The final characterization of the domains in \eqref{melwun1.188} now
follows by use of the Mellin transform near the boundary and reduces to the
same argument as in the two-dimensional case. Thus
\begin{equation}
\Dom(\Lap^{\frac n4})=\left\{u\in x^{-n/2}H^{\frac n2}_{\bo}(X),\ \Lap u\in
H^{\frac n2-2}_{\bo}(X)\right\}
\label{melwun1.192}\end{equation}
from which \eqref{melwun1.188} follows as before.

For $p\geq n/4$, the fact that $\Dom(\Lap^p) \supset x^{-\frac n2+2p}
\bH{2p} (X)$ follows from \eqref{melwun1.51} and \eqref{melwun1.338}
follows by similar arguments as those above, with the upper bound on $\delta$
arising from the first positive eigenvalue of the boundary Laplacian.
\end{proof}

The identification, \eqref{melwun1.338}, of the domain of the powers just
larger than $\frac n4$ gives a convenient mapping property independent of
dimension and of the conic metric involved
\begin{equation}
\Lap:x^{2\delta}H_{\bo}^{\frac n2+2\delta}(X)+\bbC\longrightarrow
x^{-2+2\delta}H^{\frac n2 -2+2\delta}_{\bo}(X),\ \delta >0\Msmall
\label{melwun1.339}\end{equation}
which has null space $\bbC$ and closed range which is a complement to $\bbC.$
Indeed $\Lap-\lambda$ also defines such a map for any $\lambda \in\bbC$
and this map is an isomorphism if $\lambda$ is not in the spectrum.

\section{Statement of the theorems}
\label{section:statement}
Let $\E_s$ denote the energy space, of order $s\in\bbR,$ of Cauchy data
for the wave equation:
$$
\E_s = \dcal_s \oplus \dcal_{s-1}\Mwhere \dcal_s = \Dom (\Lap^{s/2})
$$
as discussed above. If $\ocal$ is an open set containing a component of
$\pa X$, we will denote by $\dcal_s(\ocal)$ the corresponding local space
which is well defined since $\dcal_s$ reduces to $H^s(X)$ locally away from
the boundary. We define $\E_s(\ocal)$ similarly.

The Cauchy problem for the wave equation
\begin{gather}
(D_t^2-\Lap)u(t)=0,
\label{15.5.2000.1}\\
u(0)=u_0,\ D_tu(0)=u_1
\label{15.5.2000.3}\end{gather}
has a unique solution 
\begin{equation}
u\in\cC^0(\bbR;\dcal_{s})\cap\cC^1(\bbR;\dcal_{s-1})
\label{15.5.2000.2}\end{equation}
for all $(u_0,u_1)\in \E_s.$ Similarly the inhomogeneous forcing problem 
\begin{equation}
(D_t^2-\Lap)v(t)=f,\ f\in\CmI(\bbR;\cD_{s-1}),\ f=0\Min t<0
\label{melwun1.332}\end{equation}
has a unique forward solution $v\in\CmI(\bbR;\cD_{s})$ with $v=0$ in $t<0.$

If the conic manifold is complete rather than compact, away from the conic
ends, then a similar result holds. In the general case of a manifold with
conic ends the wave equation \eqref{15.5.2000.1} has a unique solution
with compact support in a finite interval $[-T,T]$ in place of $\bbR$ in
\eqref{15.5.2000.2} provided the initial data has compact support. Since
the results below are all local near the boundary and the general case can
be reduced to this one, for simplicity of presentation we consider only the
case of a compact conic manifold.

\begin{definition}\label{def:admissible}
An \emph{admissible} solution to \eqref{15.5.2000.1} is one of the form
\eqref{15.5.2000.2} for some $s \in \RR$ with the equation holding in
$\CmI(\bbR;\dcal_{t})$ for some $t\in\bbR.$ 
\end{definition}
\noindent
We deal here with admissible solutions to the wave equation,
\ie\ solutions corresponding to the Friedrichs realization of the Laplacian.

The time-translation invariance of the wave equation means that if $u$ is
an admissible solution and $e\in\CmIc(\bbR)$ then $e(t)* u$ is also an
admissible solution. We shall use this below to decompose solutions into
positive and negative parts by choosing a decomposition 
\begin{multline}
\delta (t)=e_+(t)+e_-(t)+e_\infty(t),\ e_\pm\in\mathcal{S}'(\bbR),\\
e_\infty\in\mathcal{S}(\bbR),\ \WF(e_\pm)=\{(0,\pm\infty)\}\subset S^*\bbR.
\label{melwun1.387}\end{multline}
The corresponding decomposition of an admissible solution $u$ is then 
\begin{equation}
u=u_++u_-+u_\infty,\ u_\infty\in\CI(\bbR;\dcal_{\infty}),\
u_\pm=e_\pm(t)* u,
\label{melwun1.388}\end{equation}
where all three terms are admissible solutions. Typically we choose $e_+$
and $e_-$ to have Fourier transforms supported in $[1,\infty)$ and
$(-\infty,-1].$

Global regularity theory for the wave equation shows that the strongly
continuous group of bounded operators defined by \eqref{15.5.2000.1} --
\eqref{15.5.2000.2}
\begin{equation}
U(t):\dcal_{1}\oplus L^2_g\longrightarrow \dcal_{1}\oplus L^2_g
\label{15.5.2000.4}\end{equation}
satisfies 
\begin{equation}
U(t):\E_s \longrightarrow \E_s \ \forall\ s\in\bbR,
\label{15.5.2000.5}\end{equation}
by continuous extension for $s<0.$ The solution to \eqref{melwun1.332} is
then given by Duhamel's principle.

It is useful at various point in the discussion below to change the degree
of regularity of an admissible solution; this can always be accomplished by
convolution in $t.$

\begin{definition}\label{Theta-op}
Let $\Theta_s$ be the operator on $\RR$ given by
\begin{equation}
\kappa (\Theta_s)(t,t') = \psi(t-t') \kappa(\abs{D_t}^s)(t,t')
\label{Theta}
\end{equation}
where $\psi (t)$ is a smooth function of compact support, equal to one near
$t=0$ and $\kappa$ denotes Schwartz kernel.  
\end{definition}

\begin{lemma}\label{lemma:orderchange}
Let $u \in \mathcal{C}(\RR; \E_r)$ be a solution to the wave equation.
Then $\Theta_s u \in \mathcal{C}(\RR; \E_{r-s})$ is also a solution, and
for all $s \in \RR$,
$$
\Theta_{s} \Theta_{-s} u = u \bmod \CI(\RR; \E_\infty).
$$
\end{lemma}

Away from the boundary the wave operator $D_t^2-\Lap$ is smooth with
principal symbol $\tau ^2-|\zeta|^2_z$ at the point $(t,z;\tau,\zeta)$ in
terms of the canonical coordinates associated to the coordinates $(t,z).$
H\"ormander's theorem on the propagation of wavefront set for operators
of real principal type therefore applies and shows that microlocal
regularity, in terms of the wavefront set relative to Sobolev spaces, is
contained in the characteristic variety and is constant along the null
bicharacteristic which foliate it. Thus the regularity of a solution at any
point of $T^*(\bbR\times X^{\circ})\setminus0$ is readily describable in
terms of the regularity of the initial data unless the null
bicharacteristic through the point hits the boundary at some intervening time.

We first state a diffractive theorem, which simply says that if there
are no incoming singularities at the boundary component $Y$ at some time
$\bar t$ then there are no outgoing singularities arising at that time.
Moreover this regularity can be microlocalized in $\tau$, the dual variable
to $t$.

\begin{theorem}[Diffractive regularity]\label{15.5.2000.7} If $u$ is an
  admissible solution to the conic wave equation and for some small
  $\epsilon >0,$ $s\in\bbR$, $R_{\pm,I}^\ep (\bar t,
  Y)\cap\WF^s(u)=\emptyset,$ then $R_{\pm, O}^\ep(\bar t, Y) \cap
  \WF^s(u) = \emptyset.$ If $R_{I}^\ep(\bar t, Y)\cap \WF^s (u)
  =\emptyset,$ then for some open sets $I \ni \bar t$ and $\ocal
  \supset Y,$ $u\in\ccal(I;\dcal_s(\ocal)).$
\end{theorem}
By time-reversibility of the equation, this theorem implies
Theorem~\ref{intro:diffthm}.  It is proved in
Section~\ref{section:diffractive}.

In order to state the geometric theorem, we need to introduce a
second-microlocal condition on the incoming singularities at time $t=\bar
t.$

\begin{definition}[Nonfocusing conditions]\label{NC} We say that an admissible
solution $u$ to the conic wave equation satisfies the \textbf{nonfocusing
condition} at time $\bar t$ and at boundary component $Y,$ with background
regularity $r$ and relative regularity $l,$ if for some positive integer $N$
and small $\epsilon >0$
\begin{gather}
\WF^{r}(u)\cap R_{\pm, I}^\ep(\bar t, Y) = \emptyset\Mand 
\label{melwun1.379}\\
\WF^{r+l}\left((1+\Lap_Y)^{-N}u\right)\cap R_{\pm, I}^\ep (\bar t,Y)
= \emptyset.
\label{melwun1.334}\end{gather}
If, for some $N=N(k)\in \NN,$
\begin{equation}
\WF^{r}\left((x D_x + (t-\bar t)D_t)^p(1+\Lap_{Y})^{-N} u\right)
\cap R_{\pm, I}^\ep (\bar t,Y) = \emptyset,\ p\le k,
\label{melwun1.359r}\end{equation}
we say that $u$ satisfies the \textbf{radial regularity condition}
to order $k$ (and with regularity $r).$ The combination of \eqref{melwun1.379}
with \eqref{melwun1.359r}, but for regularity $r+l$:
\begin{equation}
\WF^{r+l}\left((x D_x + (t-\bar t)D_t)^p(1+\Lap_{Y})^{-N} u\right)
\cap R_{\pm, I}^\ep (\bar t,Y) = \emptyset,\ p\le k,
\label{melwun1.359}\end{equation}
will be called the \textbf{conormal nonfocusing condition} (to order $k$
with background regularity $r$ and relative regularity $l$).
\end{definition}

These conditions measure the extent to which the solution contains a wave
collapsing radially onto the boundary, up to relative regularity order $l.$
The strongest, conormal, version asserts that after smoothing $u$ in the
tangential variables a conormal estimate at the surface $x=\bar t-t$ holds
to order $k.$ As an important example, the fundamental solution $\sin
t\sqrt\Lap/\sqrt\Lap$ with pole close to the boundary satisfies the
conormal nonfocusing condition with $r<-n/2+1$ and $l<(n-1)/2$ for any $p$
(see \S\ref{section:divthm1}).

\begin{theorem}[Geometric propagation]\label{15.5.2000.8} Let $u$ be an
admissible solution to the conic wave equation.  If $p \in R_{\pm, O}^\ep
(\bar t, Y),$ $\Gamma^\ep (p) \cap \WF^{r'}(u) = \emptyset$ for small
$\epsilon >0$ and $u$ satisfies the nonfocusing condition of
Definition~\ref{NC} with regularity $r+l>r',$ then $p \notin
\WF^{r'}(u)$.
\end{theorem}
This theorem represents a sharpening of Theorem~\ref{melwun1.318}.  We
prove a weaker version of it in Section~\ref{section:genprop} and
then obtain the full theorem in Section~\ref{section:conormal}.

The ``edge structure'' on the product $\bbR_t\times X$ is discussed in the
next section and in particular the related scale of weighted edge Sobolev
spaces is defined there. These spaces are appropriate for the description
of the boundary regularity of admissible solutions to the wave equation.  A
crucial role in our proof of the geometric propagation theorem above is
played by the following result on decay relative to weighted edge Sobolev
spaces.

\begin{theorem}[Division theorem]\label{melwun1.34} If $u$ is an
admissible solution to the conic wave equation satisfying the nonfocusing
condition in Definition~\ref{NC} with $r+l< n/2$ then there are open
sets $I \ni \bar t$ and $\ocal \supset Y$ in $X$ such that
\begin{equation}
u\in x^{r+l-n/2}H^{r+l-k}_{\eo}(I \times \ocal).
\label{melwun1.35}\end{equation} If $u$ satisfies the conormal nonfocusing
condition to order $k=1$ and if $r+l\leq 1$ then in addition
$$
D_t u\in x^{r+l-n/2-1/2-\ep}H^{r+l-k-1}_{\eo}(I \times \ocal)\ \forall\ \ep>0.
$$
\end{theorem}
We prove this theorem in Section~\ref{section:divthm1}.

One can use energy conservation and the diffractive theorem to show that
any solution to \eqref{15.5.2000.1}--\eqref{15.5.2000.2} is in $x^{s-n/2}
H_{\eo}^{s}$ as long as $-n/2<s<n/2$ and that any solution satisfying
$\WF^r(u)\cap R^\ep_{\pm, I}(\bar t)=\emptyset$ is in $x^{r-n/2} \eH{r}
(\RR\times X),$ locally in time near $t=\bar t$, as long as $-n/2<r<n/2.$ The
nonfocusing condition thus leads to a stronger growth estimate than is given
by energy estimates alone.

The fact that the division theorem yields a stronger result in the presence
of a (first-order) conormality assumption is closely related to the fact
that if a solution to the wave equation is conormal with respect to the incoming
surface $x-t=\bar x$ in $t<0$ it is conormal with respect to the
corresponding outgoing surface $x=\bar x + t$ in $t>0.$ This follows in
turn from the corresponding result on the Cauchy problem.

\begin{theorem}[Conservation of conormality]\label{thm:conormality} Let
$(u_0, u_1)\in\E_s$ be conormal with respect to the hypersurface $\{x=\bar
x\}$ for some $\bar x>0$ sufficiently small, then in $x>0,$ the solution to
\eqref{15.5.2000.1}--\eqref{15.5.2000.3} is conormal with respect to
\begin{equation*}
\{x=\abs{\bar x-t}\}\cup \{x=\bar x+t\}
\end{equation*}
for $0<t<2\bar x.$
\end{theorem}
The proof of this theorem appears in Section~\ref{section:radial}.

A microlocalized version of Theorem~\ref{thm:conormality} allows us to
show the conormality of the diffracted front, subject to radial regularity
of the incident wave.

\begin{theorem}[Conormality of the diffracted front]\label{melwun1.367}
Let $u$ be an admissible solution to the conic wave equation satisfying the
radial regularity condition in Definition~\ref{NC} at $\bar t, Y$ to every
order $k \in \NN_0$ and suppose that for some small $\epsilon >0$ and some
$p\in R^\ep_{O}(\bar t, Y),$ $\Gamma^\ep (p) \cap \WF(u)=\emptyset$
then, microlocally near the outgoing bicharacteristic through $p,$ $u$ is
conormal with respect to $x=t-\bar t.$
\end{theorem}
The proof of this theorem is in Section~\ref{section:genprop}.
As a corollary of this result we may refine Theorem~\ref{15.5.2000.8},
concluding that if in addition the full conormal nonfocusing condition,
\eqref{melwun1.359}, holds for all $k,$ then the solution is conormal to
the surface $x=t-\bar t$ near $p.$

Although it follows from the results above, we nonetheless restate
part of the result on the fundamental solution discussed in the
introduction; in fact, the proof of this result, which occurs
Section~\ref{section:conormal}, is crucial in the proof of the full
version of Theorem~\ref{15.5.2000.8}.
\begin{theorem}[Fundamental solution]\label{melwun1a.374} Let $E _{\bar
m}$ be the fundamental solution to the conic wave equation with pole
at $\bar m=(\bar x, \bar y)\in X^\circ.$ If $\bar x$ is sufficiently
small then $E_{\bar m}$ is conormal with respect to $\{x+\bar x =t\}$
away from the wave cone emanating from $\bar m$, and is of Sobolev
order $\frac12-\delta$ there for any $\delta >0.$
\end{theorem}
This is equivalent to Theorem~\ref{melwun1.374} of the Introduction.

\section{Edge pseudodifferential calculus}
\label{section:calculus}

The edge calculus of pseudodifferential operators was introduced by
Mazzeo~\cite{Mazzeo4} as a class of operators on any compact manifold $M$
with boundary having a fibration $\pa M,$ $\pi: \pa M \to N$
with fiber $F.$ In this paper, $M=\RR\times X$, $N=\RR$, and $\pi$ is the
product fibration
$$
\pi:\RR\times \pa X \to \RR,\ \pi(t,p)=t.
$$
The noncompactness of $M$ in the situation at hand necessitates only minor
changes to the calculus, namely keeping supports proper. Our operators
generally have compactly supported kernels. In this section, we discuss the
edge calculus in the general setting, as this involves no increase in
complication over the special case of $\RR\times X.$ Although there is no
treatment in \cite{Mazzeo4} of edge microsupport or wavefront set, the
properties of these objects follow easily from the properties of the
calculus discussed in \cite{Mazzeo4} much as the properties of the
conventional wavefront set and microsupport follow from the properties of
the pseudodifferential calculus on closed manifolds. First we give a
brief synopsis of the edge calculus and its properties.

A Lie algebra of \ci\ vector fields $\Ve (M),$ associated to $\pi,$ is given by
$$
V \in \Ve \Longleftrightarrow V \text{ is tangent to the fibers of }
\pi \text{ at } \pa M.
$$
If $(x,y,z)$ are local coordinates with $x$ a defining function for $\pa
M$, $y$ coordinates on $N$ lifted and extended to functions on $M$,
and $z$ restricting to coordinates on the fibers, then $\Ve(M)$ is
locally spanned over $\CI (M)$ by
$$
x D_x,\ x D_y\Mand D_z.
$$
Thus, there exists a vector bundle $\Te M$ (the \emph{edge tangent bundle}) such that
$$
\Ve(M) = \CI (M; \Te M).
$$
Let $\Testar M$ (the \emph{edge cotangent bundle}) denote the dual of $\Te M$;
sections of $\Testar M$ are then locally spanned over $\CI(M)$ by $dx/x,$
$dy/x$ and $dz.$ By the \emph{edge cosphere bundle} we mean the quotient
$$
\Sestar M =\left(\Testar M\setminus0\right)/\RR_+.
$$

 There is a canonical bundle map $\Testar M \to T^* M$ since edge vector
fields are smooth up to $\pa M.$ This map is an isomorphism over $M^\circ$
so over the interior $\Testar M$ is a symplectic manifold, with symplectic
form given as usual by the exterior derivative of its canonical one-form;
this form becomes singular at the boundary.

For $k\in \NN$, let $\Diffe k(M)$ be the space of differential operators
spanned over $\CI(M)$ by operators $V_1\dots V_m$, $V_i \in \Ve (M)$,
$m\leq k$.  There exists a canonical (principal) symbol map, $\esigma[m],$
associating to $P\in\Diffe m (M)$ the polynomial function on the fibers of
$\Testar M$ extending the usual symbol map over the interior; it gives a
surjective map
$$
\esigma[m]:\Diffe m (M)\longrightarrow \text{homogeneous polynomials of
degree $m$ on } \Testar M
$$
with null space precisely $\Diffe{m-1} (M).$ 

In the particular case of interest for this paper, $\Ve (\RR\times X)$ is
locally spanned by the vector fields $x D_x$, $x D_t$, $D_y$ over $\CI
(\RR\times X).$ If we let
$$
\xi \frac{dx}x + \lambda \frac{dt}x + \eta \cdot dy
$$
be the canonical one-form on $\Testar (\RR\times X)$, then naturally
$$
\esigma[1] (x D_x) = \xi,\ \esigma[1] (x D_t) = \lambda,\ \esigma[1]
(D_{y_i}) = \eta_i
$$
are dual coordinates on the fibers of $\Testar M.$
 
The edge calculus of pseudodifferential operators, defined in
\cite{Mazzeo4}, arises as a microlocalization of $\Diffe{*}(M).$  Let
$\dCI(M)$ denote the space of smooth functions on $M$ vanishing, with all
derivatives, to infinite order at $\pa M$, and let $\CmI(M)$ be the dual
to the corresponding space of densities, $\dCI(M;\Omega).$ The space
$\ePs{*}(M)$ is a graded algebra of operators on $\CmI (M)$; as we
will frequently use \emph{weighted} edge operators, we will depart from the
notation of \cite{Mazzeo4} in carrying along the weight as an index in the
calculus.  Thus the bigraded space of operators $\ePs{m,l}(M)$ enjoys the
following properties:
\begin{itemize}
\item
$\ePs{m,l}(M)$ is a graded $*$-algebra.
\item
$\ePs{m,l}(M) = x^l \ePs{m,0}(M)$ (and the latter space is the space
denoted $\ePs m(M)$ in \cite{Mazzeo4}).
\item
$x^l \Diffe m (M) \subset \ePs{m,l} (M)$ for all $m\in \NN$, $l\in \ZZ$.
\item
The maps $\esigma[m]$ extend to
$$
\esigma[m,l]: \ePs{m,l} (M)\longrightarrow x^l \left[ S^m_{\phg} (\Testar M)/
S^{m-1}_{\phg} (\Testar M)\right];
$$
the range space for $\esigma$ can be conveniently identified with
$\CI(\Sestar M).$
\item The symbol map is a homomorphism of $*$-algebras.
\item
The sequence
$$
0\longrightarrow \ePs{m-1,l}(M)\longrightarrow \ePs{m,l}(M)\longrightarrow
x^l \left[ S^m_{\phg} (\Testar M)/ S^{m-1}_{\phg} (\Testar
M)\right]\longrightarrow 0
$$
is exact and multiplicative.
\item
If $A\in\ePs{m,l}(M)$ and $B\in\ePs{m',l'}(M)$ then
$$ 
\esigma[m+m'-1,l+l']([A,B])=\frac 1i \{\esigma[m,l](A), \esigma[m',l'](B)\},
$$
where the Poisson bracket is computed with respect to the singular symplectic
structure on $\Testar M$ described above.
\item
$\esigma (x^{-l} A x^l) = \esigma (A)$ for all $A \in
\ePs{\mbox{}} (M)$, $l\in \RR$.
\end{itemize}

In the case $M=\RR\times X,$ the elements of $\ePs{m,l}(\RR\times X)$ may
be represented locally in the form
\begin{multline}
A u(t,x,y)=(2\pi)^{-n-2}x^l\int e^{i(s-1)\xi+(y-y')\eta+iT\lambda}\\
b(t,x,y,s,y',T,\xi,\eta,\lambda)
u(xs,y',t-xT)x^{2+l} ds\, dy'\, dT, 
\label{melwun1.335}\end{multline}
where $b$ is a classical (polyhomogeneous) symbol of order $m$ with $\xi,$
$\eta,$ $\lambda$ as fiber variables; $\esigma[m](A)$ is then the equivalence
class of $b.$

An operator $A\in \ePs{m,l}(M)$ is said to be \emph{elliptic} at $p \in
\Sestar M$ if $\sigma (A)$ has an inverse in $x^{-l} [S^{-m}_{\phg}
(\Testar M)/S^{-m-1}_{\phg}(\Testar M)]$, locally near $p.$ Assuming that
the orders are clear we may suppress indices and so $p \in \Elle (A).$
There is a related notion of microsupport for edge pseudodifferential
operators, corresponding to the fact that the composition of operators
gives an asymptotically local formula for the amplitude, $b$ in
\eqref{melwun1.335}. If $A\in\ePs{m,l}(M),$ the \emph{microsupport} of $A$,
$\eWFprime(A)$, is the closed subset of $\Sestar M$
given locally by $\esssupp b$, the conic support of $b;$ it has the
following properties:

\begin{itemize}
\item
$\eWFprime(AB)\subset \eWFprime(A)\cap\eWFprime(B).$
\item
$\eWFprime(x^{-s} A x^s)=\eWFprime(A)$ for all $s \in \RR$.
\item
If $A \in \ePs{m,l} (M)$ and $p\in\Elle(A)$ there exists
$Q\in\ePs{-m,-l}(M)$ such that
$$
p \notin \eWFprime (QA-I)\cup \eWFprime(AQ-I).
$$
\item
If $A \in \ePs{m,l}(M)$ and $\eWFprime(A) = \emptyset$ then $A \in
\ePs{-\infty,l}(M)$---note that this is not a totally residual operator.
\end{itemize}

There is a continuous quantization map (by no means unique)
$$
\Ope: x^l S^m_{\phg}(\Testar M) \to \ePs{m,l}(M)
$$
which satisfies 
\begin{gather*}
\esigma[m,l] (\Ope(a)) = [a] \in x^l S^{m}_{\phg}(\Testar M)/S^{m-1}_{\phg}(\Testar M)\
\forall\ a \in x^l S^m_{\phg} (\Testar M)\Mand\\
\eWFprime \Ope(a) \subset \esssupp(a).
\end{gather*}

Associated with the edge calculus there is a scale of Sobolev spaces. For
integral order these may be defined directly. Thus for $k\in\bbN$
and any $s\in\bbR$ we set
\begin{multline}
H^{k,s}_{\eo}(\bbR\times X)=\{u\in x^sL^2_{\bo,\loc}(\bbR\times X);\\
(x D_t)^rP_{k-r}u\in  x^sL^2_{\bo,\loc}(\bbR\times X)\ \forall\
P_{k-r}\in\Diffb{k-r}(X)\Mand 0\le r\le k\},\ k\in\bbN.
\label{melwun1.32}\end{multline}
For negative integral orders we can similarly define
\begin{multline}
\Mfor k\in-\bbN,\ H^{k,s}_{\eo}(\bbR\times X)\ni u\Longleftrightarrow
u\in \CmI(\bbR\times X) \\ \Mand \exists\
u_{i,r}\in x^sL^2_{\bo,\loc}(\bbR\times X),\ 
P_{k-r,i}\in\Diffb{k-r}(X),\ i=1,\dots,N,\\
\Mwith u=\sum\limits_{i=1}^N\sum\limits_{r=0}^k(x D_t)^rP_{k-r,i}u_{i,r}.
\label{melwun1.33}\end{multline}
For general orders, the edge Sobolev spaces can be defined using the calculus.

\begin{definition}
$u \in \eH{m,l} (M) \Longleftrightarrow \ePs{m,-l}(M)\cdot u\subset\bL(M).$
\end{definition}
\noindent Note that we have chosen to weight these Sobolev spaces with
respect to the b-weight, not the metric weight.  Note also the change of
sign on $l$.  Since $\eH{m,l}(M) = x^l \eH{m,0}(X)$, we will often use the
notation $\eH{m}(M) = \eH{m,0}(M)$ and write the $x$-weight explicitly.  In
the case of interest in this paper, when $M=\RR\times X$ is noncompact, we
will consider only edge Sobolev spaces \emph{local in $t$}, without writing
this explicitly. The corresponding $L^2$-based edge wavefront set plays a
fundamental result below.

\begin{definition}
If $u \in \eH{-\infty, l}(M)$ then $\eWF[m,l](u) \subset\Sestar M$ is
defined by the condition that $p \notin \eWF[m,l](u)$ iff there exists $A \in
\ePs{m,-l}(M)$, elliptic at $p$, such that $Au \in \bL(M).$
\end{definition}

The usual properties carry over to these spaces:
\begin{itemize}
\item $\eWF[m,l](M)$ is closed.
\item For all $l \in \RR$, $\eWF[m,l](u) \cap \Sestar M^\circ = \WF^m(u)$
(recall that $\Sestar M^\circ$ and $S^* M^\circ$ are canonically
isomorphic).
\item
$
\bigcap_{m,l} \eH{m,l}(M) = \dCI(M),\quad \bigcup_{m,l} \eH{m,l}(M) = \CmI(M).
$
\item
Complex interpolation holds:
\begin{equation*}
[\eH{m,l}(M), \eH{m',l'}(M)]_\theta=
\eH{\theta m+ (1-\theta) m',\ \theta l + (1-\theta) l'}(M).
\end{equation*}
\item
If $A \in \ePs{m,l} (M)$ then $A: \eH{m',l'} (M)\longrightarrow
\eH{m'-m,l'+l}(M)$.
\item
For $m\leq m'$,
$$
\eWF[m,l](u) \subset \eWF[m',l](u).
$$
\item
If $u \in \eH{-\infty, l'}(M)$ and $A \in \ePs{k,l}(M)$ then 
\begin{equation*}
\eWF[m,l+l'](A u) \subset \eWFprime(A)\cap \eWF[m+k,l'](u).
\end{equation*}
\item
If $u \in \eH{-\infty, l'}(M)$ and $A \in \ePs{k,l}(M)$ then
$$
\eWF[m+k,l'](u)\backslash \eWF[m,l+l'](A u)\subset (\Elle A)^\complement.
$$
\end{itemize}

We now prove a less standard result.

\begin{proposition}
If $u \in \eH{-\infty,l}(M)$ and $l'\leq l$ then for all $\theta \in [0,1]$
and $m,m'\in\bbR,$
$$
\eWF[\theta m+(1-\theta) m', \theta l+(1-\theta) l'](u)\subset
\eWF[m',l'](u)\cap \eWF[m,l](u).
$$
\label{prop:wfinterp}\end{proposition}
\begin{proof}
Suppose $p \notin \eWF[m',l'](u)\cap\eWF[m,l](u).$  Then there exists $A'
\in \ePs{m',-l'}(M)$ and $A\in\ePs{m,-l}(M)$ both 
elliptic at $p $, such that $Au,$ $A'u\in \bL.$ From elliptic regularity it
follows that if $B\in\ePs{0,0}(M)$ has $\eWFprime$ concentrated near $p$
then $Bu\in\eH{m,l}(M)\cap\eH{m',l'}(M).$  By interpolation of weighted
edge Sobolev spaces it follows that $Bu \in \eH{\theta(m-m'),\theta(l-l')}(M)$
and the result follows.
\end{proof}

We will require some results about edge-regularity of solutions to the wave
equation.

\begin{proposition}\label{prop:sobolevrelationship}
For all $p\in \RR_+,$ 
\begin{gather*}
x^p \eH{p} (\RR\times X) = H^p_{\loc}(\RR; \bL(X))\cap
x^p \Lloc(\RR; \bH{p}(X))\Mand\\
x^{-p} \eH{-p} (\RR\times X) =
H^{-p}_{\loc} (\RR; \bL(X)) + x^{-p}\Lloc(\RR; \bH{-p}(X)).
\label{melwun1.365}\end{gather*}
\end{proposition}

\begin{proof}
All spaces are local in $t$ (by definition); multiplying by any
$\phi\in\CIc(\bbR)$ it suffices to assume that supports are
compact. Fourier transformation in $t$ coupled with an interpolation
argument shows that for $p\in\bbN$
$$
H^p_{\loc}(\RR; \bL(X)) \cap x^p \Lloc(\RR; \bH{p}(X))=\bigcap_{j=0}^p x^j
H^{p-j}_{\loc}(\RR; \bH{j}(X)).
$$
The latter space is equal to $x^p \eH p (\RR\times X)$ as defined in
\eqref{melwun1.32}, proving the first result for $p \in \NN$; the second
follows by duality.  The results for general $p\in \RR$ follow by
interpolation.
\end{proof}

\begin{proposition}\label{prop:domainsinedgeterms}
If $u$ is a solution to the wave equation in $\Lloc(\RR; \dcal_s)$, with
$\abs{s} <n/2$, then $u \in x^s \eH{s-n/2}(\RR\times X).$
\end{proposition}

\begin{proof}
Certainly $u \in H^{2r}_{\text{loc}}(\RR;\dcal_{s-2r})$ for all $r \in \ZZ.$
So, by interpolation, this holds for all $r \in \RR.$ By
Lemma~\ref{melwun1.187}, $u \in \Lloc(\RR; x^{s-n/2} \bH{s}(X)) \cap 
H^s_{\text{loc}}(\RR; \bL(X))$, hence the result follows from
Proposition~\ref{prop:sobolevrelationship}.
\end{proof}

\section{Bicharacteristic flow}\label{section:flow}
The canonical one-form on $\Testar (\RR\times X)$ is
$$
\lambda \frac{dt}x + \xi \frac{dx}x +  \eta \cdot dy
$$
hence the symplectic form is
$$
\omega=\frac{d\lambda \wedge dt}x +\frac{d\xi \wedge dx}x -\frac{\lambda
  dx\wedge dt}{x^2} + d\eta \wedge dy.
$$

We can now write the symbol of the d'Alembertian
$$
p = \sigma(\Box) = \lambda^2/x^2- g(x,y,\xi dx/x +\eta\cdot dy)
=\frac{\lambda^2-\xi^2-h(x,y,\eta)}{x^2}
$$
and 
$$
p_0 =\sigma (D_t^2-\Lap_0) = \frac{\lambda^2-\xi^2-h_0(y,\eta)}{x^2}
$$
where $\Delta_0$ is defined by \eqref{modellap}. Let $\Sigma=\{p=0\}
\subset \Testar(\RR\times X)$ denote the characteristic variety of the
d'Alembertian.

Let $H_g$ and $H_{g_0}$ denote the respective Hamilton vector fields of $p$
and $p_0$ on $\Testar(\RR\times X)$, near a boundary component $Y$ of $X$.
Thus
\begin{equation}
\frac{x^2}2 H_g= \frac{x^2}2 H_{g_0}+ W = H_Y+ (\xi^2+h_0(y,\eta)) \pa_\xi + \lambda \xi
\pa_\lambda + \xi x\pa_x - \lambda x \pa_t + W,
\label{ham1}
\end{equation}
where $W$ is the Hamilton vector field of $p-p_0$, hence
\begin{equation}
\frac{x^2}2 W = \frac x2 \frac{\pa h(\eta)}{\pa x} \pa_\xi
  +(h-h_0) \pa_\xi - \h \frac{\pa (h(\eta)-h_0(\eta))}{\pa y} \cdot
  \pa_\eta + (h^{ij}-(h_0)^{ij})\eta_i \pa_{y_j},
\label{ham2}
\end{equation}
and where $H_Y$ is the Hamilton vector field in $(y,\eta)$ for
$(1/2)h_0(y,\eta)$, \ie\ is the geodesic spray in $Y.$ 

Note that $H_g$, $H_{g_0}$, and $W$ are all homogeneous of degree $1$ in
$\Testar (\RR\times X)$, and that $(x^2/2) W$ is less singular than
$(x^2/2) H_g$ at $x=0$, and vanishes at $\eta=0$.  Thus if $\Tebarstar
(\RR\times X)$ denotes the fiberwise radial compactification of $\Testar
(\RR\times X)$,
\begin{align}
(x^2/2) H_g,\ (x^2/2) H_{g_0} \in & h(\eta)^{\h}\mathcal{V}_b(\Tebarstar (\RR\times
X)\backslash 0),\\
(x^2/2) W \in x h(\eta)^\h &\mathcal{V}_b(\Tebarstar (\RR\times
X)\backslash 0).
\label{lowerorder}
\end{align}

Note also that the vector field $H_g$ is tangent to the incoming and
outgoing sets $R_{\pm, I}(Y),\ R_{\pm, O}(Y)$, which we now regard (by
homogeneity) as subsets of $\Sestar (\RR\times X)$.  These incoming and
outgoing manifolds are the interiors of smooth manifolds with boundary in
$\Sestar (\RR\times X)$, and we define their boundaries as follows, with
$\IC$ and $\OG$ standing for incoming and outgoing manifolds respectively.
\begin{definition}
Let
\begin{align*}
\IC_{\pm}(\bar t, \bar y) &= \overline{R_{\pm, I}^\ep(\bar t,\bar y)}\cap \Sestar[\RR\times \pa X]
(\RR\times X)\\
\OG_{\pm}(\bar t, \bar y) &= \overline{R_{\pm, O}^\ep(\bar t,\bar y)}\cap \Sestar[\RR\times \pa X]
(\RR\times X)
\end{align*}
with the same convention for omitted indices as was used for $R^\ep_{\pm,
\bullet}(\bar t, \bar y)$.
\end{definition}

For $Y$ a boundary component of $X$, we define a map
$$
\Upsilon : \Testar (\RR\times X)\backslash 0 \supset U \to Y
$$
constant in the fibers, which is approximately invariant under the flow of
$H_g$ and which will serve as a useful localizer.  The subset $U$ on which
$\Upsilon$ is defined is a conic neighborhood of $\IC$.  The map $\Upsilon$ is
constructed as follows: consider the data $(x,\Pi)$ of
Theorem~\ref{thm:normalform} as identifying a neighborhood of $Y$ in $X$
with a neighborhood of $x=0$ in the \emph{model cone} $\tX=\RR_+\times Y$, which
we now equip with the \emph{model metric} $g_0.$ Then for any point $q$
near $\IC$ in $\Testar(\RR\times \tX)$, set
\begin{equation}
\Upsilon (q) = \lim_{s\to s_\infty} \pi_Y \exp_q (s (x^2/2)H_{g_0}),
\label{upsilon}\end{equation}
where
$$
s_\infty = h_0(\eta(q))^{-\h} \left((\sgn\theta) \frac\pi 2 - \arctan\theta
\right),\text{ with } \theta= \frac{\xi(q)}{h_0(\eta(q))^\h},
$$
and where $\pi_Y$ is the projection onto the factor $Y$ of $\tX$.  As will
be shown below, the signs are chosen so that $\Upsilon(q)$ is the limit of
the projection on $Y$ of the unique geodesic through $q$, as it heads
toward the \emph{large end} of the model cone $\tX$, \ie\ the point at
infinity from which the geodesic emanated.  To see that $\Upsilon$ is
well-defined and smooth, note that under the flow along $(x^2/2) H_{g_0}$,
$\xi''=2\xi\xi'$, hence
\begin{equation}
  \begin{aligned}
    \xi(s)& =C\tan (Cs+\theta)\,\,\,&  \lambda(s)& =D\sec (Cs + \theta)\\
    x(s)& = E\sec (Cs + \theta)\,\,\,&   t(s)& = -E \tan (Cs + \theta)+F\\
    h_0(\eta)&=G.
  \end{aligned}
\label{ode:soln}\end{equation}
Since $\xi'=\xi^2+h_0(\eta)$, we compute $C=h_0(\eta(q))^\h$ and
$\theta=\arctan \xi(q)/h_0(\eta(q))^\h$, so that
$$
s_\infty = C^{-1} \left( (\sgn\theta) \frac\pi 2-\theta \right)
$$
depends smoothly on $q$ if $U$ is chosen small enough ($\sgn \theta$ is constant
on components of $U$), and as $s\to s_\infty$, $x\to +\infty$ and is
strictly increasing on $s\in [0, \infty)$; simultaneously, $t\to
\pm\infty$.  Thus, since $(y,\eta)$ are undergoing geodesic flow,
$$
\Upsilon (q) = y(\exp_{(y(q),\eta(q))} s_\infty H_Y)
$$
(where $H_Y$ is geodesic flow on $Y$ with metric $h_0$) is manifestly
smooth in $q$.  By the definition as a limit along flow-lines of $H_{g_0}$, we
also have
$$
\Upsilon_* \frac{x^2}2 H_{g_0} = 0.
$$

Now we turn to the perturbed flow $H_g$: Let $V$ denote the rescaled vector field
$(x^2/2)(\lambda^2+\xi^2+h(\eta))^{-1/2} H_g$ on $\Sestar (\RR\times X)$,
and let the flow along $V$ be parametrized by $s$.  Note that under the
bicharacteristic flow of $(x^2/2) H_g$ on $\Testar (\RR\times X)$, as $\abs
t \to \infty$, $s \to C^{-1} (\pm \pi/2-\theta)$, hence $\xi/\lambda \to
\pm C/D$; moreover, $\abs\xi,\ \abs\lambda \to \infty$ in this limit, while
$\eta$ remains bounded.  Thus we have established:
\begin{lemma}\label{lemma:intcurves}
Every maximally extended integral curve of $V$ over $\pa X$ contains in its
closure exactly one point in $\IC$ and one in $\OG$; the former lies over
the point $\Upsilon (p) \in Y$, for any $p$ along the integral curve.
\end{lemma}

\section{Construction of symbols of test operators}

We now write down the symbols of the operators to be used in the commutator
estimates in \S\ref{section:propagation}.  As usual, we work in a product
neighborhood of a boundary component $Y$ as described in
Theorem~\ref{thm:normalform}.  For points $y_1,$ $y_2\in Y$, we let
$d(y_1,y_2)$ denote distance with respect to the metric $h_0.$

Let $\chi \in \CI (\RR)$ vanish for $x<0$, be equal to $1$ for $x>1$, and
be nondecreasing, with smooth square root, such that $\chi'$ also has
smooth square root.  Choose $\psi(x) \in \CIc (\RR)$ to be equal to $1$ at $x=0$,
be supported in $(-1, 1),$ with derivative supported in $(-1,-1/2)\cup
(1/2,1),$ and to be the square of a smooth function; let $(\sgn x) \psi'(x)$
also be the square of a smooth function.  For positive constants $\epsilon
_i$ let $\psi_i(x) = \psi(x/\ep_i)$ and $\chi_i(x)=\chi(x/\ep_i)$.

First we consider test symbols at incoming radial points.

Given $m,l \in \RR$, $\bar y \in Y$, and $\bar x \in \RR_+$, define a
nonnegative symbol on $\Testar X$ by setting
\begin{multline*}
(\ain[m,l,\pm])^2 = \chi(\pm \lambda) \chi(\pm \xi) \chi_1(x-x_0+
vt)\chi_1(-x+x_1- v't) \\ \,\,\,\,\cdot\psi_2((d(\Upsilon,\bar y)^2-\delta x)_+)
\psi_3(h(\eta)^\h/\abs\lambda) \chi_4(t+\ep_4) \psi_5
(p(t,x,y,\lambda,\xi,\eta)/\lambda^2)(\pm\lambda)^m x^l.
\end{multline*}
where $v<1<v'$, $x_0<\bar x<x_1$, $\delta>0$, and we have written
$\Upsilon=\Upsilon(t,x,y,\lambda,\xi,\eta)$ for the map defined by
\eqref{upsilon}.  We will assume that $x_1$ is sufficiently small that the
perturbation term $W=H_g-H_{g_0}$ has the property
$$
(x^2/2) W = h_0(\eta)^\h ( \ocal(x) \pa_\xi+ \ocal(x) \pa_\eta + \ocal(x) \pa_y)
$$
(from \eqref{ham2}) where
\begin{equation}
\text{all } \ocal(x) \text{ terms are bounded by } 10^{-2} \text{ when } x<x_1
\label{xsmall1}
\end{equation}
and furthermore that
\begin{equation}
(x/2) \frac{\pa \log h}{\pa x} <10^{-2} \text{ for } x<x_1.
\label{xsmall2}
\end{equation}
Since $(x^2/2) H_{g_0} d(\Upsilon, \bar y)^2$ vanishes, \eqref{lowerorder}
yields
\begin{equation}
\abs{(x^2/2) H_g d(\Upsilon, \bar y)^2} \leq A x h(\eta)^{\h} \text{ when } x<x_1
\label{def.const}
\end{equation}
for some constant $A$.

Observe that on $\supp \ain[m,l,\pm]$, both
$$
\frac{\abs{\lambda^2-\xi^2-h(\eta)}}{\lambda^2}<\ep_5
$$
and
$$
\frac{h(\eta)^\h}{\abs \lambda}<\ep_3.
$$
Hence
\begin{equation}
\sqrt{1-\ep_5-\ep_3^2}< \abs\xi/\abs\lambda < \sqrt{1+\ep_5}
\label{blahblah}
\end{equation}
so that $\supp \ain[m,l,\pm]$ is localized arbitrarily near $R_{\pm, I}$.

We choose $\delta$ small enough that $\delta x_1<\ep_2/2$, hence $\supp
\psi_2(\cdot) \cap \{x=x_1\} \neq \emptyset$.  We now choose the other
$\ep_i$ sufficiently small and $v,v'$ sufficiently close to $1$, so that
$\ain[m,l,\pm]$ has support in an arbitrarily small neighborhood of the
closure of a bicharacteristic segment
$$
\{t= s,\ x= \bar x -s ,\ y=\bar y,\ \lambda=\pm 1,\ \xi=\pm 1,\
\eta=0; s\in [0,\bar x]\}
$$
passing through the point $q$ with coordinates $x=\bar x$, $y=\bar y$,
$t=0$, $\xi=\pm 1$, $\eta=0$ and hitting $\IC$ at time $\bar x$.  Note
that we may translate the $t$ variable freely without changing any of the
properties of $\ain$.

We now evaluate, term by term,
$$
\frac{x^2}{2} H_g (\ain[m,l,\pm])^2.
$$
The term containing
$$
\frac{x^2}2 H_g \chi_1(x-x_0+vt) = \chi_1'(x-x_0+vt)(\xi-v\lambda+\ocal(x)h_0(\eta)^\h) x
$$
has sign $\pm$ on $\supp \ain[m,l,\pm]$, since \eqref{blahblah} implies
$\abs{\xi} \geq \abs{\lambda}(1-\ep_5-\ep_3)$, hence by \eqref{xsmall1},
$\sgn(\xi-v\lambda)=\sgn \xi$ on $\supp \ain$ provided
$v<1-\ep_5-\ep_3-10^{-2} \ep_3$.

Similarly, the term involving $(x^2/2) H_g \chi_1(-x+x_1-v't)$ can also be
made to have a positive derivative if $v'$ is sufficiently greater than $1$.
Shrinking $\ep_3$ and $\ep_5$ as necessary, these conditions can always be
achieved.

The term containing $\frac{x^2}2 H_g \psi_2((d(\Upsilon,\bar y)^2-\delta x)_+)$
gives
$$
\psi_2'((d(\Upsilon,\bar y)^2-\delta x)_+) \cdot(-\delta \xi x + \ocal(x) h(\eta)^\h)
$$
where, by \eqref{def.const}, the $\ocal (x)$ stands for something bounded
by $A x$.  Hence if
\begin{equation}
\delta\abs\xi> Ah(\eta)^\h,
\label{aaa}
\end{equation}
this term has sign $\pm$; by \eqref{blahblah} and the support property
of $\psi_3(h(\eta)^\h/\abs{\lambda})$, the condition
$$
\delta>\frac{A\ep_3}{\sqrt{1-\ep_5-\ep_3^2}}
$$
suffices to ensure \eqref{aaa}.  This can be achieved by further
shrinking $\ep_3$ as necessary.

The term involving $(1/2) {x^2} H_g \psi_3(h(\eta)^\h/\abs{\lambda})$
can be evaluated using \eqref{ham1}--\eqref{ham2}, which show that $H_g h(\eta)
= \xi x \pa h/\pa x$ and hence
$$
\frac{x^2}2 H_g (h^\h/\lambda) = \frac{\xi h^\h}{\lambda}
\left(\frac{\pa h}{\pa x} \frac {x}{2h} -1\right),
$$
so by \eqref{xsmall2} this term has sign $\pm.$

The term involving $(1/2) x^2 H_g\chi_4(t+\ep_4)$, unlike those discussed
previously, has sign $\mp$.  This term is supported in a region in which we
will assume microlocal regularity: its support is in
$$
\{ \abs{t}<\ep_4,\ x\in (x_0,x_1),\ d(\Upsilon,\bar y)^2<\delta x_1+\ep_2,\
h(\eta)^{1/2}/\abs\lambda<\ep_3 \};
$$
this lies away from $\pa X$ but inside an arbitrarily small neighborhood in
$\Sestar (\RR\times X)$ of an arbitrarily specified point in $R_{\pm,
I}^\ep$.  Call this term $e$, for ``error.''

The term involving $\frac{x^2}2 H_g \psi_5$ is supported in
$p/\lambda^2>\ep_5/2$, hence vanishes identically on the characteristic
variety $\Sigma$.  Call this term $k$. The terms arising from $(1/2) {x^2} H_g
\chi(\pm\lambda)$ and $(1/2){x^2} H_g \chi(\pm\xi)$ are supported in
a compact subset of $\Testar(\RR\times X)$; let $c$ be their sum.

Finally, the factor $(\pm \lambda)^m x^l$ has derivative $(m+l) \xi (\pm
\lambda)^m$, which has sign $\pm\sgn m,$ \ie\ $\pm$ as long as $m+l> 0$;
this term, of course, has the same support as $\ain[m,l,\pm]$ itself;
denote it $\pm (a')^2.$

All the nonnegative resp.\ nonpositive terms described above can be
arranged to be squares of smooth functions resp.\ minus squares of smooth
functions. We organize the information gleaned above as follows. Let
$$
q=(\bar t,\ \bar x,\ \bar y,\ \lambda=\pm 1,\ \xi=\pm 1,\ \eta=0)
$$
be a given point in $R_{\pm, I}$ with $\bar x$ sufficiently small; let $\Omega$
be the closure of the bicharacteristic connecting $q$ to the boundary:
$$
\Omega =\{ (\bar t+s,\ \bar x - s,\ \bar y,\ \lambda=\pm 1,\ \xi=\pm 1,\
\eta=0): s\in [0, \bar x]\}.
$$
Thus we have shown:

\begin{lemma}\label{melwun1.336} Provided $m+l>0$ there exists a symbol
$\ain[m,l,\pm]$ of order $m$ and weight $l$ in $\Testar(\RR\times X)$ such that
\begin{equation}
\frac{x^2}2 H_{\bar g} (\ain[m,l,\pm])^2 = \pm (a')^2\pm\sum_j b_j^2+e+c+ k
\label{incomingconvexity}
\end{equation}
where $a'= a\cdot (\pm (m+l)\xi)^\h,$
$\supp \ain[m,l,\pm]$ is an arbitrarily small neighborhood of $\Omega,$
$\supp e$ is an arbitrarily small neighborhood of $q,$ $\supp c$ is
compact in $\Testar(\RR\times X)$ and $\Sigma\cap\supp(k)=\emptyset.$
\end{lemma}

Next we consider test symbols at outgoing radial points. Let $\chi$, $\psi$
be as above. Given $\bar x,\ \bar y,$ set
\begin{multline*}
  (\aout[m,l,\pm])^2 = \chi(\pm \lambda) \chi(\mp \xi)
  \chi_1(-x+x_1+vt)\chi_1(x-x_0- v't) 
  \\ \,\,\,\,\cdot\psi_2((d(\Upsilon,\bar
  y)^2+\delta x)_+) \psi_3(h(\eta)^\h/\abs\lambda)
  \chi_4(\ep_4-t) \psi_5
  (p(t,x,y,\lambda,\xi,\eta)/\lambda^2)(\pm\lambda)^m x^l
\end{multline*}
where $v<1<v'$ and $x_0<\bar x <x_1$.  As with $\ain$, we can choose
constants $\epsilon _i$ small enough that $\supp \aout[m,l,\pm]$ lies in a
small neighborhood of the closure of a bicharacteristic segment
$$
\Omega = \{ (t= -s,\ x= \bar x -s ,\ y=\bar y,\ \lambda=\pm 1,\ \xi=\mp 1,\
\eta=0); s\in [0,\bar x] \}
$$
passing through the point $q$ with coordinates $x=\bar x$, $y=\bar y$, $t=0$, $\xi=\pm
1$, $\eta=0$ and emanating from $\pa (\RR\times X)$ at time $-\bar x$.
Moreover if $m+l<0$, all terms in $(x^2/2) H_g \aout[m,l,\pm]$ can be arranged
to be $\pm$ squares of smooth functions, with the exception of a compactly
supported term, a term supported away from $\Sigma$, and, most importantly,
the term involving
$$
\frac{x^2}2 H_g (\psi_3(h(\eta)^\h/\abs\lambda)).
$$
Let $e$ denote this error term.  Then
$$
\supp(e)\subset\supp(\aout[m,l,\pm])\cap
\{h(\eta)^\h/\abs\lambda \in [\ep_3/2,\ep_3]\}.
$$
This is a subset of the complement of $R_{\pm, O}$ inside any given
positive conic neighborhood of $\Omega.$

More generally, let
$$
q=(\bar t,\ x=0,\ \bar y,\ \lambda=\pm 1,\ \xi=\mp 1,\ \eta=0)
$$
be any point in $\OG$; let $\Omega$ be the closure of the short
bicharacteristic extending from $q$ to $(t=\bar t+\bar x,\ x=\bar x,\
y=\bar y,\ \lambda=\pm 1,\ \xi=\mp 1,\ \eta=0)$, so
$$
\Omega =\{ (t=\bar t + s,\ x=s,\ \bar y,\ \lambda=\pm 1,\ \xi=\pm 1,\
\eta=0): s\in [0, \bar x]\}.
$$

\begin{lemma}\label{melwun1.337} Provided $m+l<0$ there is a symbol
$\aout[m,l,\pm]$ of order $m$ in $\Testar(\RR\times X)$ such that
\begin{equation}
\frac{x^2}2 H_{\bar g}(\aout[m,l,\pm])^2 =\pm (a')^2\pm\sum_j b_j^2+e+c+ k
\label{outgoingconvexity}
\end{equation}
where $\supp(\aout[m,l,\pm])$ is an arbitrarily small neighborhood of
$\Omega,$ $a'=a(\mp (m+l)\xi)^\h,$ $\supp(e)$ is contained in the
complement of $R_{\pm, O}$ in an arbitrarily small neighborhood of $\Omega,$
$\supp(c)$ is compact and $\Sigma\cap\supp(k)=\emptyset.$
\end{lemma}

\section{Propagation of edge wavefront set}\label{section:propagation}

In this section, we prove a theorem on propagation of singularities for
the edge wavefront set which is a central ingredient in the diffractive
theorem (Theorem~\ref{15.5.2000.7}), the geometric propagation theorem
(Theorem~\ref{15.5.2000.8}) and the proof of the conormal regularity of the
diffracted front in Section~\ref{section:genprop}.

\begin{theorem} \label{melwun1.a319}
For $u\in\eH{-\infty,l}(I\times[0,\ep)_x\times Y),$ a distributional
solution to the wave equation $\Box u=0,$ with $\bar t\in I\subset\bbR$ open,
the following four propagation results hold.
\renewcommand{\theenumi}{\roman{enumi}}
\begin{enumerate}
\item
If $p\in\IC(\bar t,Y),\ m>l+(n-1)/2$ and $\WF^m(u)\cap R^\ep_{\pm, I}(\bar t,
y(p))=\emptyset$ then $p\notin \eWF[m,l'](u)$ for all $l'<l.$
\label{melwun1.352}
\item
For any $m\in\bbR$, $\eWF[m,l](u)\cap\Sestar[\RR\times\pa X]
(\RR\times X)\backslash(\IC(\bar t)\cup \OG (\bar t))$
is a union of maximally extended integral curves of
$V=(x^2/2)(\lambda^2+\xi^2+h_0(\eta))^{-1/2} H_g.$
\label{melwun1.354}
\item
If $U\subset\Sestar[\RR\times \pa
X](\RR\times X)$ is a  neighborhood of $p \in \OG(\bar t,Y)$ and
$\eWF[m,l](u)\cap U\subset\OG$ then $p \notin \eWF[M,l](\Lap_Y^k u)$ for
$k\in\bbN_0$ provided $M\le m-2k$ and $M<l+(n-1)/2.$
\label{melwun1.353}
\item If $R^j u\in \eH{-\infty, l}(\RR\times X),$ with $R$ given by
\eqref{melwun1.175}, and $p\in\OG(\bar t)$ has a neighborhood
$U\subset\Sestar[\RR\times Y](\RR\times X)$ such that $\eWF[m,l](R^{j'}u)\cap
U\subset\OG,$ for $0\le j'\le j,$ then $p \notin   \eWF[M,l](R^j u)$ for
$j\in\bbN_0$ provided $M\le m-j$ and $M<l+(n-1)/2.$
\label{fourthpart}
\end{enumerate}
\end{theorem}

\begin{remark} This theorem correspond to propagation into, within and out
of the boundary. The first two parts, together with the fact that
$\eWF[m,l](u)$ is closed, can be combined with
Lemma~\ref{lemma:intcurves} to conclude that for $p\in \Sestar[\RR \times
\pa X] (\RR\times X)\backslash \OG$, if $\WF^m (u)\cap R_{\pm, I}^\ep
(t(p), \Upsilon(p))=\emptyset$ then $p\notin \eWF[m,l'](u),$ for $l'<l.$

The conclusion of \eqref{melwun1.353} implies, by closedness of the edge
wavefront set and H\"ormander's Theorem, that in fact all of $R_{\pm,
  O}^\ep(\bar t, y(p))$ is absent from $\WF^m(u).$ This third part of the
theorem is trivial, however, when it is applied with $l=s-n/2$ to a
solution in $\dcal_s.$ Therefore this part of the theorem (and the fourth
part likewise) is useless in the absence of a ``division theorem'' yielding
better $x$ decay of $u$ than is given by energy estimates.
\end{remark}

\begin{proof} First consider \eqref{melwun1.352}. By assumption, $u\in
\eH{q,l}(\bbR\times X),$ locally near $\{\bar t\}\times Y$ for some $q\in
\RR$ and $p \in \IC(\bar t,Y).$ We shall prove the following
statement:
\begin{equation}
\begin{gathered}
\text{If } m'>l'+n/2-1,\ u \in \eH{-\infty, l'}(I\times[0,\ep)\times Y),\ p \notin
\eWF[m',l'](u),\\ \text{ and } R_{\pm, I} (\bar t, y(p)) \cap \WF^{m'+1/2}(u)
=\emptyset\text{ then } p \notin \eWF[m'+1/2, l'](u).
\end{gathered}
\label{firstassertion}
\end{equation}

To do so, choose 
\begin{equation}
A_\delta \in \ePs {m',-l'-n/2+1} (\RR \times X) = \Ope \big[\psi_\delta
(\lambda,\xi,\eta)\ain[m',-l'-n/2+1,\pm]\big],
\label{melwun1.358}\end{equation}
where $\ain$ is constructed in Lemma~\ref{melwun1.336} and
$$
\psi_\delta(\lambda, \xi,\eta)= \psi((\lambda^2+\xi^2+h(\eta))\delta)
$$
with $\psi(x)$ smooth, equal to $1$ for $x<1/2$ and $0$ for $x>1.$
Choosing the supports sufficiently small, we see that the error term $e$ in
\eqref{incomingconvexity} will have support in the complement of
$\WF^{m'+1/2}(u).$ Note that $H_g \psi_\delta$ is supported in
$\{1/(2\delta)\leq\lambda^2+\xi^2+h(\eta)\leq 1/\delta\}$.  

Thus if $A'_\delta$ has symbol $\psi_\delta a'$ with $a'$ as in
\eqref{incomingconvexity},
\begin{equation}
[\Box, A_\delta^*A_\delta] =\pm (A'_\delta)^* (A'_\delta) \pm \sum_j B_{\delta,
j}^* B_{\delta,j} + E_\delta + K_\delta + R_\delta + S_\delta
\label{comm}
\end{equation}
where 
\begin{itemize}
\item
$\eWFprime (K_\delta) \cap \Sigma=\emptyset$ 
\item
$A'_\delta \in \ePs{m'+\h,-l'-\frac n2}(\RR \times X)$ with $A'_\delta \to
A'\equiv A'_0$ in $\ePs {m'+\h+\ep, -l'-\frac n2}(\RR \times X)$ for all
$\ep>0$ as $\delta \to 0$.
\item
$E_\delta \in \ePs{2m'+1, -2l'-n}(\RR \times X)$ is bounded in $\delta$ with
$\eWFprime(E_\delta)$ uniformly bounded away from $\pa X$ and contained in
$(\WF^{m'+1/2} u)^\complement,$ 
\item
$R_\delta$ bounded in $\ePs{2m', -2l'-n}(\RR \times X)$  
\item
$S_\delta$ is bounded in $\ePs{2m'+1, -2l'-n}(\RR \times X)$ and, as
$\delta\to 0,$ converges to $0$ in $\ePs{2m'+1+\ep, -2l'-n}(\RR \times X)$
for all $\ep>0$
(this is the term whose symbol involves $H_g \psi_\delta$).
\end{itemize}

Applying \eqref{comm} to $u$ and pairing with $u$ \emph{with respect to the
inner product on $L^2_g$} yields
\begin{equation}
\norm{A'_\delta u}_g^2 -\ang{S_\delta u, u}_g \leq \abs{\ang{E_\delta u,u}_g}+\abs{\ang{K_\delta
u,u}_g}+\abs{\ang{R_\delta u,u}_g};
\label{pairing}
\end{equation}
the integration by parts is justified since $u \in \eH{-\infty,
l'}(\RR\times X)$ and $\Box A_\delta^* A_\delta u \in \eH{\infty,
-l'-n}(\RR\times X)$.  All terms on the right-hand side are bounded
uniformly as $\delta\to 0$.  A weak convergence argument now shows that
$\norm{A'_0 u}_g<\infty$, hence $\eWF[m'+\h, l'](u)=\emptyset$ is disjoint
from the elliptic set of $A'_0;$ the shift by $n /2$ in the $x$ weight here
comes from the difference between $L^2_g$ and $L^2_{\bo}.$ This proves
\eqref{firstassertion}.

If $q> l+(n-2)/2,$ then iterative application of \eqref{firstassertion}
proves that $p\notin \eWF[m,l](u)$ directly. If $q\leq l+(n-2)/2$, however,
a further argument is needed.

Supposing $l_0 = \sup \{r;p \notin \eWF[m, r](u)\}<l,$ we wish to arrive at
a contradiction. We will employ an interpolation argument illustrated in
Figure~\ref{figure:interpolation}.  We have already shown that $l_0\geq q-(n-2)/2$ and by
hypothesis, $u\in\eH{q,l}(\bbR\times X)$ (at least locally), so 
Proposition~\ref{prop:wfinterp} shows that 
\begin{equation}
p\notin \eWF[\theta q + (1-\theta) m, \theta l+(1-\theta) l_0](u)
\ \forall\ \theta \in [0,1].
\label{melwun1.357}\end{equation}
In particular, since $l<m-(n-1)/2$ (by the hypothesis of the theorem)
and $q\leq l+n/2-1$,
$$
\theta'=(m-l_0-n/2+1)/((m-l_0-n/2+1)+(l-q+n/2-1)) \in [0,1].
$$
If $m'=\theta' q + (1-\theta') m$ and $l'=\theta' l+(1-\theta')l_0,$ then
$m'=l'+n/2-1$, and by \eqref{melwun1.357}, $p\notin \eWF[m', l'-\ep](u)$ for
all $\ep>0.$ Hence applying \eqref{firstassertion} iteratively shows that $p
\notin \eWF[m,l'-\ep](u)$, with $l'-\ep>l_0$, which is the desired
contradiction. Thus \eqref{melwun1.352} is proved.

\begin{figure}[ht!]
\setlength{\unitlength}{0.00035in}
\newcommand{\radius}{200}
\begingroup\makeatletter\ifx\SetFigFont\undefined%
\gdef\SetFigFont#1#2#3#4#5{%
  \reset@font\fontsize{#1}{#2pt}%
  \fontfamily{#3}\fontseries{#4}\fontshape{#5}%
  \selectfont}%
\fi\endgroup%
{\renewcommand{\dashlinestretch}{30}
\begin{picture}(6952,6660)(0,-10)
\put(6683,4212){\blacken\ellipse{\radius}{\radius}}
\put(6683,4212){\ellipse{\radius}{\radius}}
\put(83,5412){\blacken\ellipse{\radius}{\radius}}
\put(83,5412){\ellipse{\radius}{\radius}}
\put(6683,5412){\blacken\ellipse{\radius}{\radius}}
\put(6683,5412){\ellipse{\radius}{\radius}}
\put(5183,4512){\blacken\ellipse{\radius}{\radius}}
\put(5183,4512){\ellipse{\radius}{\radius}}
\put(6683,4512){\blacken\ellipse{\radius}{\radius}}
\put(6683,4512){\ellipse{\radius}{\radius}}
\path(83,6612)(83,12)(6683,12)
\path(683,12)(6683,6012)
\dashline{60.000}(83,5412)(6683,4212)
\dottedline{45}(5183,4512)(6683,4512)
\put(233,5562){$(q,l)$}
\put(4558,4737){$(m',l'-\ep)$}
\put(6383,5000){$(m,l)$}
\put(6383,3837){$(m,l_0)$}
\put(3083,2262){$l''=m''-n/2+1$}
\put(6308,162){$m''$}
\put(233,6537){$l''$}
\end{picture}
}
\caption{The interpolation argument in part one of
Theorem~\ref{melwun1.a319}: we begin with global regularity of order
$(q,l)$ and microlocal regularity of order $(m,l_0)$.  Interpolation gives
microlocal regularity of order $(m',l')$ and iterative application of
\eqref{firstassertion} is used to move along the horizontal line and obtain
microlocal regularity of order $(m,l'-\ep)$ with
$l'-\ep>l_0$.\label{figure:interpolation}}
\end{figure}
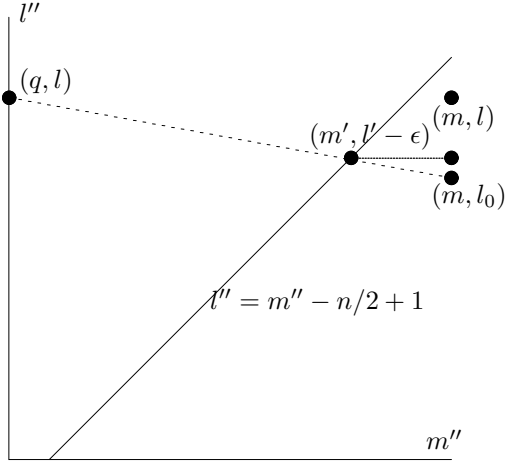

To prove \eqref{melwun1.354}, \ie\ to show that regularity propagates
across $\Sestar[\pa X] (\RR\times X)$ up to (but not including) points in
$\OG$, we appeal to the standard proof of H\"ormander's theorem on
propagation of singularities for operators of real principal type by use of
positive commutator estimates. This applies microlocally near all
characteristic points where the rescaled Hamilton vector field
$V=(x^2/2)(\lambda^2+\xi^2+\abs{\eta}^2)^{-1/2} H_g$ is non-zero, hence away
from $\IC \cup \OG.$ See \cite{Melrose43} for an analogous discussion in
the context of the scattering calculus; note that in that case the
propagation result may be reduced to H\"ormander's theorem whereas in this
case it is an analogue of it.

Now consider \eqref{melwun1.353}, first for $k=0.$ In this case we may
simply suppose that $M=m<l+(n-1)/2.$ We proceed much as in
the proof of \eqref{melwun1.352} above. The result follows by
iterative application of the following assertion:
\begin{multline}
\text{If } m'<l+n/2-1,\ u \in \eH{-\infty, l}(I\times X)\ \text{and } p \notin
\eWF[m',l](u)\text{ then}\\
p\notin\overline
{\Sestar[\RR\times \pa X](\RR\times X)\cap\eWF[m'+1/2,l](u)\setminus\OG}
\Longrightarrow p\notin \eWF[m'+1/2, l](u).
\label{secondassertion}
\end{multline}
Note that the final hypothesis on $p$ here is equivalent to the existence
of a neighborhood $U$ of $p$ such that $U\cap\eWF[m'+1/2,l](u)\subset\OG.$

To prove \eqref{secondassertion}, choose $A_\delta \in \ePs{m',-l-n/2+1}
(\RR \times X)$ as in \eqref{melwun1.358} with $\ain[m',-l-n/2+1,\pm]$
replaced by $\aout[m',-l-n/2+1,\pm],$ where $\aout$ is constructed in
Lemma~\ref{melwun1.337} and supports are kept small, corresponding to the
implicit neighborhood, $U,$ in \eqref{secondassertion}.  If 
$A'_\delta$ has symbol $\psi_\delta a'$ with $a'$ as in 
\eqref{outgoingconvexity}, then 
\begin{equation*}
[\Box, A_\delta^*A_\delta] =\pm (A'_\delta)^* (A'_\delta) \pm \sum_j B_{\delta,
j}^2 + E_\delta + K_\delta + R_\delta+ S_\delta
\label{melwun1.361}\end{equation*}
where
\begin{itemize}
\item
$\eWFprime(K_\delta)\cap \Sigma=\emptyset$
\item
$A'_\delta \in \ePs {m'+\h,-2l-n}(\bbR\times X)$ and $A'_\delta \to
A'\equiv A'_0$ in $\ePs {m'+\h+\ep,-l-\frac n2}(\bbR\times X)$ for all
$\ep>0$ as $\delta \to 0$.
\item
$E_\delta \in \ePs{2m'+1, -2l-n}(\bbR\times X)$ and
$\eWFprime(E_\delta)\subset U \backslash \OG(\bar t)$, uniformly in $\delta$.
\item
$R_\delta$ is uniformly bounded in $\ePs{2m', -2l'-n}(\bbR\times X)$
\item
$S_\delta$ is bounded in $\ePs {2m'+1, -2l-n}(\bbR\times X)$ and, as
$\delta\to 0$, converges to $0$ in $\ePs{2m'+1+\ep, -2l-n}(\bbR\times X)$
for all $\ep>0$.
\end{itemize}
Hence 
\begin{equation}
\norm{A'_\delta u}_g^2 -\ang{S_\delta u, u}_g \leq \abs{\ang{E_\delta u,u}_g}+\abs{\ang{K_\delta
u,u}_g}+\abs{\ang{R_\delta u,u}_g},
\label{melwun1.362}\end{equation}
so $\eWF[m'+\h, l](u)$ is disjoint from the elliptic set of $A'.$ This proves
\eqref{secondassertion}.

Finally consider \eqref{melwun1.353}.  We work by induction, assuming that
the result is known for all nonnegative integers smaller than a given $k.$
Note that $\Lap_Y$ is an edge (pseudo)differential operator, so the
conclusion of \eqref{melwun1.353} for positive $k$ is only non-trivial if
$m\ge l+(n-1)/2.$ In particular $p\notin\eWF[m',l](\Lap_Y^k)$ for $m'$
sufficiently negative. We therefore prove the following analogue of
\eqref{secondassertion} with $u$ replaced by $\Lap_Y^ku:$
\begin{multline}
\text{If } m'<l+n/2-1,\ u \in \eH{-\infty, l}(I\times X)\ \text{and } p \notin
\eWF[m,l](u)\text{ then}\\
p\notin\overline
{\Sestar[\RR\times \pa X](\RR\times X)\cap\eWF[m'+1/2,l](\Lap_Y^ku)\setminus\OG}
\Longrightarrow p\notin \eWF[m'+1/2, l](\Lap_Y^ku),
\label{melwun1.364}\end{multline}
proceeding much as before.

To prove \eqref{melwun1.364}, note that Lemma~\ref{lemma:commutators} gives
a distributional equation for $u_{k}=\Lap_Y^ku$ of the form
\begin{equation}
\begin{aligned}
\Box u_{k}+\left(Q_kD_x+\frac1xP_{k}\right)u=0,\Mwhere
&P_{k}\in\CI([0,\ep);\Diff{2k+1}(Y)),\\
&Q_k \in \CI([0,\ep); \Diff{2k-1}(Y)).
\label{melwun1.360}
\end{aligned}
\end{equation}
Now applying the test operator $A_\delta$ and pairing with $u_k$
gives an estimate similar to \eqref{melwun1.362} with an extra term: 
\begin{multline}
\norm{A'_\delta u_k}_g^2 - \ang{S_\delta u_k, u_k}_g \leq \abs{\ang{E_\delta u_k,u_k}_g}+\abs{\ang{K_\delta
u_k,u_k}_g}\\
+\abs{\ang{R_\delta u_k,u_k}_g}
+\abs{\ang{A_\delta \Lap_Y^ku,A_\delta B_ku}_g}.
\label{additionalterm}
\end{multline}

The first term on the left can be reorganized, modulo terms uniformly controlled
by the inductive hypothesis, to bound a positive multiple of $\|A'_\delta
u\|_{H^{2k}(Y)},$ the tangential Sobolev norm of $A'_\delta u.$ We now show
that the last term on the right can be estimated by a small multiple of
this norm, modulo the inductive bounds.  Note that $A_\delta =G_\delta
A'_\delta+C_\delta+D_\delta$ where $G_\delta\in\ePs{-1/2,1}(\bbR\times X)$ is
uniformly bounded as $\delta\downarrow0,$ $C_\delta$ is lower order and
$D_\delta$ is supported in the region of known regularity. Modulo terms
bounded by the inductive hypothesis or by the hypothesis of 
\eqref{melwun1.364}, we can thus write the last term as
$$
\ang{\Lap_Y^k A_\delta' u, B_k G_\delta^* G_\delta A_\delta' u}_g.
$$
We may rewrite $B_k$ as a sum of terms $x^{-1} C_k S_i$ with $C_k \in
\Diff{2k}(Y)$ and where $S_i$ are smooth b-vector fields.  Hence modulo
controllable terms the inner product above is estimated by a sum of terms
of the form
$$
\norm{\Lap_Y^k A_\delta' u}\norm{x^{-1} C_k S_i G_\delta^* G_\delta A_\delta' u}_g.
$$
The latter norm is (again modulo known terms) controlled by $\ep
\norm{A'_\delta u}_{H^{2k}(Y)}$, with the $\ep$ coming from the small
support in $x.$ Thus the last term in \eqref{additionalterm} can be
absorbed in the left, giving the inductive estimate and proving
\eqref{melwun1.364}.

To prove \eqref{fourthpart}, we appeal once again to
Lemma~\ref{lemma:commutators} and proceed as with \eqref{melwun1.353}.
\end{proof}

\begin{theorem}\label{thm:propagationthrough} Let $u(t)$ be the solution of
the Cauchy problem \eqref{15.5.2000.1}--\eqref{15.5.2000.2} and suppose that
$u\in\eH{-\infty,l}(I\times X)$ for some open $I\ni \bar t$ and some
$l\in\bbR$ then 
\begin{multline*}
p \in R_{\pm,O}^\ep(\bar t),\ \Gamma^\ep(p)\cap\WF^{m}(u)=\emptyset\Mfor
m<l+(n-1)/2\Mand\ep>0\\
\Longrightarrow p\notin\WF^{m-\delta}(u)\ \forall\ \delta>0.
\label{melwun1.356}\end{multline*}
\end{theorem}

\begin{proof} By \eqref{melwun1.352} and \eqref{melwun1.354} of
Theorem~\ref{melwun1.a319}, there is no edge wavefront set of order $(m,
m-(n-1)/2-\delta)$ along all bicharacteristics in $\Sestar[\RR \times \pa X]
(\RR\times X)$ terminating at $t=\bar t,$ $y=y(p).$ Hence by 
\eqref{melwun1.353} of Theorem~\ref{melwun1.a319}, closedness of edge
wavefront set, and H\"ormander's theorem, $p\notin \WF^{m-2\delta}(u).$
\end{proof}

\section{The rescaled FBI transformation}
\label{section:FBI}
To analyze the boundary behavior of solutions to the wave equation, we
employ a rescaled version of the one-dimensional FBI
(``Fourier-Bros-Iagolnitzer'') transformation in the time variable. The FBI
transform is a variant of the Bargmann transform which was employed by Bros
and Iagolnitzer \cite{Iagolnitzer1} to study microlocal regularity in the
analytic setting. Its properties were further elaborated by Sj\"ostrand
\cite{Sjostrand1}. For our purposes, the FBI transform could be dispensed
with in favor of a composition of localization in time and Fourier
transform. Use of the FBI transform, however, seems more likely to admit
generalization. The \ci\ approach to the FBI used here follows the spirit
of \cite{Wunsch-Zworski2}, which deals with compact manifolds.  A partial
FBI transform is also used in a related analysis of the wave equation by
G\'erard and Lebeau \cite{MR93f:35130}.

Consider the complex phase function $\phi (t,\tau,t') = i (t-t')^2
\ang{\tau}/2 +(t-t')\tau.$ The associated FBI transform applied just to the
time variable is
$$
T u(t,\tau,x,y) = \int_X e^{i \phi(t,\tau,t')} a(t,t',\tau) u(t',x,y)\,
dt'
$$
where $a$ is a polyhomogeneous symbol of order $1/4$ with proper support in
$t,t'.$ It is an elementary consequence of the stationary phase lemma,
demonstrated in \cite{Wunsch-Zworski2}, that $T^* T$ is a
pseudodifferential operator on $\bbR$ of order $0$ and hence is
$L^2$-bounded. It is elliptic at $(\bar t,\pm\infty)\in S^*\bbR$ if $a$ is
elliptic at $t=t'=\bar t,$ $\tau=\pm\infty.$ Thus the ellipticity of $a$ at
$(\bar t,\infty)$ implies that there exists $G \in \Psi^0(\RR)$ with $(\bar 
t,\infty)\notin\WF'(G)$ such that
\begin{equation}
T^*T=\Id+G,\Mso
Tu=0\Min\tau>1\Mnear\{t=\bar t\}\Longrightarrow u=-Gu\Mnear\{t=\bar t\}.
\label{melwun1.308}\end{equation}
The basic intertwining property of $T$ corresponds to the boundedness of
\begin{equation}
T D_t -\tau T: L^2(\bbR)\longrightarrow \ang{\tau}^{\h} L^2(\bbR^2).
\label{melwun1.380}\end{equation}

We shall choose $a$ to have support in $\tau>1,$ the region
$\tau<-1$ being handled by reflection in $t.$ Fixing a product
decomposition near the boundary and inserting a cutoff $\chi(x)$ with
$\chi(x)=1$ for $x$ sufficiently small and with support in the product
neighborhood of a boundary component $Y,$ we scale the $x$-variable and
define
\begin{equation}
Su(t,\tau,\tx,y)=\int\chi(\frac{\tx}{\tau}) e^{i \phi(t,\tau,t')} a(t,t',\tau)
u(t',\frac{\tx}{\tau},y)\, dt'.
\label{S-new-def}\end{equation}
If we let $\tX$ denote the model cone $[0,\infty) \times Y$ then 
\begin{equation}
S:\CmI(\RR\times X)\longrightarrow\CmI(\RR\times[1,\infty)\times \tX).
\label{melwun1.159}\end{equation}
In fact the assumed properness of the support in $t,t'$ ensures that $Su$
is \ci\ in $\tau$ and $t.$ However, it is the growth in $\tau,$ uniformly
in $\tx,$ that will interest us here. In view of the assumption on the
support of $a$ is we will use a decomposition of the type described in
\eqref{melwun1.388} of $u\in\CmI(\bbR\times X)$ under convolution in
$t$ and examine $u_+.$

Taking into account the scaling, the boundedness of $T$ implies that $S$ is
bounded on the scale-invariant spaces in the $x$ variable
\begin{equation}
S:\bLt(\RR \times X)\longrightarrow \bLt(\RR\times[1,\infty)\times\tX)
\label{melwun1.160}\end{equation}
where the measure on the left is $dx\, dt\, dy/x$ and on the right
$dt\,d\tau\,d\tx \, dy/\tx$ and we assume that $a$ is a function of $t-t'$
only to get the global estimate. The scale-invariance of the measue shows
that 
\begin{equation}
S^*S=\chi (x)G
\label{melwun1.389}\end{equation}
where $G$ is the pseudodifferential operator on $\bbR$ in \eqref{melwun1.308}.

The simplest intertwining properties of $S$ follow directly from
\eqref{S-new-def}:
\begin{equation}
\begin{gathered}
S\circ Qu=Q\circ Su,\ Q\in\Diff*(Y)\Mand\\
S(xu)=\frac{\tx}{\tau}Su.
\end{gathered}
\label{melwun1.200}\end{equation}
It follows from the second of these that if $f\in\dCI(X)$ then composition
with $f$ as a multiplication operator gives
\begin{equation}
S\circ f:\bLt(\RR \times X) \longrightarrow \bigcap\limits_k
\left(\frac{\tx}{\tau}\right)^k\bLt(\RR\times[1,\infty)\times\tX).
\label{melwun1.203}\end{equation}
Note that $\tx/\tau$ is bounded on the support of $Su$ and $\tau>1$ by
assumption. Since we will only use the FBI transform to examine behavior
in Taylor series at the boundary, we shall denote the space of operators with the
property \eqref{melwun1.203} as $E^\infty$ and consider these as error
terms. The presence of the localizing function $\chi$ in the definition of
$S$ means for instance that
\begin{equation}
S \circ xD_x=\tx D_{\tx}\circ S+S'\circ f',\
f'\in\dCI(X),\ S'\circ f'\in E^\infty.
\label{melwun1.204}\end{equation}
Here, $S'$ is an operator of the same form as $S$ with a different cutoff $\chi.$

Scaling \eqref{melwun1.380} and using \eqref{melwun1.200} we conclude that
\begin{equation}
S \circ (x D_t)-\tx S:\bLt(\bbR\times X)\longrightarrow \tau^{-\h}\ang{\tx}
\bLt(\bbR\times[1,\infty)\times\tX).
\label{melwun1.201}\end{equation}

The scaling of $x$ makes $\tilde x$ a global variable on
$\RR_+$, and we are primarily concerned with behavior as
$\tx\to \infty.$ The structure which arises at infinity here
corresponds to the $\RR_+$-homogeneous metric in \eqref{melwun1.163}. This
is a conic metric on $\tX=[0,\infty)\times Y,$ but its uniform behavior
near infinity is rather different from its behavior near $\tx=0.$ Consider
the inversion
\begin{equation}
\tX^\circ=(0,\infty)\times Y\ni(\tx,y)\longmapsto (1/\tx,y)\in
(0,\infty)\times Y=\tX^\circ
\label{melwun1.342}\end{equation}
which induces an isomorphism of $\CmI(\tX).$ Under this transformation, the
metric becomes a scattering metric in the sense of \cite{Melrose43} at the
boundary $w=1/\tx=0,$ which is to say it is a particular type of
asymptotically locally Euclidean metric.

In \cite{Melrose43} the associated
compactly-supported Sobolev spaces $H^m_{\scat,c}(\tX)$ are defined for any
manifold with boundary. It is natural then to introduce Sobolev spaces
which are of ``b-type'' near $\tx=0$ and of ``scattering type'' near
$\tx=\infty,$ based however on the $\RR_+$-invariant measure $d\tx \,
dy/\tx$ (note that this convention differs from the weight used in
\cite{Melrose43}).  Thus if we choose $\phi\in\CIc([0,\infty)),$ with
$\phi(\tx)=1$ near $\tx=0$ then we may define
\begin{multline}
\bscH m(\tX)=\{u\in\CmI(\tX);\phi(\tx)u\in H^m_{\boc}(\tX),\\
(1-\phi)(1/w)u(1/w,y)\in H^m_{\scatc}(\tX)\},\ m\in\bbR.
\label{melwun1.341}\end{multline}
Since both the ``b'' and the ``sc'' Sobolev space reduce to the standard
Sobolev spaces in the interior, this is independent of the choice of $\phi.$ 
We shall also employ weighted versions of these spaces which we will generally
write in terms of the weights $\tx$ and $\ang\tx.$ The former is
a boundary defining function near $0$ but is also, near infinity, of the
form of $1/w$ where $w=0$ is a defining function for inverted infinity. On
the other hand $1/\ang\tx$ is just a defining function for inverted
infinity. Thus we consider the weighted spaces
\begin{equation}
\tx^l\ang\tx^k\bscH m(\tX),\ l,k,m\in\bbR.
\label{melwun1.343}\end{equation}
Directly from the definition, these scales of Hilbertable spaces satisfy complex
interpolation. 

The spaces which arise here correspond to functions of $\tau$ and $t$ with
values in these weighted Sobolev spaces. We are primarily interested in global
behavior in $\tau$ but local behavior in $t.$ We therefore define 
\begin{multline}
\Lloc(\RR_t \times [1,\infty]; \bscH {m}(\tX))=
\big\{u\in\CmI(\bbR\times\bbR\times\tX);u=0\Min\tau<1\Mand\\
\phi (t)u\in L^2(\bbR\times\bbR;\bscH {m}(\tX))\
\forall\ \phi \in\CIc(\bbR)\big\}.
\label{melwun1.381}\end{multline}
More generally weights will be written out as in \eqref{melwun1.343}. Note
that the closed bracket at $\tau=\infty$ in \eqref{melwun1.381} is intended
to indicate that these spaces are indeed global in $\tau.$ 

\begin{lemma}\label{melwun1.205} For any $m,$ $l\in\bbR$
\begin{equation}
S:x^l\eH m (\RR\times X)
\longrightarrow
\bigcap_{\alpha\in [0,m]} \tau^{-l} \tx^{l} \ang{\tx}^{-\alpha}
\Lloc(\RR_t \times [1,\infty]; \bscH {m-\alpha}(\tX)).
\label{melwun1.202}\end{equation}
Conversely provided $\chi '\in\CI(X),$ $\chi '\chi =\chi '$ and $m\ge0,$
\begin{multline}
\Mif u\in\bLt(\bbR\times X)\Mthen\\ Su\in
\bigcap_{\alpha\in [0,m]} \tau^{-l} \tx^{l} \ang{\tx}^{-\alpha}
\Lloc(\RR_t \times [1,\infty]; \bscH {m-\alpha}(\tX))\Mnear\{t=\bar t\}
\Longrightarrow\\
\chi 'u_+=u_1+u_2,\ u_1\in x^l\eH m
(\RR\times X),\ u_2\in\CI(\bbR;\bLt(X))\Mnear\{t=\bar t\}.
\label{melwun1.207}\end{multline}
\end{lemma}

\begin{proof} 
For positive integral $m,$ $u\in x^l\eH m (\RR\times X)$ with support near
the boundary if $(xD_x)^k(xD_t)^pQu\in x^l\bLt(\bbR\times X)$ for all
$Q\in\Diff*(Y)$ and all $k+p+\ord(Q)\le m.$ The continuity estimates in
\eqref{melwun1.202} then follow from \eqref{melwun1.160},
\eqref{melwun1.200}, \eqref{melwun1.203}, \eqref{melwun1.204} and
\eqref{melwun1.201}. Similarly, for negative integral $m,$ $u\in\eH m
(\RR\times X)$ may be written as a finite superposition of edge operators
of order $-m$ applied to elements of $\bLt(\bbR\times X)$ and continuity
follows similarly. Complex interpolation on both sides then gives the
general case of \eqref{melwun1.202}.

To see the partial converse, \eqref{melwun1.207}, first replace $u$ by
$u_+$ as in \eqref{melwun1.388}. Since the amplitude of $S$ is assumed to
be supported in $\tau>1$ this simply changes $S$ by a similar operator with
rapidly decreasing amplitude. Thus we may suppose that the condition holds
for $u_+.$ The invertibility of $S$ in \eqref{melwun1.389} shows that
$S^*S\equiv\chi$ modulo a term arising from $G;$ applying the arguments above
to $S^*$ therefore gives \eqref{melwun1.207}.
\end{proof}

It follows from \eqref{melwun1.202} that if $\phi\in\CIc(\bbR)$ is $1$ near
$0$
\begin{multline}
(1-\phi(\tx))S\circ B:
x^lH^m_{\eo}(\bbR\times X)\longrightarrow\\
\tau^{-l}\ang\tx^{-M} \bscH M(\bbR\times[1,\infty)\times\tX)\ \forall\ M,\
  B\in\Psi^{-\infty}_{\eo,c}(\bbR\times X).
\label{melwun1.209}\end{multline}
Thus $Su$ is rapidly decreasing in the sense of Schwartz near infinity.
Thus the boundary part of the scattering wavefront set, in
the sense of
\cite{Melrose43}, of $Su$ is determined by the edge wavefront set, over the
boundary, of $u.$ In fact we need to use a notion of scattering wavefront
set in a uniform sense in $t$ and $\tau.$ Such uniform versions of the
scattering wavefront set are also discussed in \cite{Melrose43}. As there,
we are only really interested in the scattering wavefront set near $\tx=\infty.$

On a compact manifold with boundary, $X,$ the scattering wavefront set is a
subset of the boundary of the radial compactification of the intrinsic scattering
cotangent bundle. This bundle reduces to the usual cotangent bundle in the
interior but has basis at any boundary point $p$ the differentials $dx/x^2$ and
$dy_j/x,$ where where $x$ is a local boundary defining function and $y=0$
at $p,$ are boundary coordinates. The part of the 
scattering wavefront set over the boundary can be defined explicitly as
follows. If $(\bar\zeta,\bar\eta)\in \Tscstar_pX$ is the coordinate
representation of a general point in the fiber over $p,$
$\bar\zeta\frac{dx}{x^2}+\bar\eta\cdot\frac{dy}x\in\Tscstar_pX$ then
$(\bar\zeta,\bar\eta)\notin\scWFnoargs(u)$ if for some $\psi \in\CI(X)$
supported in the coordinate patch and with $\psi(p)\not=0,$  
\begin{equation}
\int e^{-i\frac\zeta x+i\frac{\eta \cdot y}x}\psi(x,y)u(x,y)\,\frac{dx}x\,dy
\Mis\CI\Mnear (\bar\zeta ,\bar\eta ).
\label{melwun1.383}\end{equation}
Note that changing to the singular variables $1/x,y/x$ identifies the
complement of the boundary in a neighborhood of $p$ with an open set in
$\bbR^n.$ Then the formal integral in \eqref{melwun1.383} may be identified
with the Fourier transform. Since $u$ is an extendible (hence Schwartz)
distribution, \eqref{melwun1.383} is well defined; this condition
is also independent of the choice of coordinates.

The \emph{uniform} version of the scattering wavefront we need is obtained
by strengthening \eqref{melwun1.383}.

\begin{definition}\label{melwun1.385} If $v\in
\tau^{-l}\Lloc(\bbR\times[1,\infty];\ang{\tx}^{p}\bscH {m}(\tX))$ for
some $m,$ $l$ and $p$ then
\begin{equation}
(\bar t,\infty;\bar y,\bar\zeta,\bar\eta)\notin
\scbWF(v;\tau^{-l}\Lloc(\bbR\times[1,\infty];\ang{\tx}^{q}\bscH {*}(\tX)))
\label{melwun1.386}\end{equation}
if $v=v_1+v_2$ where, for some $m,$ $v_1\in
\tau^{-l}\Lloc(\bbR\times[1,\infty];\ang{\tx}^{q}\bscH {m}(\tX))$ and 
\begin{equation}
\int e^{-i\zeta\tx+i{\eta \cdot y}\tx}\psi(\frac1\tx,y)u(\tx,y)
\,\frac{d\tx}\tx\,dy\in\tau^{-l}\Lloc(\bbR\times[1,\infty];\CI(\bbR^n))\Mnear
(\bar t,\bar\zeta,\bar\eta ),
\label{melwun1.384}\end{equation}
where $\psi(x,y)\in\CIc([0,\infty)\times Y)$ is equal to $1$ near $(0,\bar y).$
\end{definition}
\noindent
Note that the regularity, $m,$ here is irrelevant since it is the boundary
part of the scattering wavefront set we are examining.

We will show that the wavefront relation of $S$ is related to
\begin{equation}
\cR=\{(\xi',y,\mu,t;\xi,y',\eta,t',\lambda);\lambda>0,\ t=t',\ y=y',\
(\xi',\mu)=\lambda^{-1}(\xi,-\eta)\}.
\label{melwun1.217}\end{equation}

\begin{proposition}\label{melwun1.211} If $u\in x^lH^{-\infty}_{\eo,c}(\bbR\times
X)$ has $\eWF[\infty,l](u)\subset\{\lambda>\}$ then
\begin{multline}
\scbWF(Su;\tau^{-l}\Lloc(\bbR\times [1,\infty];\ang{\tx}^{l-m}\bscH{*}(\tX)))\\
\subset\cR\circ
\left(\eWF[m,l](u)\cap\Sestar[\bbR\times \pa X](\bbR\times X)\right),
\label{melwun1.212}\end{multline}
in the sense of Definition~\ref{melwun1.385}. 
\end{proposition}

\begin{proof} If $u\in x^lH^m_{\eo}(\bbR\times X)$ then, by
Lemma~\ref{melwun1.205} the scattering wavefront set of $Su$ relative to
$\tau^{-l}\Lloc(\bbR\times [1,\infty];\ang{\tx}^{l-m}\bscH{*}(\tX))$ is
empty, since $Su$ lies in the space $\tau^{-l}\Lloc(\bbR\times
[1,\infty];\ang{\tx}^{l-m}\bscH{0}(\tX)).$ Thus it suffices to prove
\eqref{melwun1.212} for $m=\infty.$

It is therefore enough to consider $u\in x^lH^m_{\eo,c}(\bbR\times X),$ for
some $m,$ having support in the product neighborhood of the boundary and
having scattering wavefront set contained in $\{\lambda>0\}.$ Then, using a
partition of unity in the edge calculus, we may decompose
\begin{equation*}
u=u'+\sum\limits_{j}B_ju
\end{equation*}
where $u'\in x^lH^\infty_{\eo,c}$ and the $B_j$ have small wavefront sets,
in $\lambda>0.$ Again by Lemma~\ref{melwun1.205} the term $u'$
produces a term which is Schwartz near $\tx=\infty$ as a function with
values in the weighted $L^2$ space $\tau^{-l}\Lloc(\bbR\times[1,\infty)).$
This corresponds to the absence of any scattering wavefront set. So we may
replace $u$ by $Bu,$ with an edge pseudodifferential operator,
$B\in\Psi_{\eo}^0(\bbR\times X),$ with essential support concentrated near
some boundary point $(\bar t,0,\bar
y,\bar\lambda,\bar\eta,\bar\xi).$ We may also
suppose that $B$ has its support near the boundary component $Y.$ Since the edge
wavefront set is a conic notion this means that we may take
$\bar\lambda=1.$ In local coordinates
\begin{multline}
Bu(t,x,y)=(2\pi)^{-1}\int e^{i((s-1)\xi+(y-y')\eta+T\lambda)}\\
b(t,x,y,s,y',T,\xi,\eta,\lambda)
u(xs,y',t-xT)x^2\, ds\, dy'\, dT\, d\xi\, d\eta\, d\lambda,
\label{melwun1.208}\end{multline}
where, at the expense of another error of order $-\infty$ in the edge
calculus, we may assume that
\begin{equation}
\begin{aligned}
\text{the amplitude } b &\text{ has small conic support (as } \lambda\to\infty
\text{) around }\\
t=\bar t,\ x=0,\ &y=\bar y,\ s=1,\ y'=\bar
y,\ T=0,\ \xi=\bar\xi\lambda\ \text{and }\eta=\bar\eta\lambda.
\end{aligned}
\label{smallsupport}\end{equation}

Since we are interested in the scattering wavefront set of the image we can
localize near $\tx=\infty;$ the assumption on $b$ means the image is
already effectively localized in $y$ near $\bar y$ and in $t$ near $\bar
t.$ Thus, in local coordinates based at $\bar y=0,$ we take the
(normalized) Fourier transform of $v=(1-\phi(\tx))SBu,$
\begin{equation}
\mathcal{F} v(\xi',\mu;\tau)=
\int e^{i\tx\xi'+\tx y\cdot\mu}v(\tx,y;\tau)\tx^{n-1}d\tx dy.
\label{melwun1.210}\end{equation}
We then localize near some point using a cutoff $\tilde \phi (\tx, y, \xi',
\mu)$ to examine the regularity; such regularity allows us to conclude
absence of scattering wavefront set at points in the support of $\tilde
\phi$, where $\xi'$ and $\mu$ are coordinates in the fibers of the
scattering cotangent bundle determined by the canonical one-form
$$
\xi' d\tx + \mu \cdot d(\tx y).
$$

Finally, then, we arrive at the composite operator 
\begin{multline}
Gu=\mathcal{F}{(1-\phi(\tx))SBu},\\
Gu=\int e^{i\psi} c
u(\frac{\tx s}{\tau},y',t'-\frac{\tx T}{\tau})
\tau^{-2}\tx^{n+1}ds\, dy'\, dT\, dt'\, d\tx\, dy\, d\xi\, d\eta\, d\lambda,\Mwhere\\
c=c(\tx,y,\xi',\mu,t,t',\tau,
\frac{\tx}{\tau},s,T,\xi,\eta,\lambda)\\
=\tilde\phi(\tx,y,\xi',\mu)
a(t,t',\tau)
b(\frac{\tx}{\tau},y,t',s,T,\xi,\eta,\lambda)\Mand\\
\psi= \tx\xi '+\tx y\cdot\mu+(t-t')\tau+\frac i2\ang\tau(t-t')^2
+(s-1)\xi+(y-y')\cdot\eta+T\lambda.
\label{melwun1.213}\end{multline}

First suppose that $e(\tx,\lambda)$ is a conic cutoff keeping
$|\tx-\lambda|\le\delta \ang{\tx,\lambda}$ asymptotically for some small
$\delta >0$ and equal to one on a smaller region of this type. Inserting
$1-e(\tx,\lambda)$ into the integral we may integrate by parts using the
vector field 
\begin{equation}
V_1= \pa_T+\frac{\tx}{\tau}\pa_{t'}\text{ satisfying } V_1\psi=\lambda
-\tx + \frac{i \tx\ang\tau}{\tau}(t'-t).
\label{melwun1.214}
\end{equation}
Since $t-t'$ is small, for large $\tau$ this is elliptic in $\lambda$ and
$\tx$ on the support of the amplitude. This gives an operator of the form
\eqref{melwun1.213} but with amplitude arbitrarily decreasing in
$\tx,\lambda.$ Rapid decrease in $\lambda$ implies rapid decrease in $\eta$
and $\xi$ and hence again effectively replaces $B$ by an operator of
arbitrarily low order. Hence this term also has no uniform scattering
wavefront set.
 
Thus, by inserting the cutoff $e$ we may assume that $\tx\sim\lambda$ on
the support of the amplitude. Now, consider the vector fields 
\begin{equation}
V_2=\pa_s-\frac{\tx}s\pa_{\tx}-\frac{\tx}{\tau s}T\pa_{t'}\Mand
W=\pa_y
\label{melwun1.215}\end{equation}
which annihilate $u(\frac{\tx s}{\tau},y',t'-\frac{\tx T}{\tau})$ and satisfy
\begin{equation}
V_2\psi=\xi-\frac{\tx\xi'}s-\frac{\tx y\cdot \mu}{s}+\frac{\tx T}s
\left(1-i\frac{\ang\tau}{\tau}(t'-t)\right)\Mand W\psi=\eta+\tx\mu.
\label{melwun1.216}\end{equation}
These functions, with $t-\bar t,$ and $V_1\psi,$ define $\cR.$

Assuming that the support of $\tilde\phi$ in \eqref{melwun1.213} is
disjoint from the image of the essential support of $B$ under $\cR,$
integration by parts now allows the amplitude,
$c,$ in \eqref{melwun1.213} to be replaced by a symbol of arbitrarily low
order in $\xi,$ $\eta$ and $\lambda.$ Simpler versions of the arguments
used to prove \eqref{melwun1.202} and the boundedness of edge
pseudodifferential operators show that this leads to a map
\begin{equation}
x^lH^m_{\eo,c}(\bbR\times X)\longrightarrow
\CI(\bbR^n\times\bbR_t;\tau^{-l}L^2([1,\infty)).
\label{melwun1.218}\end{equation}
This proves \eqref{melwun1.212}.
\end{proof}

We shall apply the FBI transform to solutions of $\Box u=0.$ Since $S$ is a
Fourier integral operators, with complex phase function, the transform
satisfies a model equation involving the Laplacian, $\Lap_0,$
on the model cone $\tX=\RR_+ \times Y$ with respect to the the product
metric $dx^2+x^2 h_0(y,dy).$

\begin{lemma}\label{melwun1.311} There are operators $S_1$ and $S_2$ of
the same form as $S$ with amplitudes of order $\frac14$ and
$L\in\Diffb2(X)$ such that
\begin{equation}
(\Lap_0-1)S =-\tau^{-2}S\Box+\tau^{-\frac12}S_1+\tau^{-1}\tx^{-1}S_2 L.
\label{melwun1.313}\end{equation}
\end{lemma}

\begin{proof}
Computing directly
$$
\Box \psi u =\psi \Box u - 2 (D_x \psi) (D_x u) - (D_x^2 \psi) u.
$$
Let $\tilde S$ be defined as $S$ was but with
cutoff $\tilde\psi\in\CI(X),$ still supported in the product neighborhood
of the boundary, with $\tilde\psi\psi=\psi.$ Then
\begin{equation}
\tS(\Box\psi u)=S(\Box u)+\tS(L_1u)\Mwith\ L_1\in\Diffb1(X).
\label{melwun1.312}\end{equation}
Now the commutation relation for $D_t$ shows that
$$
\tS D_t^2 \psi u =\tau^2 \tS \psi u +\tau^{\frac32}S'u
$$
where $S'$ is of the same form as $S$ but with a different amplitude of
order $\frac14.$
Since $S \Lap_0=\Lap_0 S,$ $\Lap\psi-\Lap_0\psi \in x^{-1}
\Diffb2(X)$ and $\tS(\psi u)=Su$ we conclude that 
\begin{equation*}
(1-\Lap_0)S(u)=\tau^{-2}\tS \Box\psi
u-\tau^{-\frac12}S'u+\tS(\Lap-\Lap_0)\psi u,
\end{equation*}
yielding \eqref{melwun1.313}.
\end{proof}

\section{Reduced normal operator} \label{section:reduced}

The reduced normal operator of the conic wave operator, as an element of
the weighted edge calculus, is $1-\Lap_0$ where $\Lap_0$ is the Laplacian
for the tangent product-conic metric on $\tX=[0,\infty)_{\tx}\times Y.$
Regularity results obtained in subsequent sections depend on \emph{global}
invertibility properties for this operator. It is important to understand
that its behavior at infinity in the product cone is quite different from
its behavior near $\tx=0.$ In fact, near $\infty$, the metric
\eqref{melwun1.163} is a special case (again the model case) of a
scattering metric considered in \cite{Melrose43}. The local and
microlocal estimates in \cite{Melrose43} are combined below with the
analysis of the domain for conic metrics above to get the desired
invertibility estimates, which are of a standard type in scattering theory.

The analysis of the Friedrichs extension of $\Lap_0$ on $\tX$ proceeds very
much as above, using also the fact that the metric is complete at
infinity. Thus the scattering Sobolev spaces are, near infinity, the
natural metric Sobolev spaces and the domains of powers of $\Lap_0$ are
closely related to the ``b-sc'' Sobolev spaces defined in \eqref{melwun1.341}.

\begin{proposition}\label{melwun1.344} For the Friedrichs extension of
$\Lap_0,$ the Laplacian for the metric \eqref{melwun1.163}, 
\begin{equation}
\begin{gathered}
\Dom(\Lap_0^{\frac p2})=\tx^{-\frac n2+p}\ang\tx^{-p}\bscH p(\tX),\
-\frac n2<p<\frac n2,\\
\Dom(\Lap_0^{\frac n4+\frac\delta2})=\tx^{\delta}\ang\tx^{-\delta-\frac n2}\bscH
{\frac n2+\delta}(\tX)+\bbC\phi,\ 0<\delta<\delta_0
\end{gathered}
\label{melwun1.345}\end{equation}
for $\delta_0>0$ sufficiently small and with $\phi\in\CIc(\tX)$ identically
equal to $1$ near $\tx=0.$
\end{proposition}

The spectrum of $\Lap_0$ is the whole of $[0,\infty).$ Just as in
\eqref{melwun1.339} in the compact-conic case, $\Lap_0$ defines a
continuous map, 
\begin{equation}
\Lap_0:x^{\delta}\ang\tx^{-\delta-\frac n2}\bscH{\frac
n2+\delta}(\tX)+\bbC\phi\longrightarrow
x^{-2+\delta}\ang\tx^{2-\delta-\frac n2}
\bscH{\frac n2 -2+\delta}(\tX),\ \delta >0\Msmall,
\label{melwun1.346}\end{equation}
however this is never an isomorphism. On the other hand, the spectral family
$\Lap_0-\lambda,$ for $\lambda\in\bbC\setminus[0,\infty)$ does define a
continuous map as in \eqref{melwun1.346} which is an isomorphism.

One of the standard results of scattering theory, proved in
\cite{Melrose43} in the context of scattering metrics, is the \emph{limiting
absorption principle}. This asserts the existence of the limit of the
operator $(\Lap_0-\lambda\pm i\gamma)^{-1}$ for $\lambda\in(0,\infty)$ as
$\gamma\downarrow0.$ Of course the limit cannot exist as a bounded
operator inverting the resolvent, \ie\ cannot be defined on the range space
of \eqref{melwun1.346}. However it does exist on a somewhat smaller space.

\begin{proposition}\label{melwun1.347} For $\gamma>0$ the resolvent 
\begin{multline}
(\Lap_0-1\pm i\gamma)^{-1}:\tx^{-2+\delta}\ang\tx^{2-\delta-\frac n2}
\bscH{\frac n2 -2+\delta}(\tX)\longrightarrow\\
\tx^{\delta}\ang\tx^{-\delta-\frac n2}\bscH{\frac n2+\delta}(\tX)+\bbC\phi
\label{melwun1.348}\end{multline}
restricts to an operator  
\begin{multline}
(\Lap_0-1\pm i\gamma)^{-1}:x^{-2+\delta}\ang\tx^{2-\delta-\frac n2-s}
\bscH{\frac n2 -2+\delta}(\tX)\longrightarrow\\
x^{\delta}\ang\tx^{-\delta-\frac n2-s}\bscH{\frac n2+\delta}(\tX)+\bbC\phi
\label{melwun1.349}\end{multline}
for any $s$ and as an operator 
\begin{multline}
(\Lap_0-1\pm i\gamma)^{-1}:\tx^{-2+\delta}\ang\tx^{2-\delta-\frac
n2-\frac12-\epsilon} 
\bscH{\frac n2 -2+\delta}(\tX)\longrightarrow\\
\tx^{\delta}\ang\tx^{-\delta-\frac n2+\frac12+\epsilon}\bscH{\frac
n2+\delta}(\tX)+\bbC\phi,\ \Mforany\epsilon >0,
\label{melwun1.350}\end{multline}
the strong limit exists, as $\gamma\downarrow0.$
\end{proposition}

\begin{proof} The only difficulty with the convergence of the resolvent is
  related to the large, \ie scattering, end of the cone. Modulo compact
  errors the analysis in \cite{Melrose43} therefore applies and gives the
  same continuity properties for the limit of the resolvent as in the
  scattering case itself.
\end{proof}

Thus we have two limiting operators $(\Lap_0-1\pm i0)^{-1}$ with the
mapping property \eqref{melwun1.350}. For these operators it is not
possible to take $\epsilon<0$ in the domain space, nor in the range space
without further restriction. However there are ``larger'' spaces on which the
limiting resolvent is defined; these are fixed in terms of conditions on the
scattering wavefront set.

\begin{proposition} \label{melwun1.324} Let $\xi'$ be the dual variable to
  $d\tx$ in the scattering cotangent bundle.  Let $U_{\pm}$ be an open
neighborhood of the part of the radial set for $\Lap_0-1$ on which
$\pm\xi'>0$.  Then for any $m\in\bbR$ and $k'>k>-\frac12$
\begin{multline}
(\Lap_0-1\pm i0)^{-1}:\left\{f\in\tx^{-2+\delta}\ang\tx^{2-\delta-\frac n2+k}
\bscH{\frac n2 -2+\delta}(\tX);\scWFnoargs(u)\cap U_{\pm}=\emptyset\right\}
\longrightarrow \\
\left\{u\in \tx^{\delta}\ang\tx^{-\delta-\frac n2+k'}
\bscH{\frac n2+\delta}(\tX)+\bbC\phi;\scWFnoargs(u)\cap U_{\pm}=\emptyset\right\}
\label{melwun1.351}\end{multline}
and this is a two sided inverse of $\Lap_0-1$ on the union over $k$ and
$k'$ of these spaces.
\end{proposition}

\begin{proof} This just follows by combining the results of
\cite{Melrose43} with the analysis above of the domain of $\Lap_0.$ 
All the estimates of \cite{Melrose43} concern the behavior as
$\tx\longrightarrow \infty$ and the appearance of the conic boundary at
$\tx=0$ makes essentially no difference to the argument. This result could
also be proved using the analysis of the forward fundamental solution of
the $\RR_+$-invariant conic wave operator from \cite{Cheeger-Taylor2}.
\end{proof}

\section{Diffractive theorem}
\label{section:diffractive}

\begin{proposition}\label{melwun1.a316} Suppose
$u\in\cC(\bbR,\dcal_{\frac n2 -\delta }),$ for some $\delta >0,$ is an
admissible solution to the conic wave equation and
\begin{equation}
\WF(u)\cap R^{\epsilon}_{I}(0)=\emptyset
\label{melwun1.317}\end{equation}
for some $\epsilon>0.$ Then for any $\phi\in\CIc(\bbR\times X)$ with
support in a sufficiently small neighborhood of $\{0\}\times\pa X$
\begin{equation}
\phi u\in x^{-\delta'}H^{\infty}_{\bo}(\bbR\times X),\ \forall\ \delta'>0,
\label{melwun1.a318}\end{equation}
so $u$ is conormal near the boundary for $|t|$ small.
\end{proposition}

\begin{remark}
Note that, since $u$ is a solution of the wave equation,
\eqref{melwun1.a318} is equivalent to the condition that for all 
$k$, $\phi(x)D_t^ku$ is an $L^2$ function of $t$ near $t=0$ with values in
$\bigcap_s \dcal_s$, provided $\phi \in\CI(X)$ is chosen to have support
sufficiently near the boundary, \ie\ the conclusion of
Theorem~\ref{15.5.2000.7} holds with $s=\infty$.
\end{remark}

\begin{proof} First we replace $u$ by $u_+$ in terms of the decomposition
  \eqref{melwun1.388}, this has the same regularity properties but its
  Fourier transform in $t$ has support in $\tau>1.$ The treatement of $u_-$
  is similar. Now we proceed in three steps.

First, the hypothesis \eqref{melwun1.317} allows
Theorem~\ref{melwun1.a319} to be applied directly. Since
$u\in\cC(\bbR,\cE_{\frac n2})\subset
x^{-\delta}L^2_{\bo,\loc}(\bbR\times X),$ the first part of
Theorem~\ref{melwun1.a319} shows that 
\begin{equation}
\IC(0)\cap\eWF[m,-\delta](u)=\emptyset\ \forall\ m,\ \forall\ \delta>0.
\label{melwun1.321}\end{equation}

Our initial goal is to apply the identity in Lemma~\ref{melwun1.311} to
show that
\begin{equation}
\phi D_t^ku\in x^{-\delta}H^{-\infty}_{\eo}(\bbR\times X)
\label{melwun1.320}\end{equation}
for all $k$ and a cutoff as in the statement above.

The initial hypothesis on $u$ implies in particular that $u\in
x^{-\delta}H^{\frac n2-\delta}_{\eo}(\bbR\times X)$ near $\{0\}\times\pa
X.$ Applying Lemma~\ref{melwun1.205}, with the cutoff in the definition of
$S$ supported sufficiently near $x=0$ we therefore conclude that
\begin{equation}
Su\in\bigcap_{\alpha \in [0, n/2-\delta]} \tx^{-\delta}\ang{\tx}^{-\alpha}
\tau^{\delta}L^2_{\loc}\left([1,\infty)_{\tau}\times\bbR_t; \bscH{\frac
n2-\delta-\alpha}(\tX)\right) \Mnear\bar t=0.
\label{melwun1.322}\end{equation}
Similar conclusions apply to the terms on the right
in \eqref{melwun1.313}, with $S_1u$ having the regularity property
\eqref{melwun1.322} and
\begin{multline}
S_2(Lu)\in\bigcap_{\alpha \in [0, n/2-\delta-2]}
\tx^{-\delta}\ang{\tx}^{-\alpha} \ang\tx^{\frac n2-\delta}
L^2_{\loc}\left([1,\infty)_{\tau}\times\bbR_t;\bscH{\frac
n2-\delta-2-\alpha}(\tX)\right)\\ \Mnear\bar t=0.
\label{melwun1.323}\end{multline}
It is also the case that near $\tx=0,$ $Su$ is in the domain of
$\Lap_0^{\frac n2+\delta}$ if $\delta >0$ is small enough. The identity
\eqref{melwun1.313} then holds, near $\tx=0,$ in the sense of the domain of
$\Lap_0^{\frac n2-2+\delta},$ which is just a weighted b-Sobolev space.

Of fundamental importance is the estimate on the scattering wavefront set
which follows from \eqref{melwun1.321}, together with the initial
assumption \eqref{melwun1.317}, by use of
Proposition~\ref{melwun1.211}. Namely, the incoming wavefront set for all
terms, $Su,$ $S_1u$ and $S_2(Lu),$ computed with respect to the spaces in
\eqref{melwun1.322} and \eqref{melwun1.323} is absent---these functions are
microlocally rapidly decreasing as $\tx\to\infty$ in this uniform
sense. Now, Lemma~\ref{melwun1.324} shows that solving \eqref{melwun1.313}
gives an ``improvement'' in the estimate \eqref{melwun1.322} by a factor of
$\tau^{-\frac12}$ at the expense of more growth as $\tx\to\infty.$ Since we
may invert $S,$ up to errors which involve a smoothing operator in $t,$ we
may iterate this argument, reapplying Theorem~\ref{melwun1.a319} to obtain
absence of scattering wavefront set with the new $\tau$ weight and using
the improved estimate on $Su$ to estimate $S_1u$ and $S_2(Lu).$ Such
iteration yields an estimate
\begin{equation}
\begin{aligned}
Su\in \bigcap_{\alpha} \tx^{-\delta}\ang\tx^{N}\tau^{-M}\bscH{\frac
n2-\delta}(\bbR\times[1,\infty)\times\tX) \Mnear\bar t=0,\ \forall
M\Mwith N=N(M).
\end{aligned}
\label{melwun1.325}\end{equation}
Using the reverse regularity estimate in Lemma~\ref{melwun1.205} we
conclude that \eqref{melwun1.320} does indeed hold.

Finally we apply Theorem~\ref{melwun1.a319} again to deduce the claimed
regularity \eqref{melwun1.317}. Indeed, we may apply the first three parts
of the theorem, in succession, to $D_t^ku$ for any integer $k.$ As before
we find that over $t=0$ and the boundary
\begin{equation}
\eWF[M,-\delta](u) \subset\OG(0),\ \forall\ M.
\label{melwun1.326}\end{equation}
Now we may use the outgoing propagation result to conclude that
$$\OG(0)\cap\eWF[\frac {n-1}{2}-\delta,-\delta](D_t^ku)=\emptyset$$ for any
$\delta>0.$ Since this is fixed finite regularity, certainly implying that
$D_t^ku\in L^2_g(\bbR\times X)$ locally for all $t$ derivatives,
\eqref{melwun1.a318} follows.
\end{proof}

\begin{proof}[Proof of Theorem~\ref{15.5.2000.7}] 
We use the decomposition \eqref{melwun1.338} to assume without loss of
generality that there are no singularities in $\{\tau<0\}$.  By
time-translation invariance, we may assume that $\bar t=0.$ Choose a
cutoff $\psi(t)\in\CI(\bbR)$ which is $1$ in $t>-\frac12\delta$ and
$0$ in $t<-\delta$ for $\delta>0$ which will be chosen small.  Then
$v=\psi(t)u$ satisfies
\begin{equation}
\Box v=f=-2i\psi'(t)D_tu-\psi''(t)u.
\label{melwun1.327}\end{equation}
Now, if $\delta>0$ is small enough, the hypothesis on the wavefront set of
$u$ means that $\WF^{r-1}(f)$ is disjoint from all incoming rays
arriving at the boundary at time $0.$ Thus, it may be divided into two
pieces $f=f_1+f_2$, with both $f_i$ supported in
$-\frac14\delta>t>-2\delta$, $f_1$ supported away from the boundary and
satisfying
\begin{equation}
f_1\in \cC^p(\bbR;\dcal_{r-1-p}),\ \forall\ p\in\bbN
\label{melwun1.328}\end{equation}
and with $\WF(f_2)$ disjoint from the incoming cone at $t=0.$ It follows
that the forward solution to $\Box v_1=f_1$ is in $H^r(\bbR\times X)$ away
from the boundary and is in $\mathcal{C}(\RR; \dcal_r)$ near it. Thus it
suffices to consider $\Box v_2=f_2$ which is a solution near $t=0$ with no
incoming singularities at all at time $t=0;$ by the uniqueness of solutions it is
equal to $u-v_1$ near $t=0.$ Now $v_2$ may be extended to a global solution
$u'$ which is equal to $v_2,$ and hence to $u-v_1$ near $t=0$ and this too
has no incoming singularities at time $0.$ After a finite amount of
smoothing in $t$ using Lemma~\ref{lemma:orderchange},
Proposition~\ref{melwun1.a316} applies to $u',$ showing that it has no
outgoing singularities at time $0$ (and indeed lies in $\mathcal{C}(\RR;
\dcal_\infty)$) locally.
\end{proof}

\section{Propagation of tangential regularity}
\label{section:tangential}

As a prelude to the division theorem, we will prove that regularity of
solutions in the tangential, \ie\ boundary, directions is conserved under
time-evolution.  This result represents a global version, in energy spaces,
of the microlocal estimates in Theorem~\ref{melwun1.a319}, part
\ref{melwun1.353}.

Recall that for a boundary component $Y$ of $X$, $\Lap_Y$ denotes the
Laplace-Beltrami operator induced on $Y$ by the metric $h_0$.  Set
$$
Y_s = (1+\Lap_Y)^{s/2},\ s \in \RR.
$$
These tangential pseudodifferential operators act naturally on the
boundary, but may be viewed as acting on the fibers of the product
decomposition in Theorem~\ref{thm:normalform}. Thus they act on functions
or distributions on $(a,b)_t\times\pa X\cap [0,\ep)_x.$

Recall that $\E_s$ denotes the energy space $\dcal_s \oplus \dcal_{s-1}$
where $\dcal = \Dom (\Lap^{s/2}).$ We write $\E = \E_1$ as well as
$\dcal=\dcal_1$ for convenience. $\E_1$ is a Hilbert space with the norm
$$
\norm{(u,v)}_\E^2 = \norm{u}^2_{L^2_g} + \norm{du}^2_{L^2_g}+ \norm{v}^2_{L^2_g},\
(u,v)\in\E.
$$
Note that if we let
$$
M = \begin{pmatrix} 0 & 1 \\ \Lap & 0 \end{pmatrix}
$$
denote the infinitesimal generator of $U(t)$, $M$ is not selfadjoint with
respect to the norm on $\E$ as defined here (owing to the
$\norm{u}_{L^2}^2$ term).

For use in the sequel, we note the following energy estimate: For an
admissible solution $(u, D_t u)$ with $(u,
D_t u) \in \CI (\RR; \E_\infty)$, and an operator $Q$ such that
\begin{equation}
Q: \CI(\RR; \dcal_\infty) \to \CI(\RR; \dcal_{n/2})
\label{hypothesis}\end{equation}
set
$$
\Psi(t) =
\begin{pmatrix} Q u(t)\\ D_t (Q u(t))\end{pmatrix}.
$$
We can then compute
\begin{equation}
\begin{aligned}
\frac 12 \frac d{dt} \norm{\Psi(t)}^2_\E &=
\Re \ang{\begin{pmatrix} i D_t (Q u) \\ i D_t^2 (Q u) \end{pmatrix}, \Psi}_\E
\\
&= \Re\ang{i M \Psi, \Psi}_{\E} + \Re \left\langle{\begin{pmatrix} 0 \\
i [\Box,Q] u
\end{pmatrix}, \Psi}\right\rangle_{\E}\\
&=\Re \ang{i [\Box, Q] u, D_t (Q u)}_{L^2_g}+ \Re \ang{i D_t (Q u), Q u}_{L^2_g}.
\end{aligned}
\label{general:energycons}
\end{equation}

Let $\ocal$ denote a product neighborhood of the boundary and let
$\Ec(\mathcal{O})\subset\E$ and $\dcalc(\ocal) \subset \dcal$ be the
subspaces with compact supports in $\mathcal{O}.$

\begin{lemma}\label{24.10.2000.1} 
If $B(x)$ is a tangential pseudodifferential operator of order $1$
depending smoothly on $x$ and $B(0)\bbC=0$ then
$x^{-1}B(x):\dcalc(\mathcal{O})\longrightarrow L^2(X)$.
\end{lemma}

\begin{proof} We know that operators of the form $x^{-1}V$ with $V$ a
tangential vector field are bounded in this way, as is $B(x)$ itself.
The generalized inverse of the tangential Laplacian satisfies
$E\Lap_Y=\Id-\Pi_0$ where $\Pi_0$ is orthogonal projection onto the
constants.  We have $B(x)\Pi_0=x\tilde B(x)$, hence
\begin{equation}
B(x)=B(x)E\Lap_Y+B(x)\Pi_0=\sum\limits_{j}B_j(x)V_j+x\tilde B(x)
\label{24.10.2000.2}\end{equation}
where $V_j$ are vector fields tangent to $\pa X$, the $B_j$ are of order
$0$ (and hence bounded on $L^2)$ and $\tilde B(x)$ is again of order $1$ and
smooth in $x.$ Thus
\begin{equation}
x^{-1}B(x)=\sum\limits_{j}B_j(x)x^{-1}V_j+\tilde B(x)
\label{24.10.2000.3}\end{equation}
is bounded as claimed by \eqref{Friedrichs}.
\end{proof}

Let $\ocal' \subset \ocal$ be a smaller product neighborhood of $\pa X$.

\begin{proposition}
For all $\abs{t}$ sufficiently small and $s \in \bbR$,
$$
Y_s U(t) Y_{-s}: \Ec(\ocal') \longrightarrow \E
$$
is bounded.
\label{prop:tangential}
\end{proposition}

\begin{proof}
As $(\dCI (X)+\bbC) \oplus \dCI(X)\subset \E_\infty$ is dense in $\E,$ it
suffices to prove the relevant estimate for Cauchy data in $\E_{\infty}.$

Let $\Phi(t)=(u(t), D_t u(t)) = U(t) (u_0, u_1)$ be the
solution to the Cauchy problem \eqref{15.5.2000.1}--\eqref{15.5.2000.3}
with $\supp u_0,\ \supp u_1 \subset \ocal'.$ Then there is an open interval
$I$ containing $0$ such that for $t\in I,$ we have $\supp u(t),$ $\supp D_t u(t)
\subset \ocal.$ The hypothesis \eqref{hypothesis} is satisfied for
$Q=Y_s$, hence
\begin{equation}
\frac 12 \frac{d}{dt} \norm{Y_s \Phi(t)}^2_\E = \Re \ang{Y_s u_t, Y_s u} +
\Re \ang{[\Lap, Y_s] u, Y_s u_t}
\label{energyest}
\end{equation}
 for all $t\in I$.  By Lemma~\ref{lemma:morecommutators} and
Lemma~\ref{24.10.2000.1}, $\norm{[\Lap, Y_s]u}$ is bounded by a multiple of
$\norm{(u,0)}_\E$, locally near $Y$.  Hence
\begin{equation*}
\frac 12 \frac{d}{dt} \norm{Y_s \Phi(t)}_{\E}^2\le
C\norm{Y_s \Phi(t)}_{\E}^2
\label{melwun1.174}\end{equation*}
from which the boundedness follows. 
\end{proof}

\section{Global propagation of conormality}
\label{section:radial}
The tangential regularity discussed in the previous section is the main
step to showing that incoming conormal waves with respect to the
surface $R_{\pm, I}$ propagate through the boundary to be conormal on the
outgoing radial surface $R_{\pm, O}$. To prove this we need further to show that
regularity is preserved under the repeated action of the radial vector field 
\begin{equation*}
R=x D_x + (t-\bar t) D_t.
\end{equation*}
where $t=\bar t$ is the time at which the cone hits the boundary. Using
time-translation invariance we may always assume $\bar t=0.$ 

To begin, note that we can combine identities from
Lemmas~\ref{lemma:commutators} and \ref{lemma:morecommutators} to obtain,
for any fixed $s,k$,
\begin{equation}
[\Box, Y_s R^k] = \sum_{j=0}^{k-1} c_j  Y_{s+1} R^j \Box
+ \sum_{j=0}^{k-1} (a_j D_x + x^{-1} P_j)  Y_{s+1} R^j
 +  (R D_x + x^{-1} S) Y_s R^k
\label{formalcommutator}\end{equation}
where $c_j \in \CC$, $a_j \in \CI([0, \ep) \times Y)$ and $R \in
\CI([0,\ep) \times Y; \Psi^{-1}(Y))$, and where $S,\ P_j\in \CI([0,
\ep)\times Y; \Psi^1(Y))$ annihilate constants at $x=0$.  Since the
operator $Y_s R_k$ satisfies the hypothesis \eqref{hypothesis},
\eqref{formalcommutator} makes sense when applied to an element of
$\CI(\RR; \dcal_\infty)$, with the equality holding in $\CI(\RR;
\dcal_{n/2-2})$. Thus using Lemma~\ref{24.10.2000.1}, the energy identity
gives
\begin{multline}
\frac 12 \frac d{dt} \norm{(Y_i R^j u, D_t Y_i R^j u)}_\E^2 \\ \leq C
\norm{(Y_i R^j u, D_t Y_i R^j u)}_{\E} \sum_{i'+j'\leq i+j} \norm{(Y_{i'} R^{j'} u, D_t
Y_{i'} R^{j'} u)}_{\E}.
\label{estimate1}
\end{multline}

\begin{proposition}\label{prop:conormality}
If $U(t)$ is the solution operator to the Cauchy problem
\eqref{15.5.2000.1}--\eqref{15.5.2000.2} then for a product 
boundary neighborhood $\ocal,$ small time $T$ and each $k \in \NN,$
\begin{equation}
\sum\limits_{i+j\le k}\|Y_{i}R^jU(t)(u_0, u_1)\|_{\E}\le C \sum\limits_{i+j\le k}
\|Y_{i}R^j(u_0, u_1)\|_{\E},\ |t|\le T,\ (u_0, u_1) \in\Ec(\ocal).
\label{melwun1.179}\end{equation}
\end{proposition}

\begin{proof}
We apply \eqref{estimate1} inductively. Note that
$$
\sum\limits_{i+j\le k}\|Y_{i}R^jU(t)(u(t), D_t u(t))\|_{\E} \sim
\sum\limits_{i+j\le k}\|\left(Y_{i}R^j u(t), D_tY_i R^ju(t)\right)\|_{\E}
$$
since
$$
[D_t, Y_i R^j] = \sum_{l=0}^{j-1} c_l Y_i R^l D_t = Y_i
\sum_{l=0}^{j-1} c_l' D_t R^l.\qed
$$
\noqed
\end{proof}

\begin{proof} [Proof of Theorem~\ref{thm:conormality}]

By results going back at least to Hadamard (see also \cite{Hormander1,
Duistermaat-Hormander1}), under the hypotheses of the theorem, there exists
$\ep>0$ such that for $0<t<\ep$, $u(t)$ and $D_t u(t)$ are conormal
distributions in $\RR\times X$ at $\{x=\bar x-t\}\cup \{x=\bar x+t\}$ with
respect to $H^s(\bbR\times X)$ and $H^{s-1}(\bbR\times X)$
respectively. The distribution conormal to $\{x=\bar x+t\}$ remains
conormal for small positive time, as it does not  reach $\pa X;$ we thus
assume without loss of generality that this component vanishes. Hence for
all $s,$ $k$
$$
Y_s R^k (u, D_t u) \in \ccal(I_0; \E_s)
$$
for some interval $I_0$ containing $0$ but not necessarily $\bar x$.

Note that with $\Theta_s$ from Definition~\ref{Theta-op}, $[\Theta_s, R]$
is a properly supported pseudodifferential operator of order $s$ on $\RR$, hence
\begin{equation}
[\Theta_s, R]=Q\Theta_s +E
\label{thetacommutator}
\end{equation}
where $Q \in \Psi^0(\RR)$ and $E \in \Psi^{-\infty} (\RR)$.  Since
convolution with a properly supported time-translation invariant
pseudodifferential operator of order $0$ on $\RR_t$ maps finite energy
solutions to finite energy solutions, iterative application
\eqref{thetacommutator} shows that for all $i$ and $j$
$$
Y_i R^j (\Theta_{s-1} u,\Theta_{s-1} D_t u) \in
\ccal(I; \E),
$$
where $I\subset I_0$ is a time interval containing $0$.

Now by Proposition~\ref{prop:conormality}, $Y_i R^j \Box^k \Theta_{s-1}
(u,D_t u) \in \E$ for all $i,j,k\in \NN$ and for all $t<T$.  The symbols of
$Y_1$, $R$, and $\Box$ are defining functions for the conormal bundle to the
hypersurface $x=\abs{t-\bar x}$.  Hence $(\Theta_{s-1} u(t),\Theta_{s-1} D_t u(t))$ is
conormal to $x=\abs{t-\bar x}$ in $(H^1, L^2)$.  By Lemma~\ref{lemma:orderchange},
$(u(t), D_t u(t))$ is conormal to $x=\abs{t-\bar x}$ in $(H^s,
H^{s-1})$.  The theorem then follows by restriction to fixed $t$.
\end{proof}

\section{Proof of the division theorem}
\label{section:divthm1}

In this section, we prove Theorem~\ref{melwun1.34}.  We begin with a
preparation argument, allowing us to replace our hypothesis of regularity
at $R_{\pm, I}^\ep$ with global regularity of the same order.  As in the
proof of Theorem~\ref{15.5.2000.7}, we may assume without loss of
generality that $u$ has edge wavefront set only in $\{\tau>0\}$.

Let $i:0 \times X^\circ \hookrightarrow \RR\times X^\circ$ be the inclusion
map, and $i^*$ the induced contravariant map on cotangent bundles. Under
the assumption of the nonfocusing condition (from Definition~\ref{NC}),
there exists a microlocal neighborhood $U$ of $i^* R^\ep_{I}(\bar t)$
and $k \in \NN$ such that $\WF^{r+l} (Y_{-k} u\restrictedto_{t=0}) \cap U =\emptyset.$ We
construct a microlocalizer in such a neighborhood which preserves
tangential regularity.

\begin{lemma}\label{melwun1.366}
  Let $U \subset T^*X^\circ$ be an open conic neighborhood of $i^* R_{I}
  (\bar t).$ There exist a smaller neighborhood $V$ and an operator $H \in
  \Psi^0(X^\circ)$, with Schwartz kernel compactly supported in an
  arbitrarily small neighborhood of the diagonal, such that $\WF' H \subset
  U$, $\WF' (I-H) \subset V^\complement$, and $[Y_k, H]=0$ for all $k
  \in \ZZ$.
\end{lemma}
\begin{proof}
Let $\psi_\ep \in \CIc(\RR)$ equal $1$ for $x<\ep/2$ and $0$ for $x>\ep$. 
Let $\chi \in \CI (\RR)$ vanish for $x<0$ and equal $1$ for $x>1$.
If $\ep$ is sufficiently small then 
$$
\supp \psi_\ep(\abs{x-\bar t})
\psi_\ep(\abs\eta^2/(\ang\xi^2+\abs\eta^2)) \subset U;
$$
here $\eta$ is dual to $\pa_y$ and $\xi$ to $\pa_x$ in $T^* X^\circ$.  Let
$$
H=\psi_\ep(\abs{x-\bar t}) \psi_\ep(\Lap_Y/(\ang{D_x}^2 +
\Lap_Y)) \psi_\ep(\abs{x-\bar t}),
$$
where the function of the operator $\Lap_Y/(\ang{D_x}^2 + \Lap_Y)$ may be
constructed on the manifold $S^1 \times X$ and its Schwartz kernel cut off
in $x$ and glued into a product neighborhood of $\{x=\bar t\}.$ Then $H$
commutes with $\Lap_Y$, and hence with $Y_k$ for all $k$.
\end{proof}

Since the kernel of $H$ vanishes near the boundary it preserves the domains
of all powers of $\Lap.$ Let $(v(t), D v(t))$ be the solution to the Cauchy
problem \eqref{15.5.2000.1}--\eqref{15.5.2000.2} with initial data $(Hu(0),
H D_t u(0)).$ By construction, $\WF (u-v) \cap R_{I}(\bar
t)=\emptyset$, so by the diffractive theorem, $u-v \in x^{-\delta}
\bH\infty (I \times U)$ for some open interval $I \ni \bar t$, $U \supset
\pa X.$ Hence it suffices to prove the desired results for $v$ instead of
$u$.

Since
$$
Y_{-k} (v(0), D_t v(0)) = H Y_{-k} (u(0), D_t v(0)),
$$
if the microlocalizer $H$ is chosen concentrated sufficiently near the
incoming set, the nonfocusing condition yields
$$
Y_{-k} (v(0), D_t v(0)) \in \E_{r+l}.
$$

Proposition~\ref{prop:tangential} now yields
$$
Y_{-k} \Theta_{s-1} U(t) (v(0), D_t v(0)) \in \E,
$$
hence
$$
Y_{-k} U(t) (v(0), D_t v(0)) \in \mathcal{C}(\RR; \E_{s})
$$
and the first part of the division theorem (that under the assumption of
the nonfocusing condition) now follows by
Proposition~\ref{prop:domainsinedgeterms}.

We now prove the second part of the theorem (under the assumption of the
conormal nonfocusing condition).  For simplicity, we now translate the time
variable so that $\bar t=0$.

We have
\begin{multline*}
R Y_{-k} (v, D_t v)\restrictedto_{t=0} \\= (x
D_x H Y_{-k} u(0) + t Y_{-k}  D_t u(0),\ x D_x H Y_{-k} D_t u(0) +
t \Lap H Y_{-k} u(0)).
\end{multline*}
Since $[H, \Lap]$ and $[H, x D_x]$ are pseudodifferential operators of
order $1$ and $0$, compactly supported in $X^\circ$, this
is in $\E_{r+l}$ provided the conormal nonfocusing condition holds.

Now we prove the second part of the division theorem for the special case $r+l=1$.

By Proposition~\ref{prop:conormality}
$$
\Box R Y_{-k} v =f \in \ccal(\RR; \dcal_{-1}),
$$
hence
$$
(R Y_{-k} v (t), D_t R Y_{-k} v(t)) = U(t) [(R Y_{-k} v, D_t R Y_{-k} v)\restrictedto_{t=0}] + \int_0^t
U(t-s) (0, f) ds.
$$
The first term in the right is in $\ccal(\RR; \E_1)$ by the
conormality assumption. The second is in $\ccal^1(\RR; \E_0)$, hence in
$\ccal(\RR; \E_1).$ Thus, $R Y_{-k} v(t) \in \ccal(\RR; \dcal_1)$, \ie\
$$
D_t Y_{-k} v \in t^{-1} [\ccal(\RR; \dcal_1) + x D_x \ccal(\RR; \dcal_1)].
$$
By Lemma~\ref{melwun1.187}, we obtain for all $\ep>0$,
$$
D_t Y_{-k} v \in t^{-1} \ccal(\RR; x^{-n/2+1-\ep} \bL (X)) \subset t^{-1/2-\ep}
x^{-n/2+1-\ep} \bL(I\times X)
$$
for all $\ep>0$.
We also know a priori that
$$
D_t Y_{-k} v \in \ccal(\RR; \dcal_0) \subset t^{1/2} x^{-n/2} \bLloc(\RR\times X).
$$
Hence by interpolation,
$$
D_t Y_{-k} v \in x^{-n/2+1/2-\ep} \bL(I\times X) \text{ for all } \ep>0,
$$
proving the theorem in the special case $r+l=1$.

If $r+l<1$, we apply Lemma~\ref{lemma:orderchange} and
\eqref{thetacommutator} to conclude
\begin{equation}
\Theta_{r+l-1} D_t Y_{-k} v \in x^{-n/2+1/2-\ep} \bL (I\times X).
\label{timederiv}
\end{equation}
Thus
$$
D_t Y_{-k} v \in x^{-n/2+1/2-\ep} H^{r+l-1} (I; \bL(X)),
$$
and since $r+l-1<0$, Proposition~\ref{prop:sobolevrelationship} yields
$$
D_t v \in x^{-n/2+1/2-\ep} \eH{r+l-k-1} (I\times X).\qed
$$

\section{Consequences of the division theorem}
\label{section:genprop}

In this section, we record two consequences of the division theorem.
First, we prove a slightly weakened version of Theorem~\ref{15.5.2000.8},
the geometric propagation theorem.  Then we prove
Theorem~\ref{melwun1.367}, establishing conormality of the diffracted front.

The weakened version of Theorem~\ref{15.5.2000.8} we now establish is as
follows
\begin{partialthm}
If $u$ is an admissible solution to the conic wave equation satisfying the
nonfocusing condition in Definition~\ref{NC} at $\{\bar t\}\times Y$ for some
$r\in\RR$ and $l\in(0, n/2)$ and if $p \in R_{\pm, O}^\ep (\bar t, Y)$,
then
\begin{equation}
\Gamma^\ep (p) \cap \WF^{r'}(u) = \emptyset\Mfor r' \in (r, r+l-1/2)
\Longrightarrow p \notin \WF^{r'-\delta}(u)\ \forall \delta>0.
\label{melwun1.333}\end{equation}
If $u$ satisfies the conormal nonfocusing condition at $\{\bar t\}\times Y,$ with
$k=1$, then \eqref{melwun1.333} holds for all $r'\in (r, r+l).$
\end{partialthm}
In other words, the theorem holds almost as originally stated, subject to
limits on $l$ and to the stronger conormal version of the nonfocusing
condition, or yields one half derivative less, subject to the original
hypothesis (again for a limited range of $l).$

\begin{proof}
Let
$$
(v(t), D_t v(t)) = \Theta_{r+n/2-1} (u(t), D_t u(t)).
$$
By \eqref{thetacommutator}, the new solution $v$ satisfies the nonfocusing
condition, with the same $l.$ It also 
satisfies the conormal nonfocusing condition if $u$ does, with $r+l\leq
1$. The division theorem now implies that
there exists an interval $I \ni \bar t$, and neighborhood $\ocal \supset
\pa X$ in $X$, such that
$$
D_t v\in x^{-n+r'-r+1/2}H^{-\infty}_{\eo}(I \times \ocal),
$$
where we allow $r' \in (r,r+l-1/2)$ resp.\ $(r, r+l)$ in the case where $u$
satisfies the nonfocusing condition resp.\ conormal nonfocusing condition.
On the other hand, the incoming regularity assumption shows that $\Gamma_\ep
(p) \cap \WF^{-n/2+r'-r}( D_t v) =\emptyset$.
Theorem~\ref{thm:propagationthrough} now yields $p \notin
\WF^{-n/2+r'-r-\delta} (D_t v),$ hence $p \notin \WF^{r'-\delta}( u).$
\end{proof}

We now prove Theorem~\ref{melwun1.367} by including the effect of iterated
regularity with respect to the vector field $R$ defined in \eqref{melwun1.175}.

\begin{proof}[Proof of Theorem~\ref{melwun1.367}]
We will show that away from the geometrically
propagated rays, the solution maintains its regularity under application
$Y_i R^j$ for all $i\in 2\NN$, $j \in \NN$.  The proof thus amounts to a
microlocalized version of the argument used previously to prove
conservation of radial conormality.

Using Lemma~\ref{lemma:orderchange}, we may assume that
$r=1$.

By Proposition~\ref{prop:conormality}, for all $i\in 2\NN$, $j \in \NN$,
there exists $N=N(i,j)<0$ such that
$$
Y_N Y_i R^j (u, D_t u) \in \mathcal{C}(\RR; \E)
$$
for all $t$ small.  Thus, $Y_i R^j D_t u \in x^{-\ep} \eH{-\infty}(\RR\times X)$,
locally near $t=\bar t$.

We may apply the first two parts of Theorem~\ref{melwun1.a319} to the
solution $D_t u$ and conclude that, for all $m$ and $\ep>0,$
$\eWF[m,-\ep](D_t u)$ is disjoint from all incoming bicharacteristic
segments, into and within the boundary, with endpoints $p.$ Using the
interpolation of wavefront sets in Proposition~\ref{prop:wfinterp} we
conclude that the same holds for $\eWF[m,-\ep](Y_i R^j D_t u)$ for all $i$,
$j$ and $m.$ The third and fourth parts of Theorem~\ref{melwun1.a319}. Then
shows that $p\notin\eWF[m,-\ep](Y_i R^j u)$ provided $m<-\ep.$ This
establishes conormality.
\end{proof}

\section{Conormal Cauchy data and the proof of Theorem~\ref{15.5.2000.8}}
\label{section:conormal}

The most significant class of examples to which the conormal nonfocusing
condition applies is given by solutions $u$ with Cauchy data conormal to a
hypersurface $W$ which is at most simply tangent to the radial surfaces.

\begin{lemma}\label{lemma:conormaldata} Let $W$ be a compact smooth
hypersurface in $X^\circ$ such that, in a collar neighborhood of the boundary,
$dx$ restricted to $W$ vanishes at only finitely many points $p_j=(x_j,y_j),$
$j=1,\dots,M,$ at each of which $W$ is simply tangent to $x=x_j$ and suppose
$u_0\in H^r(X)$ is conormal with respect to $W$ then for any $N\in \NN$,
there exists $k\in\NN$ such that
\begin{equation} Y_{-k} u_0 \in H^s_{\loc}(X^\circ),\ V_1\dots V_p Y_{-k} u_0
\in H^s_{\loc}(X^\circ),\ p\le N,\ s<r,
\end{equation}
for all $V_i\in\CIc(X^\circ; TX)$ that are tangent to $\{x=x_j\}$ for all $j.$
\end{lemma} 

\begin{proof} Without loss of generality, we may localize and take $M=1.$
Since $Y_{-k}\in\Psi^{-k}(Y),$ for any $K$ we may choose $k$ so that
$\kappa(Y_{-k})(y,y')$, the Schwartz kernel of $Y_{-k}$, is in
$\ccal^{K}(Y\times Y);$ we choose $K>(n-1)/2+N$.

The distribution $u$ is of the form
$$
u(x,y) = \int c(x,y,\xi) e^{i \phi(x,y)\xi} \; d\xi,
$$
where $c\in S^{-r-1/2+\ep}(\RR\times X^\circ)$ (for all $\ep>0$) is a
symbol with one phase variable and $\phi$ is a defining function for $W$.
Hence
$$
Y_{-k} u = \int c(x,y',\xi) \kappa(Y_{-k})(y,y') e^{i \phi(x,y')\xi} \; d\xi\; dy'.
$$
If $\phi \neq 0$ or $d_y \phi \neq 0$ on $\supp c$, integration by parts
using $\phi^{-1} D_\xi$ or $(\xi \phi'_{y'})^{-1} D_{y'}$ shows that $Y_{-k} u
\in H^{K+r}(X^\circ).$  Hence may assume that, on the support of $c,$ there
is one point, which in local coordinates we may 
take to be $x=\bar x,$ $y=0,$ at which $\phi=0=d_y\phi.$ The hypothesis of simple
tangency and the Morse Lemma allow us to arrange that, locally, $\phi(x,y) =
(x-\bar x)-\sum_j \sigma_j y_j^2$ with $\sigma_j=\pm.$ The method of
stationary phase now yields
$$
Y_{-k} u = \int \ang\xi^{-(n-1)/2} \tilde c(x,y,\xi) e^{i(x-\bar x)\xi} \; d\xi
$$
where $\tilde c$ satisfies any finite number of the symbol estimates of
any order greater than $-r-1/2$ provided that $K$ is taken large
enough. Thus $Y_{-k} u \in H^{r+(n-1)/2-\ep}_{\loc}(X^\circ)$ for any
$\epsilon >0.$ Furthermore, as long as $\alpha +\abs\beta\leq N$, integration
by parts shows  that $((x-x_1) D_x)^\alpha D_y^\beta Y_{-k} u$ remains in 
this space.
\end{proof}

Such conormality for the initial data leads to solutions satisfying the
conormal nonfocusing condition.

\begin{proposition}\label{prop:smoothedconormality} Suppose that
$(u_0, u_1)\in \E_r(X)$ vanishes near $\pa X$ and $u_0$ and $u_1$ are
conormal with respect to a hypersurface $W \subset X^\circ$ as in
Lemma~\ref{lemma:conormaldata} then, for small $t \in \RR,$ the solution
$u$ to the Cauchy problem problem \eqref{15.5.2000.1}--\eqref{15.5.2000.3}
satisfies the conormal nonfocusing condition with background regularity $r$
and relative regularity $l<(n-1)/2$ to all orders $k \in \NN.$
\end{proposition}

\begin{proof} By localization of the initial data we may again assume that
there is just one point of tangency of $W$ to one radial surface, $x=c.$
Since the conormal nonfocusing condition is trivial microlocally away from
the radial directions we may assume that $u_0$ and $u_1$ are supported close
to the point of tangency. Initial data for the wave equation which is
conormal with respect to a hypersurface, $W,$ gives rise, by Huygen's
principle to a solution which is conormal, for $t>0$ small, to the union of
the two characteristic hypersurfaces $W_\pm$ through $W;$ this follows from
the original construction of Lax or from \cite{Hormander1}. Here $W_\pm$
are each tangent to the radial surface $\{x\pm t=c\}.$ Thus $W_-$ is
outgoing: the bicharacteristics forming it do not hit the boundary small
times. Now, we may apply Lemma~\ref{lemma:conormaldata} above, regarding
$t$ as a parameter, to conclude that the solution satisfies the conormal
nonfocusing condition with background regularity $r$ for any $l<(n-1)/2$ to
all orders.
\end{proof}

\begin{proof} [Proof of Theorem~\ref{melwun1a.374}]
We can write the fundamental solution in the form
$$
E_{\bar m}=\frac{\sin t\sqrt\Lap}{\sqrt\Lap} \delta_{\bar m}=U(t)(0,i\delta
_{\bar m})
$$
where $\bar m\in X^\circ$ lies sufficiently close to $\pa X.$
Then for all $\ep>0,$ the initial data
$$
(0,i\delta _{\bar m}) \in H^{-n/2+1-\ep}(X^\circ)\oplus H^{-n/2-\ep}(X^\circ)
$$
hence the solution lies is $\E_{-n/2+1-\ep}$ for all $\ep>0.$

For small $t>0$ $(E_{\bar m}(t), D_t E_{\bar m}(t))$ are conormal with
respect to the hypersurface $W$ which is the geodesic sphere of radius
$t$ around $\bar m.$ This is tangent only to the two radial surfaces
$x=d(\bar m,\pa X)\pm t$ and the tangency is certainly simple. Hence
Proposition~\ref{prop:smoothedconormality} applies and shows that $u$
satisfies the conormal nonfocusing condition. Now
Theorem~\ref{15.5.2000.7} and Partial Theorem~\ref{15.5.2000.8} shows
that for, $2d(\bar m,\pa X)>t>d(\bar m,\pa X),$ the inclusion
\eqref{melwun1.375} holds. On the diffracted front, but away from the
direct front, application of Theorem~\ref{melwun1.367}
shows that $E_{\bar m}$ is conormal and of Sobolev regularity
$-n/2+1+(n-1)/2-\delta =1/2-\delta$ for every $\delta >0.$ Iterated
regularity with respect to this space then follows by interpolation.
\end{proof}

We now use this special case to prove Theorem~\ref{15.5.2000.8}.

\begin{proof}[Proof of Theorem~\ref{15.5.2000.8}] We begin by sharpening
the regularity results of Theorem~\ref{melwun1.374} to include
uniformity and regularity in the location of the pole.

Consider the Schwartz kernel of the fundamental solution, $E(t,x,y;x',y').$
For fixed $(x',y')\in(0,\epsilon )\times\pa X$ for $\epsilon >0$ small
enough and small $t$ and $x$ this is the distribution discussed in
Theorem~\ref{melwun1.374}. Moreover, by uniqueness of the solution to the
Cauchy problem it depends continuously on $(x',y').$ It follows that the
results described there hold uniformly in $(x',y').$ In particular near some
fixed point $\bar t=\bar x+\bar x'$ with $y,y'$ such that it is not on the
direct front,  
\begin{equation*}
((t-x')\pa_t+x\pa_x)^kD_y^\alpha E\in L^2_{\loc},\ \forall\ k,\alpha .
\label{melwun1.376}\end{equation*}
Here we have given up approximately half a derivative in the $x,y$
variables, and continuity in $x',y'$ to settle for iterative regularity
with respect to $L^2$ in all variables. By the symmetry and $t$-translation
invariance of the problem $E(t,x,y;x',y')=E(-t,x',y';x,y)$ so it also
follows that 
\begin{equation*}
((t-x)\pa_t+x'\pa_{x'})^kD_y^\alpha E\in L^2_{\loc},\ \forall\ k,\alpha
\label{melwun1.377}\end{equation*}
in the same region. All these vector fields commute, so by interpolation we
deduce regularity with respect to all the vector fields
simultaneously. Now, together with the wave operator itself, in both sets
of variables, these symbols of these operators define the (two components
of) the conormal bundle to $\{t=x+x'\}.$ Thus we deduce that $E$ is, away
from the direct front and for small $t,x,x'>0,$ conormal with respect to
this hyperplane.

Although we have only shown iterative regularity in $L^2,$
iterative regularity in $H^{\frac12-\epsilon }_{\loc}(\bbR\times X\times
X)$ in this set follows by interpolation.
Hence, a fortiori, $E$ is conormal with respect to the diffracted front
with iterate regularity in $H^{\frac12-\epsilon }(\bbR\times X\times X).$

Now, consider an admissible solution to the wave equation satisfying the
hypotheses of Theorem~\ref{15.5.2000.8}. Using Theorem~\ref{15.5.2000.7}
and a partition of unity, we may assume that the Cauchy data $(u(0), D_t
u(0))$ is identically zero in a neighborhood of (the projection to $X$ of)
all points geometrically related to $p \in R_{\pm, O}^\ep(\bar t)$ and that
at non geometrically related points it is supported in a microlocal
neighborhood of $R_{\pm, I}(\bar t).$ The nonfocusing condition then
implies that there exists $N \in \NN$ such that locally near $x=\bar t,$
\begin{multline*}
u(t,x,y) =\\
\int \frac{\pa E(t,x,y,x',y')}{\pa t} u_0(x',y')
\frac{dx'}{x'}\, dy' +\int E(t,x,y,x',y') u_1(x',y') \frac{dx'}{x'}\, dy'
\label{melwun1.378}\end{multline*}
we obtain the desired boundedness, by regarding $E$ and $\pa E/\pa t$
locally near the diffracted front as the kernels of Fourier integral
operator of order $\ep$ resp.\ $1+\ep$ on $\RR_+$, smoothly parametrized
by $t$, with values in $\Psi^{-\infty}(Y).$
\end{proof}

The explicit construction of the fundamental solution in the product case
by Cheeger-Taylor \cite{Cheeger-Taylor2, Cheeger-Taylor1} can be used to
show that Theorem~\ref{15.5.2000.8} cannot be strengthened by omitting the
hypothesis of the nonfocusing condition.  In particular, in Example 4.1 of
\cite{Cheeger-Taylor2} the authors show that on $X=\RR_+ \times S^1_2$, the
cone over the circle of circumference $4\pi$ (with coordinate $\theta$),
the fundamental solution $E(t,x,\theta, x',\theta')$ has a jump
discontinuity across the diffracted wavefront, where the value of the jump
varies smoothly with $\theta, \theta'$, at least in
$\abs{\theta-\theta'}<\pi$.  Now let $\phi(\theta)$ be a smooth function
supported in $\abs{\theta}<\delta$. The function
$$
v(t,x,\theta)=\int_{-2\pi}^{2\pi} E(t,x,\theta,x',\theta')\phi(\theta')\,
d\theta'
$$
is also an admissible solution to the wave equation, lying in $\ccal(\RR;
\dcal_{1/2-\ep})$ for all $\ep>0$.  For $t>x'$, it has geometrically
propagated singularities contained in $(x=t-x', \theta \in (-\delta,
\delta) \pm \pi)$.  On the other hand, for all $\theta' \in \supp \phi$,
$E(t,x,\theta,x',\theta')$ has a jump discontinuity along $x=t-x'$,
$\abs{\theta}<\pi-\delta$, hence $v$ also has a jump discontinuity this
surface, \ie\ this diffracted singularity is in $H^{1/2-\ep}$ for all
$\ep>0$.  Hence the diffracted singularity for $v$ is no smoother than the
geometrically propagated singularities.

\bibliography{all,My}
\bibliographystyle{amsplain}

\end{document}